\documentclass[leqno]{article}
\usepackage{amssymb}
\usepackage{amsmath, enumerate, vmargin}

\makeatletter
\renewcommand\section{\@startsection {section}{1}{\z@}%
 {-3.5ex \@plus -1ex \@minus -.2ex}%
 {2.3ex \@plus.2ex}%
 {\center \normalfont\large\bfseries}}
\makeatother

\setmarginsrb{3cm}{3cm}{3cm}{3cm}{1mm}{6mm}{0mm}{10mm}

\newtheorem{thm}{Theorem}[section]
\newtheorem{prop}[thm]{Proposition}
\newtheorem{cor}[thm]{Corollary}
\newtheorem{lem}[thm]{Lemma}
\newtheorem{defi}[thm]{Definition}
\newtheorem{remark}[thm]{Remark}
\newtheorem{example}[thm]{Example}
\newtheorem{pb}[thm]{Problem}

\newenvironment{rk}{\begin{remark}\rm}{\end{remark}}

\newcommand{\real}{{\mathbb R}}
\newcommand{\nat}{{\mathbb N}}
\newcommand{\ent}{{\mathbb Z}}
\newcommand{\comp}{{\mathbb C}}
\newcommand{\un}{1\mkern -4mu{\textrm l}}
\newcommand{\FF}{{\mathbb F}}
\newcommand{\T}{{\mathbb T}}

\newcommand{\E}{{\mathcal E}}
\newcommand{\F}{{\mathcal F}}

\newcommand{\M}{{\mathcal M}}
\newcommand{\N}{{\mathcal N}}
\renewcommand{\P}{{\mathcal P}}
\newcommand{\R}{{\mathcal R}}
\renewcommand{\S}{{\mathcal S}}

\renewcommand{\a}{\alpha}
\renewcommand{\b}{\beta}
\newcommand{\g}{\gamma}
\renewcommand{\d}{\delta}
\newcommand{\D}{\Delta}
\newcommand{\e}{\varepsilon}
\renewcommand{\th}{\theta}
\renewcommand{\l}{\lambda}

\renewcommand{\e}{\varepsilon}
\newcommand{\f}{\varphi}
\renewcommand{\O}{{\Omega}}
\renewcommand{\o}{{\omega}}
\newcommand{\s}{\sigma}

\newcommand{\tr}{\mbox{\rm tr}}
\newcommand{\ot}{\otimes}

\renewcommand{\t}{\tau}
\newcommand{\8}{\infty}
\newcommand{\el}{\ell}

\newcommand{\la}{\langle}
\newcommand{\ra}{\rangle}

\newcommand{\wh}{\widehat}
\newcommand{\n}{\noindent}
\newcommand{\pf}{\noindent{\it Proof.~~}}
\newcommand{\cqd}{\hfill$\Box$}
\newcommand{\be}{\begin{eqnarray*}}
\newcommand{\ee}{\end{eqnarray*}}
\newcommand{\beq}{\begin{equation}}
\newcommand{\eeq}{\end{equation}}
\newcommand{\dis}{\displaystyle}

\numberwithin{equation}{section}

\begin{document}


\title{Noncommutative maximal ergodic theorems}
\author{Marius Junge$^1$  and Quanhua Xu}

\date{}

\maketitle


\begin{abstract}
This paper is devoted to the study of various maximal ergodic
theorems in noncommutative $L_p$-spaces. In particular, we prove
the noncommutative analogue of the classical Dunford-Schwartz
maximal ergodic inequality for positive contractions on $L_p$ and
the analogue of  Stein's maximal inequality for symmetric positive
contractions. We also obtain the corresponding individual ergodic
theorems. We apply these results to a family of natural examples
which frequently appear in von Neumann algebra theory  and in
quantum probability.
\end{abstract}

 \setcounter{page}{0}
 \setcounter{section}{-1}
 \thispagestyle{empty}


\vskip 3cm

\n{\bf Plan:}
\begin{itemize}
\item[0.]  Introduction
\item[1.]  Preliminaries
\item[2.]  The spaces $L_p(\M;\el_\8)$
\item[3.]  An interpolation theorem
\item[4.]  Maximal ergodic inequalities
\item[5.]  Maximal inequalities for symmetric
contractions
\item[6.]  Individual ergodic theorems
\item[7.] The non tracial case
\item[8.] Examples
\end{itemize}

\footnotetext[1]{
   \hskip 0.1cm The first  author was partially supported by
the National Science Foundation Foundation DMS-0301116
 }

\newpage

\section{Introduction}\label{intro}

The connection between ergodic theory and theory of von Neumann
algebras goes back to the very beginning of the theory of ``rings
of  operators''. Maximal inequalities in ergodic theory provide an
important tool in classical analysis. In this paper we prove the
noncommutative analogue of the classical Dunford-Schwartz maximal
ergodic theorem thereby connecting these different aspects of
ergodic theory.

At the early stage of  noncommutative ergodic theory, only mean
ergodic theorems have been obtained (cf. e.g. \cite{ja1, ja2} for
more information). The study of individual ergodic theorems really
took off with Lance's pioneer work \cite {la-erg}. Lance  proved
that the ergodic averages associated with an automorphism of a
$\sigma$-finite von Neumann algebra which leaves invariant a
normal faithful state converge almost uniformly. Lance's ergodic
theorem was extensively extended and improved by (among others)
Conze, Dang-Ngoc \cite{cdn}, K{\"u}mmerer \cite{kum} (see
\cite{ja1, ja2} for more references). On the other hand, Yeadon
\cite{ye1} obtained a maximal ergodic theorem in the preduals of
semifinite von Neumann algebras. Yeadon's theorem provides a
maximal ergodic inequality which might be understood as a weak
type $(1,1)$ inequality. This inequality is the ergodic analogue
of Cuculescu's \cite{cu} result obtained previously for
noncommutative martingales. We should point out that in contrast
with the classical theory, the noncommutative nature of  these
weak type $(1, 1)$ inequalities seems a priori unsuitable for
classical interpolation arguments.

Since then the problem of finding a noncommutative analogue of the
Dunford-Schwartz maximal ergodic inequalities was left open. The
main reason is that all usual techniques in classical ergodic
theory involving maximal function seem  no longer available in
the noncommutative case. In fact, this applies for the definition
of the maximal function itself. As an example, we consider
 $$a_1=\left( \begin{array}{ll}
   2 & 0 \\
   0 & 0
 \end{array}\right) ,\quad
 a_2=\left( \begin{array}{ll}
   1 & 1 \\
   1 & 1
   \end{array}\right)
 ,\quad
 a_3= \left( \begin{array}{ll}
   0 & 0 \\
   0 & 1
      \end{array}\right).$$
Then there is no $2\times 2$ matrix $a$ such that
 $$\la\xi, a\xi\ra=\max\big\{\la\xi, a_1\xi\ra,\; \la\xi,
 a_2\xi\ra\; \la\xi, a_3\xi\ra\big\}  $$
holds for all $\xi\in\el_2^2$.

However, this obstacle has been overcome recently in the theory of
noncommutative martingale inequalities.  In fact, most of the
classical martingale inequalities have been successfully
transferred to the noncommutative setting. These include
Burkholder inequalities on conditioned square function
\cite{jx-burk}, Burkholder-Gundy inequalities on square function
\cite{px-BG}, Doob maximal inequality \cite{ju-doob}, Rosenthal
inequalities on independent random variables \cite{jx-ros} and
boundedness of martingale transforms \cite{ran-mtrans}. See the
survey \cite{xu-martsurv} for the state of the art regarding this
theory.

Let us point out that this new development of noncommutative
martingale inequalities is  inspired and motivated by interactions
with operator space theory. For instance, the formulation of the
noncommutative Doob maximal inequality was directly derived from
Pisier's theory of vector-valued noncommutative $L_p$-spaces
\cite{pis-ast}.

Following the well-known analogy between martingale theory and
ergodic theory, we show  that the techniques developed for
noncommutative martingales can be used to prove the
noncommutative maximal ergodic inequalities as well.

\smallskip

To state our main results we need some notation. Let $\M$ be a
semifinite von Neumann algebra equipped with a normal semifinite
faithful trace $\t$. Let $L_p(\M)$ be the associated
noncommutative $L_p$-space.  Let $T: \M\to\M$ be a linear map
which might satisfy some of the following properties
 \begin{itemize}
 \item[{\rm (0.I)}]  $T$ is a contraction on $\M$:
$\|Tx\|_\8\le \|x\|_\8$ for all $x\in \M$.
\item[{\rm (0.II)}]  $T$ is positive:
$Tx\ge 0$ if $x\ge 0$.
\item[{\rm (0.III)}] $\t\circ T\le\t$:
$\t(T(x))\le\t(x)$ for all $x\in L_1(\M)\cap\M_+$.
\item[{\rm (0.IV)}]  $T$ is symmetric relative to $\t$:
$\t(T(y)^*x)=\t(y^*T(x))$ for all $x,y\in L_2(\M)\cap\M$.
 \end{itemize}

The properties (0.I), (0.II)  and (0.III) will be essential for
what follows. If $T$ satisfies these properties, then $T$
naturally extends to a contraction on $L_p(\M)$ for all $1\le
p<\8$ (see Lemma \ref{extension} below). The extension will be
still denoted by $T$. If $T$ additionally has (0.IV), then its
extension is selfadjoint on $L_2(\M)$.  We consider the ergodic
averages of $T$:
 $$M_n(T)=\frac{1}{n+1}\,\sum_{k=0}^n T^k\,, \quad n\in\nat.$$
The following is one of our main results:

\begin{thm}\label{max+} Let $1<p<\8$ and $T$ be a linear map
satisfying {\rm (0.I)} - {\rm (0.III)} above.
 \begin{enumerate}[{\rm i)}]
 \item For any  $x\in L_p(\M)$ with $x\ge 0$ there is $a\in L_p(\M)$ such
that
 $$\forall\;n\in\nat,\ M_n(T)(x)\le a
 \quad\mbox{and}\quad\|a\|_p\le C_p\|x\|_p\,,$$
where $C_p$ is a positive constant depending only on $p$.
Moreover, $C_p\le C
 p^2(p-1)^{-2}$  and $(p-1)^{-2}$ is the optimal order of $C_p$ as $p\to 1$.
 \item If additionally $T$ satisfies
{\rm (0.IV)}, then for any  $x\in L_p(\M)$ with $x\ge 0$ there is
$a\in L_p(\M)$ such that
$$\forall\;n\in\nat,\ T^n(x)\le a
 \quad\mbox{and}\quad\|a\|_p\le C'_p\|x\|_p\,.$$
 \end{enumerate}
 \end{thm}

Part i) above is the noncommutative analogue of the classical
Dunford-Schwartz theorem in commutative $L_p$-spaces (cf.
\cite{dunford}). Note that the optimal order of the constant $C_p$
above is different from that in the commutative case, which is
$(p-1)^{-1}$ as $p\to 1$. Part ii) is the noncommutative analogue
of Stein's maximal ergodic inequality (see \cite{st-lp}). Note
that in the case where $\t$ is normalized (i.e. $\t(1)=1$), the
following weak form of part i) was obtained in \cite{gog}: Given
$\e>0$ such that $p-\e>1$ and $x\in L_p(\M)$  ($x\ge0$) there is
$a\in L_{p-\e}(\M)$ such that
 $$\forall\;n\in\nat,\ M_n(T)(x)\le a
 \quad\mbox{and}\quad\|a\|_{p-\e}\le C_{p,\e}\|x\|_p\,.$$

As in the commutative case, Theorem \ref{max+} also holds for all
elements of $L_p(\M)$ (not only the positive ones). This requires
an appropriate definition of the space $L_p(\M;\el_\8)$ in the
noncommutative setting (see section~2 for more details). On the
other hand, by discretization, we have a similar theorem for
semigroups.

The proof of Theorem \ref{max+}, i) relies on Yeadon's weak type
(1,1) maximal ergodic inequality already quoted before (see also
Lemma \ref{yeadon} below). As in the commutative case, the main
idea is  to interpolate this weak type (1,1) inequality with the
trivial case $p=\8$. Additional complications are due to the fact
that the weak type (1,1) estimate does not provide a majorant
$-a\le M_n(T)(x)\le a$ such that $a$ is in weak $L_1$. In our
proof of the noncommutative version of the classical Marcinkiewicz
theorem (see Theorem \ref{interpolation} below) we first establish
an intermediate inequality using noncommutative Lorentz spaces.
Then we use the real interpolation method. We should emphasize
that contrary to the classical situation, this interpolation
theorem is not valid for the spaces $L_p(\M;\el_\8)$ themselves
but only for their positive cones.

For the proof of part ii) of Theorem \ref{max+}, we adapt Stein's
arguments in \cite{st-lp} to the  noncommutative setting.

\medskip

As usual, the maximal ergodic inequalities in Theorem \ref{max+}
imply the corresponding pointwise ergodic theorems. The arguments
are standard in the tracial case. However, the non tracial case
requires additional work (see  section \ref{haagerup}). Our
approach to the individual ergodic theorems seems new. In order to
ensure pointwise convergence we use the space  $L_p(\M;c_0)$ which
is the closure of finite sequences in $L_p(\M;\ell_{\infty})$
($p<\8$). The main step towards individual ergodic theorems is
contained in the following result:

\begin{thm}\label{maxc00}
Let $1< p<\8$ and $T$  satisfy {\rm(0.I) - (0.III)}. Let $F$ be
the projection onto the fixed point subspace of $T$ considered as
a map on $L_p(\M)$. Then $\big(M_n(x)-F(x)\big)_n\in L_p(\M;c_0)$
for any $x\in L_p(\M)$.  If additionally, $T$ has {\rm (0.IV)},
then $\big(T^n(x)-F(x)\big)_n\in L_p(\M;c_0)$.
\end{thm}

Let us mention here one application to free group von Neumann
algebras for illustration. Let $\FF_n$ be a free group on $n$
generators, and let $VN(\FF_n)$ be its von Neumann algebra
equipped with the canonical normalized trace $\t$. Let $\l$ be the
left regular representation of $\FF_n$ on $\el_2(\FF_n)$. Recall
that $VN(\FF_n)$ is the von Neumann algebra on $\el_2(\FF_n)$
generated by $\{\l(g) :  g\in\FF_n\}.$ Let $| \cdot |$ denote the
length function on $\FF_n$ (relative to a fix family of $n$
generators). Haagerup \cite{haag-free} proved that the map $T_t$
defined by $T_t(\l(g))=e^{-t|g|}\l(g)$ extends to a completely
positive map on $VN(\FF_n)$. It is easy to check that $T_t$
possesses all properties (0.I) - (0.IV). Consequently, $T_t$
extends to a positive contraction on $L_p(VN(\FF_n))$ for all
$1\le p<\8$ (still denoted by $T_t$). It is then clear that
$(T_t)$ is a symmetric semigroup on $L_p(VN(\FF_n))$ and strongly
continuous for $p<\8$.  Thus applying Theorem \ref{max+}, ii) to
this semigroup, we obtain the following result formulated in terms
of (bilateral) almost uniform convergence:

\begin{thm}\label{free}
Let $1<p\le\8$. Then for any $x\in L_p(VN(\FF_n))$ with $x\ge0$
there is $a\in L_p(VN(\FF_n))$ such that
 $$\forall\;t>0,\ T_t(x)\le a
 \quad\mbox{and}\quad\|a\|_p\le C'_p\|x\|_p\,.$$
Consequently, $\lim_{t\to 0} T_t(x)=x$ bilaterally almost
uniformly $($and almost uniformly if  $p>2)$ for any $x\in
L_p(VN(\FF_n))$.
\end{thm}

The notion of (bilateral) almost uniform convergence is a
noncommutative analogue of the notion of almost everywhere
convergence. We refer to section~6 for the relevant definitions.
Note that when $n=1$, $VN(\FF_1)$ is just $L_\8(\T)$, where $\T$
is the unit circle and $T_t$ becomes the usual Poisson semigroup.
Thus the theorem above is the free analogue of the classical
radial maximal inequality and of the radial pointwise convergence
theorem about the Poisson integral in the unit disc.

\medskip

Let us end this introduction with a brief description of the
organization of the paper. The first six sections concern solely
the semifinite case. After a preliminary section, we give some
elementary properties on the vector-valued noncommutative
$L_p$-spaces $L_p(\M; \el_\8)$ in section \ref{Lpel}. These
vector-valued $L_p$-spaces were first introduced by Pisier
\cite{pis-ast} for injective von Neumann algebras and then
extended to general von Neumann algebras by the first named author
in \cite{ju-doob}. They provide  the main tool of this paper.

Section~3 is devoted to the noncommutative analogue of the
classical Marcinkiewicz interpolation theorem. This is the most
technical result of the paper. It seems reasonable to expect
further applications in the noncommutative setting.

Section~4 contains our first maximal ergodic theorems. The main
result there is   Theorem \ref{max+}, i). This is an immediate
consequence of the previous interpolation theorem.

Section~5 deals with the maximal inequalities when assuming the
symmetry condition (0.IV). In particular, we prove Theorem
\ref{max+}, ii). Our proof requires Stein's interpolation
technique using fractional averages (which makes it quite
involved).

In Section~6 we study the individual ergodic theorems. In
particular, we prove Theorem \ref{maxc00} above.

The objective of section~7 is to extend all previous results to
the general (non-tracial) von Neumann algebras by a reduction
argument. This argument is based on an important (unfortunately)
unpublished result due to Haagerup \cite{haag-red}. Let us mention
that the arguments for pointwise convergence in  Haagerup
$L_p$-spaces are  usually more delicate than their semifinite
counterparts. However, our new approach presented in section~6
permits us to give a unified treatment of both cases.

Section~8 presents some natural examples to which our theory
applies. These include the free products of completely positive
semigroups, the Poisson semigroup of a free group (which yields
Theorem \ref{free} above) and the $q$-Ornstein-Uhlenbeck
semigroups.

\medskip

The main results of this paper were announced in
\cite{jx-maxnote}.

\section{Preliminaries}\label{preli}

The noncommutative $L_p$-spaces used in this paper are most of the
time those based on semifinite von Neumann algebras, except those
in the last two sections. Thus in this preliminary section  we
concentrate ourselves only on the semifinite noncommutative
$L_p$-spaces.  There are numerous references for these spaces. Our
main reference is \cite{fk}.  The recent survey \cite{px-survey}
presents a rather complete picture on noncommutative integration
and contains a lot of references.

Let $\M$ be a semifinite von Neumann algebra equipped with a
normal semifinite faithful trace
$\t$. Let $\S_+$ denote the set of all $x\in \M_+$ such that
$\t({\rm supp}\,x)<\infty$, where ${\rm supp}\,x$ denotes the
support of $x$. Let $\S$  be the linear span of $\S_+$. Then $\S$
is a w*-dense $\ast$-subalgebra of $\M$. Given $0< p<\8$, we
define
 $$\|x\|_p=\big[\t(|x|^p)\big]^{1/p},\quad x\in\S,$$
where $|x|=(x^*x)^{1/2}$ is the modulus of $x$. Then
$(\S,\;\|\cdot\|_p)$ is a normed (or quasi-normed for $p<1$)
space, whose completion is the noncommutative $L_p$-space
associated with $(\M, \t)$, denoted by $L_p(\M,\t)$ or simply by
$L_p(\M)$. As usual, we set $L_\infty(\M,\t)=\M$ equipped with the
operator norm.

\medskip

The elements in $L_p(\M)$ can be viewed as closed densely defined
operators on $H$ ($H$ being the Hilbert space on which $\M$ acts).
We recall this briefly. Let $L_0(\M)=L_0(\M,\t)$ denote the space
of all closed densely defined operators on $H$ measurable with
respect to $(\M, \t)$.  For a measurable operator $x$ we define
its  generalized singular numbers by
 $$\mu_t(x)=\inf\big\{\l>0: \t\big(\un_{(\l,\8)}(|x|)\big)\le t\big\},
 \quad t>0.$$
Let
 $$V(\e,\d)=\{x\in L_0(\M): \mu_\e(x)\le\d\}.$$
Then $\{V(\e,\d): \e>0,\d>0\}$ is a system of neighbourhoods at
$0$ for which $L_0(\M)$ becomes a metrizable topological
$\ast$-algebra. The convergence with respect to this topology  is
called the convergence in measure. Moreover,  $\M$ is dense in
$L_0(\M)$.

The trace $\t$ extends to a positive tracial functional on the
positive part $L_0^+(\M)$ of $L_0(\M)$, still denoted by $\t$,
satisfying
 $$\t(x)=\int_0^\infty\mu_t(x)dt,\quad x\in L_0^+(\M).$$
Then for $0<p<\infty$,
 $$L_p(\M)=\big\{x\in L_0(\M): \t(|x|^p)<\infty\big\}$$
and for $x\in L_p(\M)$
 $$\|x\|_p^p=\t(|x|^p)=\int_0^\8(\mu_t(x))^p\,dt\,.$$
More generally, we can define the noncommutative Lorentz space
$L_{p,q}(\M)$:
 $$L_{p,q}(\M,\t)=\big\{x\in L_0(\M): \|x\|_{p,q}<\infty\big\},$$
where
 $$\|x\|_{p,q}=\Big(\int_0^\8\big(t^{\frac{1}{p}}\mu_t(x)
 \big)^q\,\frac{dt}{t}\Big)^{1/q}$$
for $q<\8$ and with the usual modification for $q=\8$. The
positive cone of $L_{p,q}(\M)$ is denoted by $L^+_{p,q}(\M)$

\medskip

As the commutative $L_p$-spaces, the noncommutative $L_p$-spaces
behave well with respect to the complex interpolation method and
the real interpolation method (in the semifinite case).  Let
$0<\theta<1$, $1\le p_0<p_1\le\8$ and $1\le q_0, q_1, q\le\8$.
Then
 \begin{equation}\label{complex}
 L_p(\M)=\big(L_{p_0}(\M),\;L_{p_1}(\M)\big)_{\theta}
 \quad (\mbox{with equal norms})
 \end{equation}
and
 \begin{equation}\label{real}
  L_{p,q}(\M)=\big(L_{p_0, q_0}(\M),\;L_{p_1, q_1}(\M)\big)_{\theta,q}
 \quad (\mbox{with equivalent norms}),
 \end{equation}
where
 $\displaystyle\frac{1}{p}=\frac{1-\theta}{p_0}+\frac{\theta}{p_1}$,
and where $(\cdot\,,\cdot)_{\theta}$, $(\cdot\,,\cdot)_{\theta,q}$
denote respectively the complex and real interpolation methods.
Our reference  for interpolation theory is \cite{bl}.

\medskip

Let $T$ be a linear map on $L_{p}(\M)$. $T$ is called {\it
positive} if $T$ preserves the positive cone of $L_{p}(\M)$, i.e.
$x\ge 0\Longrightarrow Tx\ge0$. The following lemma is elementary
and  certainly well-known. See \cite{ye1}, where the extensions of
$T$ to all $L_p(\M)$ were obtained but not their contractivity.
The contractive extensions are, of course, important in ergodic
theory. We include a proof for completeness.

\begin{lem}\label{extension}
Let $T$ satisfy {\rm (0.I)} - {\rm (0.III)}. Then $T$ extends in a
natural way to a positive contraction on $L_p(\M)$ for all $1\le
p<\8$. Moreover, $T$ is  normal on $\M$. If $T$ additionally has
{\rm (0.IV)}, the extension of $T$ on $L_2(\M)$ is selfadjoint.
\end{lem}

\pf It is clear that $T$ extends to $L_1^+(\M)$ and $\|Tx\|_1\le
\|x\|_1$ for all $x\in L_1^+(\M)$. Then by a standard argument,
$T$ extends to a bounded map on $L_1(\M)$, still denoted by $T$,
which is of norm $\le 2$ and positive too. By duality, $S=T^*:
\M\to\M$ is a bounded positive map. Consequently,
$\|S\|=\|S(1)\|_{\8}$ (see \cite{pa-cb}). However, one easily
checks that $\|S(1)\|_{\8}\le1$. Thus $S$ is contractive, and so
is $T$ on $L_1(\M)$.  Then by complex interpolation, $T$ extends
to a contraction on $L_p(\M)$ for all $1<p<\8$.

Note that the positive map $S$ on $\M$ introduced above satisfies
the same assumptions as $T$. Thus applying the result just proved
to $S$ instead of $T$, we see that $S$ can be extended to a
contraction on $L_1(\M)$. Then a simple calculation shows that
the adjoint of this extension of $S$ on $L_1(\M)$ is equal to
$T$. Hence $T$ is normal on $\M$. The last part is clear. \cqd

\medskip

In the sequel, unless explicitly specified otherwise, $T$ will
always denote a map on $\M$ with (0.I) - (0.III). The same symbol
$T$ will also stand for the extensions of $T$ on $L_p(\M)$ given
by Lemma \ref{extension}. Let $T$ be such a map. We form its
ergodic averages:
 $$M_n(T)=\frac{1}{n+1}\,\sum_{k=0}^n T^k\,.$$
$M_n(T)$ will be denoted by $M_n$ whenever no confusion can occur.

By general ergodic theory on Banach spaces (cf. \cite{dunford}),
one sees that $T$ is mean ergodic on $L_p(\M)$ for any $1<p<\8$,
i.e. $M_n(x)$ converges to $\wh x$ in $L_p(\M)$ for all $x\in\M$.
On the other hand, in case the trace $\tau$ is finite, by a
well-known mean ergodic theorem (cf. e.g. \cite[Theorem
2.2.1]{ja1} and the references therein), $M_n(x)$ converges to
$\wh x$ with respect to the strong operator topology for every
$x\in\M$. This implies (and in fact is equivalent to) that
$M_n(x)$ converges to $\wh x $ for any $x\in L_1(\M)$.

$T$ induces a canonical splitting on $L_p(\M)$ for $1< p<\8$:
 $$L_p(\M)=\F_p(T)\oplus \F_p(T)^\perp,$$
where $\F_p(T)=\{x\in L_p(\M):T(x)=x\}$ and $\F_p(T)^\perp$ is the
closure of the image $(I-T)(L_p(\M))$. The dual space of $\F_p(T)$
coincides with $\F_{p'}(T)$, where $p'$ is the index conjugate to
$p$. If in addition $\tau$ is finite, the previous splitting is
also true for $p=1$ and $p=\8$. Note then however that
$\F_\8(T)^\perp$ is the w*-closure of the image $(I-T)(L_\8(\M))$.

Let $F_p$ be the contractive positive projection from $L_p(\M)$
onto $\F_p(T)$. Then $F_2$ is the orthogonal projection from
$L_2(\M)$ onto $\F_2(T)$ and $F_p^*=F_{p'}$ for $1< p<\8$. (This
is also true for $p=1$ if $\t$ is finite.)  Note that $F_p$ and
$F_q$ coincide on $\F_p(T)\cap \F_q(T)$ for two different $p, q$.
This allows to denote the $F_p$'s by a same symbol $F$ in the
sequel. All previous facts are elementary and well known (cf. e.g.
\cite{ye1}).

\medskip

Let us transfer the discussion above to the setting of semigroups.
We will say that a semigroup $(T_t)_{t\ge 0}$ of linear maps on
$\M$ satisfies one of the conditions (0.I)-(0.IV) if so does $T_t$
for every $t\ge 0$. All semigroups considered in this paper will
be assumed to satisfy (0.I) - (0.III).  They will be further
required to be w*-continuous on $\M$ and such that $T_0$ is the
identity. As before, such semigroups are automatically extended to
positive contractive semigroups on $L_p(\M)$ for every $1\le
p\le\8$. Note that the w*-continuity of $(T_t)$ on $\M$ implies
that $(T_t)$ is strongly (i.e. norm) continuous on $L_p(\M)$ for
$1\le p<\8$. Put again
 $$M_t=\frac{1}{t}\int_{0}^t T^s\,ds,\quad t>0.$$
Note that for notational simplicity we will use the same letter
$M$ to denote the ergodic averages for a contraction as well as
for a semigroup. The precise meaning should be clear in the
concrete context. Again, the mean ergodic theorem asserts that
$M_t(x)$ converges to $F(x)$ in $L_p(\M)$  for all $x\in L_p(\M)$
($1<p<\8$) , where $F$ stands for the projection from $L_p(\M)$
onto the fixed point space of $(T_t)$, i.e. the space $\{x\in
L_p(\M)\;:\; T_t(x)=x, \forall\;t>0\}$.

\medskip

The following result due to Yeadon \cite{ye1} will play an
important role in this paper.  $\P(\M)$ denotes the lattice of
projetions in $\M$.  Given  $e\in\P(\M)$,  set $e^\perp=1-e$.

\begin{lem}\label{yeadon}
 Let $T$ satisfy {\rm (0.I) - (0.III)}.  Let $x\in L^+_1(\M)$.
Then for any $\l>0$  there is   $e\in\P(\M)$ such that
 $$ \sup_{n\ge 0}\|e\,M_n(T)(x)\,e\|_\8\le\l \quad\mbox{and}\quad
 \t(e^\perp)\le \frac{\|x\|_1}{\l}\ .$$
\end{lem}

The reader can easily recognize that this is a noncommutative
analogue of the classical weak type (1,1) maximal ergodic
inequality. Yeadon's theorem has a martingale predecessor obtained
by Cuculescu \cite {cu}.

\section{The spaces $L_p(\M;\el_\8)$}\label{Lpel}

A fundamental object of this paper is the noncommutative spaces
$L_p(\M;\el_\8)$. Given $1\le p\le\8$, $L_p(\M;\el_\8)$ is defined
as the space of all sequences $x=(x_n)_{n\ge 0}$ in $L_p(\M)$
which admit a factorization of the following form: there are $a,
b\in L_{2p}(\M)$ and $y=(y_n)\subset L_\8(\M)$ such that
 $$x_n=ay_nb,\quad\forall\; n\ge0.$$
We then define
 $$\|x\|_{L_p(\M;\el_\8)}=\inf\big\{\|a\|_{2p}\,
 \sup_{n\ge0}\|y_n\|_\8\,\|b\|_{2p}\big\} ,$$
where the infimum runs over all factorizations as above. One can
(rather easily) check that $\big(L_p(\M;\ell_\8),\;
\|\cdot\|_{L_p(\M;\el_\8)}\big)$ is a Banach space. These spaces
are introduced in \cite{pis-ast} and \cite{ju-doob}. (In
\cite{pis-ast}, $\M$ is required to be hyperfinite.) To gain a
very first understanding on $L_p(\M;\el_\8)$, let us consider a
positive sequence $x=(x_n)$. Then one can show that $x\in
L_p(\M;\el_\8)$ iff there are $a\in L^+_{p}(\M)$ and $y_n\in
L^+_\8(\M)$ such that
 $$x_n=a^{\frac{1}{2}}\, y_n\, a^{\frac{1}{2}},\quad\forall\; n\ge0.$$
We can clearly assume that the $y_n$ are positive contractions. Thus
if $x$ has a factorization as above, then $x_n\le a$ for all $n$.
Conversely, if $x_n\le a$ for some $a\in L^+_p(\M)$, then
$x_n^{1/2}=u_na^{1/2} $ for a contraction $u_n\in\M$, and so
$x_n=a^{1/2}u_n^*u_na^{1/2}$. Thus $x\in L_p(\M;\el_\8)$. In
summary, a positive sequence $x$ belongs to $L_p(\M;\el_\8)$ iff
there is $a\in L^+_p(\M)$ such that $x_n\le a$ for all $n$, and
moreover,
 $$\|x\|_{L_p(\M;\el_\8)}=\inf\big\{\|a\|_{p}\;:\; a\in L^+_p(\M)
 \;\mbox{s.t.}\; x_n\le a,\ \forall\;n\ge0\big\}.$$

\n{\bf Convention.} The norm of $x$ in $ L_p(\M;\el_\8)$ will be
very often denoted by
 $\big\|\sup_n^+x_n\big\|_p\ .$
\medskip

We should \underline{warn} the reader that
$\big\|\sup^+_nx_n\big\|_p$ is just a notation for $\sup_nx_n$
does not make any sense in the noncommutative setting. We find,
however, that $\big\|\sup^+_nx_n\big\|_p$ is more intuitive than
$\|x\|_{L_p(\M;\el_\8)}$.

\medskip

It is proved in \cite{ju-doob} that $L_p(\M;\el_\8)$ is a dual space
for every $p>1$. Its predual is $L_{p'}(\M;\el_1)$ ($p'$ being the
index conjugate to $p$). Let us define this latter space. Given
$1\le p\le\8$, a sequence $x=(x_n)$ belongs to $L_p(\M;\el_1)$ if
there are $u_{k n},\; v_{k n}\in L_{2p}(\M)$ such that
 $$x_n=\sum_{k\ge 0}u_{k n}^*\,v_{k n}$$
for all $n$ and
 $$\sum_{k, n\ge 0}u_{k n}^*\,u_{k n}\in L_p(\M),\quad
 \sum_{k, n\ge 0}v_{k n}^*\,v_{k n}\in L_p(\M).$$
Here all series are required to be convergent in $L_p(\M)$
(relative to the w*-topology in the case of $p=\8$).
$L_p(\M;\el_1)$ is a Banach space when equipped with the norm
 $$\|x\|_{L_p(\M;\el_1)}=\inf\Big\{
 \big\|\sum_{k, n\ge 0}u_{k n}^*\,u_{k n}\big\|_p^{\frac{1}{2}}\,
  \big\|\sum_{k, n\ge
  0}v_{k n}^*\,v_{k n}\big\|_p^{\frac{1}{2}}\Big\}\,,$$
where the infimum is taken over all $(u_{k n})$ and  $(v_{k n})$
as above. It is clear that finite sequences are dense in
$L_p(\M;\el_1)$ if $p<\8$. The duality between $L_p(\M;\el_\8)$
and
 $L_{p'}(\M;\el_1)$ is given by
  $$\langle x,\; y\rangle =\sum_{n\ge 0} \t(x_ny_n).$$

As previously for $L_p(\M;\el_\8)$, it is easy to describe the
positive sequences in $L_p(\M;\el_1)$. In fact, a positive sequence
$x=(x_n)$ belongs to $L_p(\M;\el_1)$ iff $\sum_nx_n\in L_p(\M)$.
If this is the case,
 $$\|x\|_{L_p(\M;\el_1)}=\big\|\sum_{n\ge 0}x_n\big\|_p\ .$$
Compare this equality (whose member on the right has the
\underline{usual} sense) with our previous convention for the norm
in $L_p(\M;\el_\8)$. This partly justifies the intuitive notation
$\big\|\sup^+_nx_n\big\|_p\,$.

\medskip

We collect some elementary properties of these spaces in the
following proposition. We denote by $L_p(\M;\el_\8^{n+1})$ the
subspace of $L_p(\M;\el_\8)$ consisting of all finite sequences
$(x_0, x_1,\,\cdots\,,x_n, 0,\,\cdots )$. In accordance with our
preceding convention, the norm of $x$ in $L_p(\M;\el_\8^{n+1})$
will be denoted by $\|\sup^+_{0\le k\le n}x_k\|_p$.  Similarly, we
introduce the subspace $L_p(\M;\el_1^{n+1})$ of $L_p(\M;\el_1)$.

\begin{prop}\label{vectorLp0} Let $1\le p\le\8$.
\begin{enumerate}[{\rm i)}]
\item  Each element in the unit ball of $L_p(\M;\el_\8)$
$($resp. $L_p(\M;\el_1))$ is a sum of sixteen $($resp. eight$)$
positive elements in the same ball.
\item A sequence $x=(x_n)$ in $L_p(\M)$ belongs to
$L_p(\M;\el_\8)$  iff
 $$\sup_{n\ge 0}\big\|\mathop{{\sup}^+}_{0\le k\le n}x_k\big\|_p<\8.$$
 If this is the case, then
 $$\big\|\,{\sup_n}^+x_n\big\|_p=\sup_{n\ge 0}
 \big\|\mathop{{\sup}^+}_{0\le k\le n}x_k\big\|_p\,.$$
\item Let $x=(x_n)$ be a positive sequence in
$L_p(\M;\el_\8)$. Then
 $$\big\|\,{\sup_n}^+x_n\big\|_p=\sup\big\{\sum_n\t(x_ny_n)\ :\ y_n\in
 L_{p'}^+(\M)\ \mbox{and}\ \big\|\sum_ny_n\big\|_{p'}\le 1\big\}.$$
\item We have the following Cauchy-Schwarz type
inequality: For any sequences $(x_n)$ and $(y_n)$  in $L_{2p}(\M)$
 $$\big\|\,{\sup_n}^+x_n^*\,y_n\big\|_p
 \le \big\|\,{\sup_n}^+x_n^*x_n\big\|_p^{\frac{1}{2}}\,
 \big\|\,{\sup_n}^+y_n^*y_n\big\|_p^{\frac{1}{2}}\ .$$
\end{enumerate}
\end{prop}

\pf  {\rm i)} First note that both $L_p(\M;\el_\8)$ and
$L_p(\M;\el_1)$ are closed with respect to involution. Thus we
need only to consider selfadjoint elements. Let $x$ be a
selfadjoint element (i.e. $x_n=x_n^*$ for all $n$) in the unit
ball of $L_p(\M;\el_\8)$. Write a factorization of $x$:
 $$x_n=a^*y_nb\quad\mbox{with}\quad \|a\|_{2p}\le 1,\
 \|b\|_{2p}\le 1\ \mbox{and}\ \sup_n\|y_n\|_\8\le 1.$$
Then by a standard polorization argument,
 \begin{eqnarray*}x_n
 &=&\frac{1}{4}\sum_{k=0}^3i^{-k}\,(a+i^kb)^*\,y_n\,(a+i^kb)
 =\frac{1}{4}\sum_{k=0}^3(a+i^kb)^*\,z_n\,(a+i^kb)\\
 &=& \sum_{k=0}^3\big(\frac{a+i^kb}{2}\big)^*\,z_n^+\,
 \big(\frac{a+i^kb}{2}\big) \; -\;
 \sum_{k=0}^3\big(\frac{a+i^kb}{2}\big)^*\,z_n^-\,
 \big(\frac{a+i^kb}{2}\big)\,,
 \end{eqnarray*}
where
 $$z_n= \frac{i^{-k}y_n+ (i^{-k}y_n)^*}{2}\,.$$
Hence the  assertion concerning $L_p(\M;\el_\8)$ follows. The one
for $L_p(\M;\el_1)$ is proved similarly.

{\rm ii)} It is trivial that
 $$\sup_{n\ge 0}\big\|\mathop{{\sup}^+}\limits_{0\le k\le
 n}x_k\big\|_p\le \big\|\,{\sup_n}^+x_n\big\|_p\,.$$
To prove the converse, we introduce the subspace $L_p(\M;\el^0_1)$
of $L_p(\M;\el_1)$, which consists of all finite sequences $x$
admitting a factorization as in the definition above of
$L_p(\M;\el_1)$ but with finite families $(u_{k n})$ and $(v_{k
n})$ only. Note that for $p<\8$, $L_p(\M;\el^0_1)$ is dense in
$L_p(\M;\el_1)$. Now let $x$ be a sequence such that $\sup_{n\ge
0}\big\|\sup^+_{0\le k\le n}x_k\big\|_p=1$. Define $\el:
L_{p'}(\M;\el^0_1)\to\comp$ by $\el(y)=\sum_n\t(x_ny_n)$. Then
$\el$ is a continuous linear functional of norm $\le1$. Thus if
$p>1$ (i.e. $p'<\8$), by the duality result in \cite{ju-doob}
already quoted previously, $\el$ can be identified with an element
of $L_p(\M;\el_\8)$. This element must be $x$, and so we are done
in this case. It remains to consider the case $p=1$. Using a
standard Hahn-Banach argument as presented in \cite{ju-doob}, we
deduce two states $\f$ and $\psi$ on $\M$ such that
 $$|\t(x_nu^*v)|\le \big(\f(u^*u)\big)^{\frac{1}{2}}\;
 \big(\psi(v^*v)\big)^{\frac{1}{2}}\ ,\quad n\ge0,\ \forall\; u, v\in
 \M.$$
Since $x_n\in L_1(\M)\simeq \M_*$ is a normal functional, we can
replace in the inequality above $\f$ and $\psi$ by their normal
parts respectively, and so we can assume $\f$ and $\psi$ already
normal. (In fact, in the present case, one can check that the
singular parts of $\f$ and
$\psi$ are zero.) Identifying $\f$ and $\psi$ with two positive
operators $a$ and $b$ in the unit ball of $L_1(\M)$, respectively,
we rewrite the
inequality above as
 $$|\t(x_nu^*v)|\le \|ua^{\frac{1}{2}}\|_2\;\|vb^{\frac{1}{2}}\|_2
 \ ,\quad \quad n\ge0,\  \forall\; u, v\in\M.$$
Then as in \cite{ju-doob}, we find contractions $y_n\in\M$ such
that $x_n=b^{1/2}\,y_n\, a^{1/2}$. Therefore,  $x\in
L_1(\M;\el_\8)$ and $\|\sup^+_nx_n\|_1\le 1$.

Note that if additionally $x$ is positive, in the Hahn-Banach
argument above, we can use only the positive cone
$L^+_{p'}(\M;\el^0_1)$ to get a factorization of $x$ as
$x_n=a^{1/2}\,y_n\, a^{1/2}$ with $a\in L_p^+(\M)$ and $y_n$
positive contractions. See \cite{ju-doob} for more details.

{\rm iii)} For $p>1$ this is already proved in \cite{ju-doob}. For
$p=1$ this is a consequence of {\rm ii)} and the previous remark.

{\rm iv)} We use duality. Let $(u_{k n})$ and $(v_{k n})$ be two
finite families in $L_{2p'}(\M)$. Then by the Cauchy-Schwarz inequality
 \begin{eqnarray*}
 \big|\sum_{k,n}\t(x_n^*\,y_n\,u_{k n}^*\,v_{k n})\big|
 &\le&
 \big(\t\sum_{k,n}x_n\,v_{k n}^*\,v_{k n}\,x_n^*\big)^{\frac{1}{2}}\;
 \big(\t\sum_{k,n}y_n\,u_{k n}^*\,u_{k n}\,y_n^*\big)^{\frac{1}{2}}\\
 &\le& \big\|\,{\sup_n}^+x_n^*x_n\big\|_p^{\frac{1}{2}}\,
 \big\|\,{\sup_n}^+y_n^*y_n\big\|_p^{\frac{1}{2}}\,
 \big\|\sum_{k,n}v_{k n}^*v_{k n}\big\|_{p'}^{\frac{1}{2}}
 \big\|\sum_{k,n}u_{k n}^*u_{k n}\big\|_{p'}^{\frac{1}{2}}\,.
 \end{eqnarray*}
whence the desired inequality.\cqd

\medskip

\begin{rk}\label{Lplinfty} i) From the proof of part ii) above, one sees
that the infimum defining the norm $\|\sup_nx_n\|_p$ is attained
for any $x\in L_p(\M;\el_\8)$ ($1\le p\le\8$). The same proof
shows that $L_1(\M;\el_\8)$ is identified as an isometric subspace
of the dual of $L_\8(\M;\el^0_1)$.

ii) We have a statement similar to Proposition \ref{vectorLp0},
ii)  for $L_p(\M; \el_1)$. On the other hand, let $L_p(\M;c_0)$ be
the closure of finite sequences in $L_p(\M; \el_\8)$ for $1\le
p<\8$. Then one can show that the dual space of $L_p(\M;c_0)$ is
equal to $L_{p'}(\M;\el_1)$ isometrically.

 iii) Neither $L_p(\M;\el_\8)$ nor $L_p(\M;\el_1)$ is
stable under the operation $(x_n)_n\mapsto (|x_n|)_n$. Thus
$\|\sup^+_nx_n\|_p\neq\|\sup^+_n|x_n|\,\|_p$ in general.
\end{rk}

\begin{rk}\label{vectorLp1bis}
If $\N\subset \M$ is a von Neumann subalgebra such that the trace
$\t$ restricted to $\N$ is semifinite on $\N$, then we have a
natural isometric inclusion $L_p(\N)\subset L_p(\M)$. This extends
to isometric inclusions:
 \be
 L_p(\N;\ell_{\8})\subset L_p(\M;\ell_{\8})\quad\mbox{and}\quad
 L_p(\N;\ell_1) \subset L_p(\M;\ell_1) \,.
 \ee
Indeed, by the definition of $L_p(\M;,\ell_{\8})$ and
$L_p(\M;\ell_1)$, the inclusions above are contractive. On the
other hand, the duality result from the preceding proposition and
remarks imply immediately that they are both  isometric.
\end{rk}

\medskip

 \begin{rk}\label{vectorLp0bis}
 The definition of $L_p(\M;\el_\8)$ and $L_p(\M;\el_1)$ can be extended to an
arbitrary index set. Let $I$ be an index set. Then
$L_p(\M;\el_\8(I))$ and $L_p(\M;\el_1(I))$ are defined similarly
as before. For instance, $L_p(\M;\el_\8(I))$ consists of all
families $(x_i)_{i\in I}$ in $L_p(\M)$ which can be factorized as
 $x_i=ay_ib$
with $a, b\in L_{2p}(\M)$ and a bounded family $(y_i)_{i\in
I}\subset L_\8(\M)$. The norm of $(x_i)_{i\in I}$ in $L_p(\M;
\el_\8)$ is defined as the infimum
 $$\inf \|a\|_{2p}\, \sup_i\|y_i\|_\8\,\|b\|_{2p}$$
running over all factorizations as above. As before, this norm is
also denoted by
 $$\big\|\,{\sup_i}^+ x_i\big\|_p\;.$$
Again the dual space of $L_p(\M; \el_1(I))$ for $p<\8$ is  $L_{p'}(\M; \el_\8(I))$.
Proposition \ref{vectorLp0} remains true in this general setting.
 \end{rk}

\medskip
 We end this section with a simple result on complex
interpolation of these vector-valued noncommutative $L_p$-spaces.

\begin{prop}\label{c-interp-vectorLp}
 Let $1\le p_0<p_1\le\8$ and $0<\th<1$. Then we have isometrically
 $$L_p(\M;\el_1)=\big(L_{p_0}(\M;\el_1),\;
 L_{p_1}(\M;\el_1)\big)_{\th}\,, \quad
L_p(\M;\el_\8)=\big(L_{p_0}(\M;\el_\8),\;
 L_{p_1}(\M;\el_\8)\big)_{\th}\; ,$$
where $\frac{1}{p}=\frac{1-\th}{p_0} + \frac{\th}{p_1}$.
\end{prop}

\pf We use the column and row spaces $L_p(\M; \el_2^c(\nat^2))$
and $L_p(\M; \el_2^r(\nat^2))$ (cf. \cite{px-BG} for the
definition). It is known that $\{L_p(\M; \el_2^c(\nat^2))\}_{1\le
p\le\8}$ form an interpolation scale with respect to the complex
interpolation. The same is true for the row spaces. Note that by
the definition of $L_p(\M;\el_1)$,  the bilinear map
 \begin{eqnarray*}
 B\;:\;
 &L_p(\M;\el_2^r(\nat^2))&\times ~~~ L_p(\M; \el_2^c(\nat^2))
 \longrightarrow
 L_p(\M;\el_1)\\
 &\hskip 0.3cm\big(u_{k n}\big)_{k,n\ge0} &\times\hskip 0.7cm
 \big(v_{k n}\big)_{k,n\ge0}~~
 \longmapsto \big(\sum_ku_{k n}v_{k n}\big)_{n\ge0}
 \end{eqnarray*}
is contractive (in fact, $L_p(\M;\el_1)$  is just the quotient
space of $L_p(\M; \el_2^r(\nat^2))\times L_p(\M; \el_2^c(\nat^2))$
by the kernel of $B$). Thus by the complex interpolation for
bilinear maps (cf. \cite{bl}), we deduce that
 $$B:L_p(\M; \el_2^r(\nat^2))\times L_p(\M; \el_2^c(\nat^2))\to
  \big(L_{p_0}(\M;\el_1),\; L_{p_1}(\M;\el_1)\big)_\th$$
is contractive. This yields
 \beq\label{c-interp-vectorLp1}
 L_p(\M;\el_1)\subset\big(L_{p_0}(\M;\el_1),\;
 L_{p_1}(\M;\el_1)\big)_{\th}\,,\quad \mbox{a contractive inclusion}.
 \eeq

Similarly, using the complex interpolation for trilinear maps, we
obtain the following contractive inclusion
 \beq\label{c-interp-vectorLp2}
 L_p(\M;\el_\8)\subset\big(L_{p_0}(\M;\el_\8),\;
 L_{p_1}(\M;\el_\8)\big)_{\th}\;.
 \eeq
Alternatively, this can be easily proved by using directly the
factorization of elements in $L_p(\M;\el_\8)$.

Dualizing the corresponding inclusion of
(\ref{c-interp-vectorLp1}) for finite sequences, we get
 \be
 \big(L_{p_0'}(\M;\el^n_\8),\;\big(L_{p_1}(\M;\el^n_1)\big)^*\big)^{\th}
 \subset L_{p'}(\M;\el^n_\8),
 \ee
where $(\cdot,\;\cdot)^\th$ denotes Calderon's second complex
interpolation method. However,
  \be
 \big(L_{p_0'}(\M;\el^n_\8),\;\big(L_{p_1}(\M;\el^n_1)\big)^*\big)_{\th}
 \subset
 \big(L_{p_0'}(\M;\el^n_\8),\;
 \big(L_{p_1}(\M;\el^n_1)\big)^*\big)^{\th}\;.
 \ee
It follows that
 \be
 \big(L_{p_0'}(\M;\el^n_\8),\;\big(L_{p_1}(\M;\el^n_1)\big)^*\big)_{\th}
 \subset L_{p'}(\M;\el^n_\8).
 \ee
Since $L_{p_0'}(\M;\el^n_\8)\cap\big(L_{p_1}(\M;\el^n_1)\big)^*$
is dense in the complex interpolation space on the left and
 $$L_{p_0'}(\M;\el^n_\8)\cap\big(L_{p_1}(\M;\el^n_1)\big)^* L_{p_0'}(\M;\el^n_\8)\cap L_{p_1'}(\M;\el^n_\8),$$
we deduce that
 \be
 \big(L_{p_0'}(\M;\el^n_\8),\; L_{p_1'}(\M;\el^n_\8)\big)_{\th}
 \subset L_{p'}(\M;\el^n_\8),\quad \mbox{a contractive inclusion}.
 \ee
Reformulating this for the indices $p_0,p_1$, we have
 \be
 \big(L_{p_0}(\M;\el^n_\8),\; L_{p_1}(\M;\el^n_\8)\big)_{\th}
 \subset L_{p}(\M;\el^n_\8),\quad \mbox{a contractive inclusion}.
 \ee
>From this we easily get the same inclusion for infinite sequence
spaces. Indeed, let $x=(x_k)_{k\ge0}$ be an element in
$\big(L_{p_0}(\M;\el_\8),\; L_{p_1}(\M;\el_\8)\big)_{\th}$ of norm
$\le 1$ and let $x^{(n)}= (x_0, x_1, ..., x_n, 0,0,...)$. Then
$x^{(n)}\in \big(L_{p_0}(\M;\el^{n+1}_\8),\;
L_{p_1}(\M;\el^{n+1}_\8)\big)_{\th}$ and is of norm $\le 1$. Thus
$x^{(n)}\in  L_{p}(\M;\el^{n+1}_\8)$ and is of norm $\le 1$.
Consequently,
 $$\sup_n\big\|x^{(n)}\big\|_{L_{p}(\M;\el_\8)}\le 1.$$
Therefore, by Proposition \ref{vectorLp0}, ii), we deduce that $x$
belongs to the unit ball of $L_{p}(\M;\el_\8)$, and so by
homogeneity, we obtain the converse inclusion of
(\ref{c-interp-vectorLp2}).

The converse inclusion of (\ref{c-interp-vectorLp1}) can be proved
similarly. But this time instead of Proposition \ref{vectorLp0},
it suffices to  use the fact that finite sequences are dense in
$\big(L_{p_0}(\M;\el_1),\; L_{p_1}(\M;\el_1)\big)_{\th}$. We omit
the details. \cqd

\medskip

\n{\bf Remark.} In a forthcoming paper we will show that the
interpolation equalities in Proposition \ref{c-interp-vectorLp}
are no longer true for the real interpolation. This is one of the
difficulties we will encounter for proving the Marcinkiewicz type
theorem in the next section.

\section{An interpolation theorem}\label{interp}

The main result of this section is a Marcinkiewicz type
interpolation theorem  for $L_p(\M;\,\el_\8)$. It is the key to
our proof of the noncommutative maximal ergodic inequalities. We
first introduce the following notion. For every integer $n\ge0$
assume given a map $S_n: L^+_p(\M)\to L^+_0(\M)$. Set
$S=(S_n)_{n\ge 0}$. Thus $S$ is a map which sends positive
operators to sequences of positive operators. We say that $S$ is
of {\it weak type} $(p,p)$ ($p<\8$) if there is a positive
constant $C$ such that for any $x\in L_p^+(\M)$ and any $\l>0$
there is a projection $e\in \M$ such that
 $$\tau(e^\perp)\le \big[C\,\frac{\|x\|_{p}}{\l}\big]^{p}
 \quad  \mbox{and}\quad  e\big(S_n(x)\big)e\le
 \l\,,\ \ \forall\;n\ge0.$$
Similarly, we say that $S$ is of {\it type} $(p,p)$ (this time $p$
may be equal to $\8$) if there is a positive constant $C$ such
that for any $x\in L_p^+(\M)$  there is $a\in L_p^+(\M)$
satisfying
 $$\|a\|_p\le C\,\|x\|_p\quad\mbox{and}\quad S_n(x)\le a,
 \ \ \forall\;n\ge0.$$
In other words, $S$ is of  type $(p,p)$ iff
 $$\|S(x)\|_{L_{p}(\M;\,\el_\8)}\le C\,\|x\|_{p}\,,\quad
\forall\;x\in L_{p}^+(\M).$$

\medskip

\begin{thm}\label{interpolation}
 Let $1\le p_0<p_1\le\8$. Let $S=\big(S_n\big)_{n\ge 0}$
be a sequence of maps from $L_{p_0}^+(\M)+ L_{p_1}^+(\M)$ into
$L_{0}^+(\M)$. Assume that $S$ is subadditive in the sense that
$S_n(x+y)\le S_n(x)+S_n(y)$ for all $n\in\nat$. If $S$ is of weak
type $(p_0, p_0)$ with constant $C_0$ and  of type $(p_1, p_1)$
with constant $C_1$, then for any $p_0<p< p_1$, $S$ is of type
$(p,p)$ with constant $C_p$ satisfying
\begin{equation}\label{interpolation estimate}
 C_p\le C\,C_0^{1-\th}\,C_1^{\th}\,\big(\frac{1}{p_0}-
 \frac{1}{p}\big)^{-2}\,,
\end{equation}
where $\th$ is determined by $1/p=(1-\th)/p_0+\th/p_1$ and $C$ is
a universal constant.
\end{thm}

The reader can easily recognize that this result is a
noncommutative analogue of the classical Marcinkiewicz
interpolation theorem. Recall that in the classical case the
constant $C_p$ is majorized by
$C_0^{1-\th}\,C_1^{\th}\,\big(1/p_0- 1/p\big)^{-1}$, i.e. without
the square in (\ref{interpolation estimate}). We will see later
that the estimate given by (\ref{interpolation estimate}) is
optimal in the noncommutative setting. This difference indicates
that though similar in form to the classical Marcinkiewicz
interpolation theorem, Theorem \ref{interpolation} cannot be
proved by the standard argument in the commutative case. The rest
of this section is entirely devoted to its proof. In the
following $S$ will be fixed as in the theorem above; $p$ will
denote a number such that $p_0<p<p_1$, and $\th$ is determined by
$1/p=(1-\th)/p_0+\th/p_1$. $C$ will stand for a universal
constant.

\medskip

The following lemma is entirely elementary.

\begin{lem}\label{vector matrix}
 Let $(x_{ij})$ be a finite matrix of bounded operators on a
Hilbert space $H$. Let $(e_i)$ and $(f_i)$ be two sequences of
pairwise disjoint projections in $B(H)$. Then
 $$\big\|\sum_{i,j}e_i x_{ij} f_j\big\|_{B(H)}\le
 \big\|\big(\|e_i x_{ij} f_j\|_{B(H)}\big)_{i,j}\big\|_{B(\el_2)}\ .$$
\end{lem}

\pf Let $\xi,\eta\in H$. Then
 \begin{eqnarray*}
 \langle\xi,\;\sum_{i,j}e_i x_{ij}f_j\eta\rangle
 &=&\sum_{i,j}\|e_i x_{ij} f_j\|\,\|e_i\xi\|\,\|f_j\eta\|\\
 &\le&\big\|\big(\|e_i x_{ij} f_j\|_{B(H)}\big)_{i,j}\big\|_{B(\el_2)}
 \big(\sum_i\|e_i\xi\|^2\big)^{\frac{1}{2}}\;\big(\sum_j\|f_j\eta\|^2
 \big)^{\frac{1}{2}}\\
 &\le& \big\|\big(\|e_i\,x_{ij}\,f_j\|_{B(H)}\big)_{i,j}\big\|_{B(\el_2)}
 \,\|\xi\|\,\|\eta\|.
 \end{eqnarray*}
This yields the desired inequality. \cqd

\begin{lem}\label{interpolation1}
 Let $x\in L_{p_0}^+(\M)$ and $\l>0$. Then there is  $e\in\P(\M)$ such that
$$\tau(e^\perp)\le \big[C_0 \l^{-1}\|x\|_{p_0}\big]^{p_0},\quad
 \big\|\,{\sup_n}^+\big(eS_n(x)e\big)\big\|_p\le
 C(1-\frac{p_0}{p})^{-1-\frac{1}{p}}\big[C_0\|x\|_{p_0}\big]^{\frac{p_0}
 {p}}\l^{1-\frac{p_0}{p}}\,.$$
\end{lem}

\pf Fix an $x\in L_{p_0}^+(\M)$. Set $x_n=S_n(x)$ for $n\in \nat$.
Since $S$ is of weak type $(p_0,p_0)$, given $k\in\ent$ there is
$f_k\in\P(\M)$ such that
 $$\tau(f_k^\perp)\le \big[C_0 2^{-k}\|x\|_{p_0}\big]^{p_0}
 \quad\mbox{and}\quad  f_kx_nf_k\le2^k\,,\quad\forall\;n\in\nat.$$
Let
 $$g_k=\bigvee_{j\ge k}f_j^\perp\,.$$
Then $(g_k)_{k\in\mathbb Z}$ is a decreasing sequence of
projections and
 $$\tau(g_k)\le \big[C_0 2^{-k+1}\|x\|_{p_0}\big]^{p_0}\,.$$
Thus $\lim_{k\to +\8}g_k=0$. Put $g_{-\8}=\lim_{k\to -\8}g_k$.
Then $g_{-\8}\ge g_{k}\ge f_k^\perp$, and so $g_{-\8}^\perp\le
g_{k}^\perp\le f_k$ for all $k\in\ent$. Put
 $$d_k=g_{k}-g_{k+1}\quad\mbox{and}\quad e_k=\sum_{j\le k}d_j\,.$$
Then $e_k=g_{-\8}-g_{k+1}$. We claim that
 $$(g_{-\8}^\perp+e_k)\,x_n\,(g_{-\8}^\perp+e_k)=e_kx_ne_k.$$
Since $g_{-\8}^\perp\le f_k$, by the choice of $f_k$
 $$g_{-\8}^\perp x_n g_{-\8}^\perp
 =g_{-\8}^\perp \big(f_kx_n f_k\big)g_{-\8}^\perp\le 2^k
 g_{-\8}^\perp\,,\quad k\in\ent;$$
thus letting $k\to-\8$, we get $g_{-\8}^\perp x_n
g_{-\8}^\perp=0$. On the other hand,
 $$\|g_{-\8}^\perp x_n e_k\|_{\8}\le
 \|g_{-\8}^\perp x_n g_{-\8}^\perp\|_{\8}^{\frac{1}{2}}\;\|e_k x_n
 e_k\|_{\8}^{\frac{1}{2}}=0\,;$$
whence
 $$g_{-\8}^\perp x_n e_k=0=e_k x_n g_{-\8}^\perp\,.$$
Therefore our claim is proved.

Now let $0<s<1$ such that $sp>p_0$. Set
 $$b_k=\sum_{j\le k} 2^{js}d_j\,.$$
Since the $d_j$ are disjoint projections, we have
 \be
 \|b_k\|_p&=&\big(\sum_{j\le k}2^{jsp}\t(d_j)\big)^{\frac{1}{p}}\\
 &\le& \big(\sum_{j\le k}2^{jsp}\big[C_0 2^{-j+1}\|x\|_{p_0}
 \big]^{p_0}\big)^{\frac{1}{p}}\\
 &\le& C(sp-p_0)^{-\frac{1}{p}}
 \big(C_0\|x\|_{p_0}\big)^{\frac{p_0}{p}}\,2^{k(s-\frac{p_0}{p})}\,.
 \ee
On the other hand, since the support of $b_k$ is equal to $e_k$,
$e_kb_k^{-\frac{1}{2}}$ can be regarded as a well-defined
operator,  and we have
 $$e_kb_k^{-\frac{1}{2}}=\sum_{j\le k} 2^{-\frac{js}{2}}d_j\,.$$
Thus
 $$b_k^{-\frac{1}{2}}e_kx_n e_kb_k^{-\frac{1}{2}}=\sum_{i,j\le k}
 2^{-\frac{is}{2}}2^{-\frac{js}{2}}d_ix_n d_j\,.$$
Since $d_i\le g_{i+1}^\perp\le f_{i+1}$, then by the choice of
$f_i$, we get
 $$\|d_ix_n d_j\|_\8\le \|d_ix_n d_i\|_\8^{\frac{1}{2}}\,
 \|d_jx_n d_j\|_\8^{\frac{1}{2}}\le 2^{\frac{i+1}{2}}\,2^{\frac{j+1}{2}}.$$
Therefore, using Lemma \ref{vector matrix},  we deduce
 \begin{equation}\label{matrix norm}
 \big\|b_k^{-\frac{1}{2}}e_kx_n e_kb_k^{-\frac{1}{2}}\big\|_\8\le
 C(1-s)^{-1}\,2^{k(1-s)}\,.
 \end{equation}
Note that the sequence $(e_kx_n e_k)_{n\ge 0}$ admits the
following factorization
 $$e_kx_n e_k=b_k^{\frac{1}{2}}\,\big[b_k^{-\frac{1}{2}}e_kx_n
 e_kb_k^{-\frac{1}{2}}\big]\,b_k^{\frac{1}{2}}\ .$$
Combining this with the previous inequalities, we obtain
 $$\big\|\,{\sup_n}^+\big(e_kx_n e_k\big)\big\|_p
 \le C (1-s)^{-1}(sp-p_0)^{-\frac{1}{p}}\big(C_0\|x\|_{p_0}
 \big)^{\frac{p_0}{p}}\,2^{k(1-\frac{p_0}{p})}\,.$$
Thus the choice of $s=(1+p_0)/(1+p)$ yields
 $$\big\|\,{\sup_n}^+\big(e_kx_n e_k\big)\big\|_p
 \le C (1-\frac{p_0}{p})^{-1-\frac{1}{p}}\big(C_0\|x\|_{p_0}
 \big)^{\frac{p_0}{p}}\,2^{k(1-\frac{p_0}{p})}\,.$$
Given $\l>0$ we choose $k$ such that $2^{k}\le\l<2^{k+1}$. Then
$e=g_{-\8}^\perp + e_k$ is the  desired projection.\cqd

\medskip

\n{\bf Remark.}  If we simply use the triangle inequality to
majorize the norm  $\big\|b_k^{-\frac{1}{2}}e_kx_n
e_kb_k^{-\frac{1}{2}}\big\|_\8$ instead of Lemma \ref{vector
matrix}, the estimates in (\ref{matrix norm}) becomes
$(1-s)^{-2}$. This does not give the right estimate in
(\ref{interpolation estimate}).

\medskip

The following lemma is a key step towards the proof of
Theorem \ref{interpolation}.

\begin{lem}\label{interpolation2}
 For any $x\in L^+_{p_0}(\M)\cap L^+_{p_1}(\M)$
 $$\big\|\,{\sup_n}^+ S_n(x)\big\|_p
 \le C(1-\frac{p_0}{p})^{-1-\frac{1}{p}}
 \big(C_0\|x\|_{p_0}\big)^{1-\th}\big(C_1\|x\|_{p_1}\big)^{\th}\,.$$
\end{lem}

\pf  Fix $x\in L^+_{p_0}(\M)\cap L^+_{p_1}(\M)$, and set
$x_n=S_n(x)$ as before. Let $\l>0$.  Choose $e\in\P(\M)$ as in
Lemma \ref{interpolation1}. Then
 $$x_n=ex_ne+e^\perp x_ne+ex_ne^\perp+ e^\perp x_ne^\perp\,.$$
Let us first estimate $\|\sup_n e^\perp x_ne^\perp\|_p\,.$ Since
$S$ is of type $(p_1,p_1)$, there is $a\in L^+_{p_1}(\M)$ such
that
 $$\|a\|_{p_1}\le C_1\|x\|_{p_1}\quad\mbox{and}
 \quad x_n\le a,\quad\forall\;n\in\nat.$$
Thus
 $$e^\perp x_ne^\perp \le e^\perp ae^\perp,\quad\forall\;n\in\nat.$$
With $r$ determined by $1/r=1/p-1/p_1$, by the H\"older
inequality, we have
 $$\|e^\perp ae^\perp\|_p\le\big(\tau(e^\perp)\big)^{\frac{1}{r}}\,
 \|a\|_{p_1}\le \big(C_0\|x\|_{p_0}\l^{-1}\big)^{\frac{p_0}{r}}\,
 C_1 \|x\|_{p_1}\,.$$
Therefore
 $$\big\|\,{\sup_n}^+\big(e^\perp x_ne^\perp\big)\big\|_p\le
 \big(C_0\|x\|_{p_0}\l^{-1}\big)^{\frac{p_0}{r}}\,C_1 \|x\|_{p_1}\,.$$
For the two mixed terms, by the Cauchy-Schwarz inequality in
Proposition \ref{vectorLp0}, we have
 $$\big\|\,{\sup_n}^+\big(e^\perp x_ne\big)\big\|_p\le
 \big\|\,{\sup_n}^+\big(e^\perp x_ne^\perp\big)\big\|_p^{\frac{1}{2}}\,
 \big\|\,{\sup_n}^+\big(ex_ne\big)\big\|_p^{\frac{1}{2}}\,,$$
and the same inequality holds for the other mixed term. Hence, we
deduce
\begin{eqnarray*}
 \big\|\,{\sup_n}^+x_n\big\|_p
 &\le& 2\big(\big\|\,{\sup_n}^+\big(e x_ne\big)\big\|_p
 +\big\|\,{\sup_n}^+\big(e^\perp x_ne^\perp\big)\big\|_p\big)\\
 &\le& C(1-\frac{p_0}{p})^{-1-{\frac{1}{p}}}\big(C_0\|x\|_{p_0}
 \big)^{\frac{p_0}{p}}\l^{1-\frac{p_0}{p}}+
 C\big(C_0\|x\|_{p_0}\l^{-1}\big)^{\frac{p_0}{r}}\,C_1\|x\|_{p_1}\ .
\end{eqnarray*}
Choosing $\l$ such that
 $$\l^{1-\frac{p_0}{p_1}}=\big(C_0\|x\|_{p_0}\big)^{-\frac{p_0}{p_1}}\,
 C_1\|x\|_{p_1}\,, $$
we obtain the desired inequality. \cqd

\medskip

The previous lemma can be restated as follows.

\begin{lem}\label{interpolation3}
 For any $x\in L^+_{p,1}(\M)$
 $$\big\|\,{\sup_n}^+ S_n(x)\big\|_p\le C(1-\frac{p_0}{p}
 )^{-1-\frac{1}{p}}\,C_0^{1-\th}\,C_1^{\th}\,\|x\|_{p,1}\,.$$
\end{lem}

\medskip

We will need to interpolate a compatible couple of cones. We refer
to \cite{bl} for the J- and K-methods in interpolation theory for
Banach spaces.  Let $(B_0,\; B_1)$ be a compatible couple of
Banach spaces. Let $A_i\subset B_i$ be a closed cone ($i=0,1$).
Given $0<\th<1$ and $1\le q\le\8$ we can define the J-method for
the couple $(A_0, A_1)$. More precisely, $(A_0, A_1)_{\th,q;J}$
consists of all $x\in B_0+B_1$ which admit a decomposition of the
following form
 $$x=\int_0^\8 u(t) \,\frac{dt}{t}\quad(\mbox{convergence in}\;B_0+B_1)$$
with $u(t)\in A_0\cap A_1$ such that
 $$\Big(\int_0^\8 \big[t^{-\th}\,\max\big(\|u(t)\|_{B_0},\,t\,
 \|u(t)\|_{B_1}\big)\big]^q\,\frac{dt}{t}\Big)^{\frac{1}{q}}<\8.$$
We define
 $$\|x\|_{(A_0, A_1)_{\th,q;J}}=\inf\Big\{
 \Big(\int_0^\8 \big[t^{-\th}\,\max\big(\|u(t)\|_{B_0},\,t\,
 \|u(t)\|_{B_1}\big)\big]^q\,\frac{dt}{t}\Big)^{\frac{1}{q}}\Big\},$$
where the infimum runs over all decompositions of $x$ as above.

\medskip

It is clear that
 $$(A_0, A_1)_{\th,q;J}\subset (B_0, B_1)_{\th,q;J}\;, \quad
 \mbox{a contractive inclusion}.$$
But in general the norm in $(A_0, A_1)_{\th,q;J}$ is not
equivalent to that of $(B_0, B_1)_{\th,q;J}$ when restricted to
$(A_0, A_1)_{\th,q;J}$. However, this is true for a couple of
noncommutative $L_p$-spaces.

\begin{remark}\label{interpolation of Lp cones}
 The following natural inclusion
 $$\big(L^+_{p_0,q_0}(\M),L^+_{p_1,q_1}(\M)\big)_{\th,q;J}\subset
 \big(L_{p_0, q_0}(\M), L_{p_1,q_1}(\M)\big)_{\th,q;J}$$
is isometric.
\end{remark}

\pf Let $x\in
\big(L^+_{p_0,q_0}(\M),L^+_{p_1,q_1}(\M)\big)_{\th,q;J}$. Let
 $$x=\int_0^\8 u(t) \,\frac{dt}{t}$$
be a decomposition of $x$ relative to $\big(L_{p_0,q_0}(\M),
L_{p_1,q_1}(\M)\big)_{\th,q;J}$ with $u(t)\in L_{p_0,q_0}(\M)\cap
L_{p_1,q_1}(\M)$ such that
 $$\Big(\int_0^\8 \big[t^{-\th}\,\max\big(\|u(t)\|_{p_0, q_0},\,t\,
 \|u(t)\|_{p_1, q_1}\big)\big]^q\,\frac{dt}{t}\Big)^{\frac{1}{q}}<\8.$$
Then we must find a similar decomposition of $x$ with all $u(t)$
in $L^+_{p_0, q_0}(\M)\cap L^+_{p_1,q_1}(\M)$ without increasing
the integral above. Since $x\ge0$, we can assume all $u(t)$ above
selfadjoint. Decomposing $u(t)$ into its positive and negative
part, we have
 $$x=\int_0^\8 u(t)^+ \,\frac{dt}{t} - \int_0^\8 u(t)^- \,\frac{dt}{t}
 \le \int_0^\8 u(t)^+ \,\frac{dt}{t}\, .$$
Therefore there is a contraction $v\in\M$ such that
 $$x^{\frac{1}{2}}=v\,
 \Big[\int_0^\8 u(t)^+ \,\frac{dt}{t}\Big]^{\frac{1}{2}}\,,$$
and so
$$x \int_0^\8 \big[v\,u(t)^+\,v^* \big]\,\frac{dt}{t}$$
yields the desired decomposition of $x$. \cqd

\medskip

We will need the following result from \cite{holm}, which gives
the optimal estimates for the equivalence constants in
(\ref{real}). Note that this result is stated in \cite{holm} for
the commutative $L_p$-spaces only. It is easy to see that the
noncommutative result follows immediately.

\begin{lem}\label{holm}
 Let $1\le p_0\not=p_1\le\8$ and $1\le q_0,q_1,q\le\8$. Then the
equivalence constants in the following equality
 $$L_{p,q}(\M)=\big(L_{p_0, q_0}(\M),L_{p_1,
 q_1}(\M)\big)_{\theta,q; K}$$
are estimated as follows
\begin{eqnarray*}
 &&C^{-1}\,\th^{-\min(\frac{1}{q},\,\frac{1}{q_0})}\,
 (1-\th)^{-\min(\frac{1}{q},\,\frac{1}{q_1})}\,\|x\|_{p,q}
 \le \|x\|_{\th, q;K}\\
 &&\hskip 4cm \le C\,\th^{-\max(\frac{1}{q},\,\frac{1}{q_0})}\,
 (1-\th)^{-\max(\frac{1}{q},\,\frac{1}{q_1})}\,\|x\|_{p,q}\ .
\end{eqnarray*}
\end{lem}

\begin{lem}\label{interpolation4}
The norm of the following inclusion
$$L_p(\M)\subset \big(L_{p_0,1}(\M),
 L_{p_1}(\M)\big)_{\th,p;J}$$
is majorized by
 $\displaystyle C\,(1-\th)^{1- 1/ p}$.
\end{lem}

\pf This is an immediate consequence of Lemma \ref{holm} by
duality.\cqd

\medskip
Our last result in this section  concerns  the real interpolation
of the positive cones $L^+_{p}(\M;\el_\8)$ of the spaces
$L_{p}(\M;\el_\8)$. Together with Lemma \ref{interpolation3}, it
constitutes the main technical part of the proof of Theorem
\ref{interpolation}.

\begin{lem}\label{interpolation5}
We have
 $$\big(L^+_{p_0}(\M;\el_\8),
 L^+_{p_1}(\M; \el_\8)\big)_{\th,p;J}\subset L^+_{p}(\M;\el_\8)$$
and the inclusion norm is
 $\displaystyle \le C\th^{-1+1/p}(1-\th)^{-1+1/p_1}$.
\end{lem}

\pf Let $x\in \big(L^+_{p_0}(\M;\el_\8), L^+_{p_1}(\M;
\el_\8)\big)_{\th,p;J}$ of norm $<1$. Choose $u(t)\in
L^+_{p_0}(\M;\el_\8)\cap L^+_{p_1}(\M;\el_\8)$ such that
 $$x=\int_0^\8 u(t) \,\frac{dt}{t}\quad\mbox{and}\quad
 \int_0^\8 \big[t^{-\th}\,J_{t}(u(t))\big]^p\,\frac{dt}{t}<1.$$
Here we have set
 $$J_{t}(y)=\max\big(\|y\|_{L^+_{p_0}(\M;\el_\8)},\,t\,
 \|y\|_{L^+_{p_1}(\M; \el_\8)}\big).$$
In order to prove $x\in L^+_{p}(\M;\el_\8)$, we use duality. Let
$y=(y_n)_n\in L^+_{p'}(\M;\el_1)$ of norm $\le1$. Set
$a=\sum_ny_n$. Then $\|a\|_{p'}\le1$. Let $K_t$ denote the
K-functional relative to $(L_{p_0'}(\M),\,L_{p_1'}(\M))$, i.e.
$K_t(\cdot)$ is the norm of the space
$L_{p_0'}(\M)+tL_{p_1'}(\M)$. Since $a\ge 0$, for every $t>0$
there is a spectral projection $e(t)$ of $a$ such that
 $$\|e(t)a\|_{p_0'}+t^{-1}\|e(t)^\perp a\|_{p_1'}\le
 2K_{t^{-1}}(a).$$
Then
 \begin{eqnarray*}
 \langle x, y\rangle
 &=&\sum_{n}\t(x_ny_n)=\int_0^\8\sum_n\t[u_n(t)y_n]\,\frac{dt}{t}\\
 &=&\int_0^\8\sum_n\t\big[u_n(t)[e(t)y_ne(t)+e(t)^\perp y_ne(t)^\perp
 +e(t) y_ne(t)^\perp+e(t)^\perp y_ne(t)]\big]\,\frac{dt}{t}.
\end{eqnarray*}
Since $y_n$ is positive, we have
 $$e(t) y_ne(t)^\perp
 +e(t)^\perp y_ne(t)
 \le e(t)y_ne(t)+e(t)^\perp y_ne(t)^\perp\ .$$
Hence  $u_n(t)\ge 0$ implies
 $$\t\big[u_n(t)[e(t) y_ne(t)^\perp+e(t)^\perp y_ne(t)]\big]
 \le \t\big[u_n(t)[e(t)y_ne(t)+e(t)^\perp y_ne(t)^\perp]\big].$$
Therefore
 \begin{eqnarray*}
 \langle x, y\rangle
 &\le& 2\int_0^\8\sum_n\t\big[u_n(t)[e(t)y_ne(t)
 +e(t)^\perp y_ne(t)^\perp]\big]\,\frac{dt}{t}\\
 &=&2\int_0^\8\big[\langle u(t), w(t)\rangle +
 \langle u(t), v(t)\rangle\big]\,\frac{dt}{t},
\end{eqnarray*}
where
 $$w(t)=\big(e(t)y_ne(t)\big)_{n\ge 0}\quad\mbox{and}\quad
 v(t)=\big(e(t)^\perp y_ne(t)^\perp\big)_{n\ge 0}\ .$$
Note that
 $$\|w(t)\|_{L_{p_0'}(\M;\el_1)}=\big\|\sum_ne(t)y_ne(t)\big\|_{p_0'}
 =\|e(t)a\|_{p_0'}\ .$$
Similarly,
 $$\|v(t)\|_{L_{p_1'}(\M;\el_1)}=\|e(t)^\perp a\|_{p_1'}\ .$$
It then follows that
 \begin{eqnarray*}
 \langle x, y\rangle
 &\le& 2\int_0^\8\big[\|u(t)\|_{L_{p_0}(\M;\el_\8)}\,
 \|w(t)\|_{L_{p_0'}(\M;\el_1)} +\|u(t)\|_{L_{p_1}(\M;\el_\8)}\,
 \|v(t)\|_{L_{p_1'}(\M;\el_1)}\big]\,\frac{dt}{t}\\
 &\le& 4\int_0^\8 J_{t}(u(t))\, K_{t^{-1}}(a)\,\frac{dt}{t}\\
 &\le& 4\Big(\int_0^\8 \big[t^{-\th}\,J_{t}(u(t))\big]^p\,\frac{dt}{t}
 \Big)^{\frac{1}{p}}
 \,\Big(\int_0^\8 \big[t^{\th}\,K_{t^{-1}}(a)\big]^{p'}
 \,\frac{dt}{t}\Big)^{\frac{1}{p'}}\\
 &\le& 4\, \big\|a\big\|_{ (L_{p_0'}(\M),\,L_{p_1'}(\M) )_{\th,
 p';K}}\ .
 \end{eqnarray*}
By Lemma \ref{holm}, the norm of the following inclusion
 $$L_{p'}(\M)\subset (L_{p_0'}(\M),\,L_{p_1'}(\M) )_{\th, p';K}$$
is controlled by $C\th^{-1/p'}(1-\th)^{-1/p'_1}$. Hence we deduce
 $$\langle x, y\rangle \le C\,\th^{-\frac{1}{p'}}(1-\th)^{-\frac{1}{p'_1}}\ .$$
Finally, taking the supremum over all positive $y$ in the unit
ball of $L_{p'}(\M;\el_1)$, we obtain the announced result.\cqd

\medskip

 Now we are in a position to prove Theorem \ref{interpolation}.

\medskip

\n{\it Proof of Theorem \ref{interpolation}.} Fix $x\in L^+_p(\M)$
such that $\|x\|_p\le 1$. Let $p_0<q<p$. Let $\eta$ and $\f$ be
determined by $1/q=(1-\eta)/p_0+\eta/p_1$ and $(1-\f)\eta+\f=\th$.
Applying Remark \ref{interpolation of Lp cones} and Lemma
\ref{interpolation4} with $q$ in place of $p_0$, we deduce that
 $x\in \big(L^+_{q,1}(\M),L^+_{p_1}(\M)\big)_{\f,p;J}$
and a decomposition of $x$
 $$x=\int_0^\8u(t)\,\frac{dt}{t}$$
such that
 $$\int_0^\8\big[t^{-\f}
  \max\big(\|u(t)\|_{q,1}\,,\ t\|u(t)\|_{p_1}\big)\big]^p\,\frac{dt}{t}
  \le C^p\,(1- \f)^{p-1}\,.$$
Set $v(t)=u(C_0^{\eta-1}\,C_1^{1-\eta}\,t)$. Then we again have
 $$x=\int_0^\8 v(t)\,\frac{dt}{t}\,.$$
Therefore, the subadditivity of $S$ implies
 \beq\label{pont}
 S(x)\le \int_0^\8 S(v(t))\,\frac{dt}{t}\
 {\dis\mathop=^{\rm def}}\ y\,.
 \eeq
Applying Lemma \ref{interpolation3} with $q$ instead of $p$ and by
the type $(p_1,p_1)$ of $S$, we deduce
 \be
 &&\max
 \big(\|S(v(t))\|_{L^+_q(\M;\el_\8)}\,,\
 t\|S(v(t))\|_{L^+_{p_1}(\M;\el_\8)}\big)\\
 &&~~~~\le C\,(1-\frac{p_0}{q})^{-1-\frac{1}{q}}
 \max\big(C_0^{1-\eta}\,C_1^{\eta}\,\|v(t)\|_{q,1},\
  t\,C_1\|v(t)\|_{p_1}\big)\,.
 \ee
Hence
 \be
 &&\int_0^\8\big[t^{-\f}
  \max\big(\|S(v(t))\|_{L^+_q(\M;\el_\8)}\,,\
 t\|S(v(t))\|_{L^+_{p_1}(\M;\el_\8)}\big)\big]^p\,\frac{dt}{t}\\
 &&~~\le\Big[C\,(1-\frac{p_0}{q})^{-1-\frac{1}{q}}\Big]^p
 \int_0^\8\big[t^{-\f}
  \max\big(C_0^{1-\eta}\,C_1^{\eta}\,\|v(t)\|_{q,1},\
  t\,C_1\|v(t)\|_{p_1}\big)\big]^p\,\frac{dt}{t}\\
  &&~~=\Big[C\,C_0^{1-\th}\,C_1^{\th}\,
  (1-\frac{p_0}{q})^{-1-\frac{1}{q}}\Big]^p
 \int_0^\8\big[t^{-\f}
  \max\big(\|u(t)\|_{q,1},\
  t\,\|u(t)\|_{p_1}\big)\big]^p\,\frac{dt}{t}\\
 &&~~\le \Big[C\,C_0^{1-\th}\,C_1^{\th}\,
  (1-\frac{p_0}{q})^{-1-\frac{1}{q}}
  (1- \f)^{1-\frac{1}{p}}\Big]^p\,.
 \ee
It thus follows that
 $$y\in \big(L^+_{q}(\M;\el_\8), L^+_{p_1}(\M;\el_\8)\big)_{\f,p;J}$$
and
 $$\|y\|_{\f,p;J}\le C\,C_0^{1-\th}\,C_1^{\th}\,
 (1-\frac{p_0}{q})^{-1-\frac{1}{q}}
 (1- \f)^{1-\frac{1}{p}}\,.$$
Therefore, by Lemma \ref{interpolation5} (applied with $q$ and
$\f$ in place of $p_0$ and $\th$, respectively), we  deduce that
$y\in L^+_{p}(\M;\el_\8)$ and
 \begin{eqnarray*}
 \|y\|_{L_{p}(\M;\el_\8)}
 &\le& C\,C_0^{1-\th}\,C_1^{\th}\,(1-\frac{p_0}{ q})^{-1-\frac{1}{q}}
 \f^{-1+\frac{1}{p}}(1-
 \f)^{-\frac{1}{p}+\frac{1}{p_1}}\\
 &\le& C\,C_0^{1-\th}\,C_1^{\th}\,(1-\frac{p_0}{
 q})^{-1-\frac{1}{q}}\,
 (\frac{1}{q}-\frac{1}{p})^{-1+\frac{1}{p}}\,.
 \end{eqnarray*}
Choosing $q$ such that
  $$\frac{1}{p_0}-\frac{1}{q}=\frac{1}{2}\,
  \big(\frac{1}{p_0}-\frac{1}{p}\big)\ ,$$
we get
 \be
  \|y\|_{L_{p}(\M;\el_\8)}\le C\,C_0^{1-\th}\,C_1^{\th}\,
  \big(\frac{1}{p_0} -\frac{1}{p}\big)^{-2}\,.
 \ee
This last inequality, together with (\ref{pont}), implies the
desired inequality (\ref{interpolation estimate}). Thus we have
completed the proof of Theorem \ref{interpolation}.\cqd


\section{Maximal ergodic inequalities}\label{maxergodic}

The following is our main maximal ergodic inequality in
noncommutative $L_p$-spaces. Restricted to the positive cone
$L_p^+(\M)$, it becomes Theorem \ref{max+}, i). Recall that for a
map $T$ with {\rm (0.I) - (0.III)}, $T$ also denotes its
extensions to $L_p(\M)$ given by Lemma \ref{extension}.

\begin{thm}\label{max}
Let $T$ be a linear map  with {\rm (0.I) - (0.III)}. Let
 $$M_n\equiv M_n(T)=\frac{1}{n+1}\,\sum_{k=0}^n T^k.$$
Then for every $1<p\le\8$ we have
 \begin{equation}\label{max1}
 \big\|\,{\sup_n}^+M_n(x)\big\|_p\le C_p\,\|x\|_p\ , \quad\forall\;x\in L_p(\M).
 \end{equation}
Moreover, $C_p\le C\, p^2(p-1)^{-2}$  and this is the optimal
order of $C_p$ as $p\to 1$.
\end{thm}

\pf Decomposing an operator into a linear combination of four
positive ones, we can assume $x\in L^+_p(\M)$. Now consider each
$M_n$ as a map defined on $L^+_1(\M)+L^+_\8(\M)$. Then $M_n$ is
positive and additive (and so subadditive too). Let $M=(M_n)_{n\ge
0}$. Yeadon's inequality says that $M$ is of weak type (1,1). On
the other hand, $M$ is trivially of type $(\8,\8)$ with constant
1. Therefore, by Theorem \ref{interpolation}, $M$ is of type
$(p,p)$ for every $1<p<\8$ with constant $C_p$ verifying
$$
 C_p\le C\,\big(1-\frac{1}{p}\big)^{-2}.
$$
Thus we have proved (\ref{max1}). The optimality of this estimate
follows from the optimal order of the best constant in the
noncommutative Doob maximal inequality obtained in \cite{jx-const}
and the following useful lemma due to Neveu, whose validity in the
noncommutative setting was observed by Dang-Ngoc \cite{dngoc}.

\begin{lem} Let $(\M_n)_{n\ge 0}$ be a decreasing sequence of von Neumann
subalgebras of $\M$. Assume that for every $n$ there is a normal
faithful conditional expectation $\E_n$ from $\M$ onto $\M_n$ such
that $\t\circ\E_n=\t$. Let $(\a_n)$ be an incresing sequence in
$[0, 1)$ with $\a_0=0$. The the map
 $$T=\sum_{n\ge 0} (\a_{n+1}-\a_n)\E_n$$
satisfies all conditions {\rm (0.I) - (0.IV)}. Moreover, given any
$\e>0$ one can choose $(\a_n)$ and an increasing sequence $(m_n)$
of positive integers such that
 $$\sum_{n\ge 0}\big\|M_{m_n}(T) - \E_n\big\|<\e,$$
where the norm is relative to  $L_p(\M)$ for any $1\le p\le\8$.
 \end{lem}

Note that if additionally $\lim_n\a_n=1$,   $T$ preserves the
trace $\t$ for the $\E_n$ preserve $\t$. With the help of this
lemma, one sees that the noncommutative maximal inequality
(\ref{max1}) implies the noncommutative Doob maximal inequality
proved in \cite{ju-doob} and $\d_p\le C_p$, where $\d_p$ is  the
best constant in the latter inequality. On the other hand, it was
shown in \cite{jx-const} that the optimal order of $\d_p$ is
$(p-1)^{-2}$ as $p\to 1$. It then follows that the estimate for
$C_p$ in (\ref{max1}) is optimal. Theorem \ref{max} is thus
proved. \cqd

\medskip
\n{\bf Remark.} The optimal order of the constant $C_p$ in
(\ref{max1}) implies that the estimate given in
(\ref{interpolation estimate}) is the best possible as $p\to p_0$
(with $p_0=1$). Recall that in the commutative case  the best
$C_p$ in (\ref{max1}) is of order $(p-1)^{-1}$ as $p\to 1$. This
explains partly the extra (noncommutative) effort in getting
(\ref{max1}).

\medskip

We will see in section~6 that Theorem \ref{max} implies the
ergodic averages $(M_n(x))_n$ converge bilaterally almost
uniformly for any $x\in L_p(\M)$.  However, for $p>2$ the
bilateral almost uniform convergence can be improved to the almost
uniform convergence. This improvement will be a consequence of the
following corollary of Theorem \ref{max}. For the  formulation of
this result we need a further notation from \cite{Mu,dj}. Let
$2\le p\le\8$ and $I$ be an index set. We define the space
$L_p(\M; \el_\8^c(I))$ as the family of all $(x_i)_{i\in I}\subset
L_p(\M)$ for which there are $a\in L_p(\M)$ and $(y_i)_{i\in
I}\subset L_\8(\M)$ such that
 $$x_i=y_ia \quad\mbox{and}\quad \sup_{i\in I}\|y_i\|_\8<\8.$$
$\|(x_i)\|_{L_p(\M; \el_\8^c(I))}$ is then defined to be the
infumum $\{\sup_{i\in I}\|y_i\|_\8\,\|a\|_p\}$ over all
factorizations of $(x_i)$ as above. It is easy to check that $\|\
\|_{L_p(\M; \el_\8^c(I))}$ is a norm, which makes $L_p(\M;
\el_\8^c(I))$ a Banach space. Note that $(x_i)\in L_p(\M;
\el_\8^c(I))$ iff $(x_i^*x_i)\in L_{\frac{p}{2}}(\M; \el_\8(I))$.
If $I=\nat$, $L_p(\M; \el_\8^c(I))$ is simply denoted by $L_p(\M;
\el_\8^c)$

\begin{cor}\label{maxc}
 Let $T$ be as in Theorem \ref{max} and $2<p\le\8$. Then
 $$\big\|\big(M_n(x)\big)_{n\ge0}\big\|_{L_p(\M; \el_\8^c)}
 \le \sqrt{C_{p/2}}\,\|x\|_p\;,
 \quad\forall\;x\in L_p(\M).$$
 \end{cor}

\pf Let $x\in L_p(\M)$. By decomposing $x$ into its real and
imaginary parts, we can assume $x$ selfadjoint. Since $T$ is
positive, so is $M_n$ for every $n$. Thus by the classical Kadison
inequality \cite{kad-cs}, we have
 $$ (M_n(x))^2\le M_n(x^2).$$
Thus applying Theorem \ref{max} to $x^2\in L_{p/2}(\M)$ we get
$b\in L^+_{p/2}(\M)$ such that
 $$\|b\|_{p/2}\le C_{p/2}\,\|x^2\|_{p/2}\quad\mbox{and}\quad
 M_n(x^2)\le b,\quad\forall\;n\ge 0.$$
Hence $(M_n(x))^2\le b$. It then follows that for each $n$ there
is a contraction $y_n\in \M$ such that $M_n(x)=y_nb^{1/2}$. This
gives the desired factorization of $\big(M_n(x)\big)_{n\ge0}$ as
an element in $L_p(\M; \el_\8^c)$ and thus proves the corollary.
\cqd

\medskip

 The following maximal inequality for multiple ergodic averages is
an easy consequence of Theorem \ref{max}.

\begin{cor}\label{max-multi}
 Let $T_1, ...\,,T_d$ be $d$ maps satisfying {\rm (0.I) - (0.III)}. Set
 $$M_{n_1,...\,,n_d}=\big[\prod_{j=1}^d\frac{1}{n_j+1}\big]
 \sum_{k_d=0}^{n_d}
 \;\cdots\;\sum_{k_1=0}^{n_1} T_d^{k_d}\,\cdots\,T_1^{k_1}\,.$$
 Then for any $1<p<\8$
 $$
 \big\|\mathop{{\sup}^+}_{n_1,...\,,n_d}\,M_{n_1,...\,,n_d}(x)
 \big\|_p \le C_p^d\;\|x\|_p\ ,
 \quad\forall\;x\in L_p(\M)
 $$
and for $2<p<\8$
 $$
 \big\|\big(M_{n_1,...\,,n_d}(x)\big)_{n_1,...\,,n_d}
 \big\|_{L_p(\M; \el_\8^c(\nat^d))}\le
 C_{p/2}^{\frac{d}{2}}\;\|x\|_p\ , \quad\forall\;x\in L_p(\M).
 $$
 \end{cor}

\pf The first part is obtained from Theorem \ref{max} by iteration.
The second is proved in the same way as Corollary \ref{maxc}. \cqd

\medskip

By a standard discretization argument, Theorem \ref{max} and the
previous corollaries imply the following maximal ergodic
inequalities for semigroups.

 \begin{thm}\label{max-sg} {\rm i)}
Let $(T_t)_{t\ge0}$ be a semigroup satisfying the conditions {\rm
(0.I) - (0.III)}. Let
 $$M_t=\frac{1}{t}\,\int_0^t T_s\,ds,\quad t>0.$$
Then for $1<p<\8$
 \be
 \big\|\,{\sup_t}^+M_t(x)\big\|_p\le C_p\,\|x\|_p\ ,\quad
 \forall\;x\in L_p(\M)
 \ee
and for $2<p<\8$
 $$
 \big\|\big(M_t(x)\big)_{t>0}\big\|_{L_p(\M; \el_\8^c(\real_+))}\le
 \sqrt{C_{p/2}}\;\|x\|_p\ , \quad\forall\;x\in L_p(\M).
 $$
 \indent{\rm ii)} Let  $(T_t^{(1)})_{t\ge 0}, ...,
 (T_t^{(d)})_{t\ge 0}$ be $d$
such semigroups. Let
 $$M_{t_1,...\,, t_d}=\frac{1}{t_1\,\cdots\, t_d}\,
 \int_0^{t_d}T^{(d)}_{s_d}\,ds_d
 \,\cdots \,\int_0^{t_1}T^{(1)}_{s_1}\,ds_1.$$
Then for any $1<p<\8$
 \be
 \big\|\mathop{{\sup}^+}_{t_1>0,...\,,t_d>0}\,M_{t_1,\,...\,,t_d}(x)\big\|_p\le
 C_p^d\;\|x\|_p\ , \quad\forall\;x\in L_p(\M)
 \ee
 and for $2<p<\8$
 $$
 \big\|\big(M_{t_1,...\,,t_d}(x)\big)_{t_1,...\,,t_d}
 \big\|_{L_p(\M; \el_\8^c(\real^d))}\le
 C_{p/2}^{\frac{d}{2}}\;\|x\|_p\ , \quad\forall\;x\in L_p(\M).
 $$
 \end{thm}

\pf We show only the first inequality in i). Recall that  the
semigroup $(T_t)_{t\ge0}$ is strongly continuous on $L_p(\M)$,
i.e. for any $x\in L_p(\M)$ the function $t\mapsto T_t(x)$ is
continuous from $[0,\8)$ to $L_p(\M)$, and so is the function
$t\mapsto M_t(x)$. Thus to prove the first inequality in i) it
suffices to consider $M_t(x)$ for $t$ in a dense subset of
$(0,\8)$, for instance, the subset $\{n2^{-m}\;:\; m, n\in
\nat\}$. Using once more the strong continuity of
$(T_t)_{t\ge0}$, we can replace the integral defining $M_t(x)$ by
a Riemann sum. Thus we have approximately
 \be
 M_{n2^{-m}}(x)
 &=&\frac{1}{n2^{-m}}\sum_{k=0}^{n-1}\int_{k2^{-m}}^{(k+1)2^{-m}}T_s(x)\,ds\\
 &\approx& \frac{1}{n}\sum_{k=0}^{n-1}T_{k2^{-m}}(x)
 =M_{n-1}(T_{2^{-m}})(x).
 \ee
Thus by Theorem \ref{max} applied to $T=T_{2^{-m}}$, we obtain
 $$\big\|\,{\sup_n}^+M_{n2^{-m}}(x)\big\|_p\le C_p\,\|x\|_p\;.$$
Since the subsets $\{n2^{-m}\;:\;  n\in \nat\}$ are increasing in
$m$, by Proposition \ref{vectorLp0}, we get
 $$\big\|\mathop{{\sup}^+}_{m, n} M_{n2^{-m}}(x)\big\|_p\le C_p\,\|x\|_p\;.$$
This is the desired inequality.\cqd

\medskip
 It is easy to show that the ergodic averages in Theorem \ref{max-sg}
can be replaced by many other averages. Let us consider, for
instance, the Poisson semigroup subordinated to $(T_t)$:
 \beq\label{subordination}
 P_t =\frac{1}{\sqrt\pi}\int_0^\infty \frac{e^{-u}}{\sqrt
 u}T_{t^2/4u} \,du.
 \eeq
Recall that if $A$ denotes the infinitesimal generator of $(T_t)$,
then that of $(P_t)$ is $-(-A)^{1/2}$. More generally, given any
$0<\a<1$, we can consider a semigroup $(P_t)$ subordinated to
$(T_t)$ via the following formula:
 \beq\label{subordinationbis}
 P_t = \int_0^\infty
 \varphi(s)T_{t^\b s}\,ds,
 \eeq
where $\b=1/\a$ and $\f$ is the function on $\real_+$ defined by
 \be
 \varphi(s)=\int_0^\infty \exp\big[st\cos\theta -
 t^\alpha\cos(\alpha\theta)\big]\times\sin\big[st\sin\theta -
 t^\alpha\cos(\alpha\theta)+\theta\big]dt,
 \ee
$\theta$ being any number in $[\pi/2,\ \pi]$. When $\a=1/2$,
(\ref{subordinationbis}) reduces to (\ref{subordination}). Note
that the infinitesimal generator of $(P_t)$ in
(\ref{subordinationbis}) is $-(-A)^{\a}$ (see \cite[IX]{yosida}).

\begin{cor}\label{max-poisson}
 Let $(T_t)$ be a semigroup verifying {\rm
(0.I) - (0.III)}. Let $0<\a<1$ and $(P_t)$ be the semigroup
subordinated to $(T_t)$ as in $(\ref{subordinationbis})$. Then for
$1<p\le \8$
 \beq\label{max-poisson1}
 \big\|\,{\sup_t}^+ P_t(x)\big\|_p\le C_{\a,p}\,\|x\|_p\ ,\quad
 \forall\;x\in L_p(\M)
 \eeq
and for $2<p<\8$
 $$
 \big\|\big(P_t(x)\big)_{t>0}\big\|_{L_p(\M; \el_\8^c(\real_+))}
 \le C_{\a,p}'\;\|x\|_p\ , \quad\forall\;x\in L_p(\M).
 $$
\end{cor}

\pf Let us first rewrite (\ref{subordinationbis}) as
 $$P_t = t^{-\b}\,\int_0^\infty\varphi(t^{-\b}s)T_{s}\,ds
 =t^{-\b}\,\int_0^\infty\f(t^{-\b}s)\,
 \frac{d}{ds}(sM_s)\,ds.$$
Thus by integration by parts,
 $$
 P_t= t^{-\b}\int_0^\8\f' (t^{-\b}s)\,t^{-\b}s\,M_s\,ds
 =\int_0^\8s\f' (s)\,\,M_{t^\b s}\,ds.
 $$
Therefore, by Theorem \ref{max-sg}, i), for any $x\in L_p(\M)$
 \be
 \big\|\,{\sup_t}^+P_t(x)\big\|
 \le C_p\,\int_0^\8 s|\f' (s)|\,ds\,\|x\|_p.
 \ee
It is easy to see that the integral on the right is finite by
virtue of the definition of $\f$. Thus we have proved
(\ref{max-poisson1}). In the same way, we get the second
inequality. \cqd

\begin{rk}\label{yeadon-sg}
Let $(T_t)$ be a semigroup as in Theorem \ref{max-sg}, i). Using
Lemma \ref{yeadon} and the preceding discretization argument, one
can easily obtained the following weak type (1,1) inequality: for
any $x\in L_1^+(\M)$ and $\l>0$ there is a projection $e\in\M$
such that
 $$\sup_{t>0}\big\|eM_t(x)e\big\|_\8\le\l\quad\mbox{and}\quad
 \t(e^\perp)\le \frac{\|x\|_1}{\l}\;.$$
Moreover, $M_t$ in the inequality above can be replaced by $P_t$
in (\ref{subordinationbis}).
\end{rk}


\section{Maximal inequalities for symmetric
contractions}\label{symergodic}

The  main result of this section is the following, which is a
reformulation of Theorem \ref{max+}, ii) for general elements in
$L_p(\M)$.

\begin{thm}\label{smax}
Let $T$ be a linear map on $\M$ satisfying {\rm (0.I) - (0.IV)}.
Then for any $1<p<\8$ we have
 \beq\label{smax1}
 \big\|\,{\sup_n}^+T^n(x)\big\|_p\le C'_p\,\|x\|_p\ ,\quad
 \forall\;x\in L_p(\M),
 \eeq
where $C'_p$ is a constant depending only on $p$.
\end{thm}

 This is the noncommutative analogue of a classical inequality due
to Stein (cf. \cite {st-erg}; see also \cite[Chapter III]{st-lp}).
Stein's approach is via complex interpolation. The main ingredient
is the maximal ergodic inequality (\ref{max1}), which allows to
deduce similar maximal inequalities for fractional averages. We
refer to \cite{sta} (and the references therein) for more general
maximal inequalities based on Rota's dilation theorem. (We are
grateful to the referee for bringing \cite{sta} to our attention.)
Let us point out that Rota's theorem is not sufficiently
understood in the noncommutative setting and hence Starr's method
is not yet available. It is Stein's original approach that suits
well to the noncommutative setting. Thus we will follow the same
pattern as in \cite[Chapter III]{st-lp}. Throughout the remainder
of this section $T$ will be fixed as in Theorem \ref{smax}.

\medskip

We begin by introducing the fractional averages on the powers of
$T$. Given a complex number $\a $ and a nonnegative integer $n$,
set
 $$
  A_n^\a = \frac{(\a+1)(\a +2)\cdots (\a+n)}{n!}
 $$
and
 $$S_n^\a\equiv S_n^\a(T)=\sum_{k=0}^{n}A_{n-k}^{\a-1}\,T^k\,,\quad
 M_n^\a\equiv M_n^\a(T)=(n+1)^{-\a}\,S_n^\a\ .$$
The $M_n^\a$ are the so-called fractional averages of the $T^k$.
Note that $M_n^0=T^n$ and $M_n^1$ is the usual ergodic average
$M_n$ already considered before. Also if $\a$ is a negative
integer $-m$, then
 $$S_n^{-m}=\D_n^m\big((T^k)_{k\ge0}\big),$$
where $\D_n$ denotes the first difference map on sequences, i.e.
 $$\D_n(a)=a_n-a_{n-1}$$
for every sequence $a=(a_n)$. Then $\D_n^m=\D_n( \D_n^{m-1})$ is
defined by induction and is the difference map of order $m$. Here
and in the sequel we adopt  the convention that for any sequence
$(a_n)_{n\ge 0}$ we put $a_n=0$ for $n<0$. Since we will only
consider actions of $\D_n^m$ on the sequence $(T^k)_{k\ge0}$, we
will simply put
 $$\D_n^m=\D_n^m\big((T^k)_{k\ge0}\big).$$
Thus $M_n^{-m}=(n+1)^m\D_n^m$.

\medskip

We will need a generalization of Theorem \ref{smax}.

\begin{thm}\label{smaxbis}
 Let $T$ be as in Theorem \ref{smax} with the additional
assumption that $T$ is positive as an operator on $L_2(\M)$
$($i.e. $\la x,\; Tx\ra=\t(x^*Tx)\ge0$ for all $x\in L_2(\M))$.
Then for all $\a\in\comp$ and $p\in(1,\8)$ we have
 \beq\label{smaxbis1}
 \big\|\,{\sup_n}^+ M_n^\a(x)\big\|_p\le
 C_{\a, p}\,\|x\|_p,\quad\forall\;x\in L_p(\M)\,,
 \eeq
where $C_{\a, p}$ is a constant depending only on $\a$ and $p$.
\end{thm}

It is easy to see that Theorem \ref{smaxbis} implies Theorem
\ref{smax}. Indeed, applying Theorem \ref{smaxbis} to $T^2$ with
$\a=0$, we get
 \be
 \big\|\,{\sup_n}^+ T^{2n}(x)\big\|_p\le C'_p\,\|x\|_p\quad\mbox{and}
 \quad\big\|\,{\sup_n}^+ T^{2n+1}(x)\big\|_p\le C'_p\,\|Tx\|_p
 \le C'_p\,\|x\|_p
 \ee
for every $x\in L_p(\M)$. Thus (\ref{smax1}) follows.

\medskip

As already said before, our proof of Theorem \ref{smaxbis} will
follow the pattern set up by Stein. The main steps are as follows.
First, using Theorem \ref{max}, we show that (\ref{smaxbis1})
holds for all complex $\a$ whose real part is greater than 1. Then
with the help of the discrete Littlewood-Paley function, we deduce
(\ref{smaxbis1}) for $p=2$ and for all non positive integers $\a$.
It is this $L_2$ result which demands the symmetry of $T$. For
interpolation we need to modify slightly this $L_2$ result into
another one, i.e. to prove (\ref{smaxbis1}) for $p=2$ again and
for all complex $\a$ whose real part is of the form $-m+1/2$ with
$m\in\nat$. Finally, complex interpolation permits us to conclude
the proof.

\medskip

We will use the following elementary properties of the $A_n^\a$:
for all $\a\in\comp$ and $\b>-1$
 \beq\label{An}
 A_n^\a=\sum_{k=0}^nA_k^{\a-1}\,,\quad
 A_n^\a -A_n^{\a-1}=A_{n-1}^{\a}\,,\quad A_n^\b\le c_\b\, (n+1)^\b\,,
 \eeq
where $c_\b$ is a positive constant depending only on $\b$ (noting
that $A_n^\b>0$ when $\b>-1$). The reader is referred to
\cite[Vol.I, Chapter III.1]{zy} for these formulas. The following
estimates on $A_n^\a$ are also well-known.

\begin{lem}\label{A estimate}
 Let $\a=\b+i\g\in\comp$.
\begin{enumerate}[{\rm i)}]
\item If $\b=m+r$ with $m\in\ent$ and $0<r<1$, then
 $|A_n^\a|\le \exp(c_r\g ^2)\,|A_n^\b|.$
\item If $\b>-1$, then
 $|A_n^\a|\le \exp(c_\b\g^2)\,A_n^\b.$
\end{enumerate}
\end{lem}

\pf We have
 $${\rm Log}\,\frac{A_n^\a}{A_n^\b}=\sum_{k=1}^n{\rm Log}(1+i\,
 \frac{\g}{\b+k}).$$
Writing $\displaystyle \frac{A_n^\a}{A_n^\b}=e^{u+iv}$ with $u,
v\in\real$, we see that
 $$u=\sum_{k=1}^n{\rm Re}\big({\rm Log}(1+i\,\frac{\g}{\b+k})\big)
 \le \frac{\g^2}{2}\,\sum_{k=1}^n\frac{1}{(\b+k)^2}
 \le c_r\,\g^2.$$
The second part is proved similarly. \cqd

\medskip

The following is an easy consequence of Theorem \ref{max}.

\begin{lem}\label{max big 1}
 Let $\a=\b+i\,\g$ with $\b>1$.
 Then for  all $1<p<\8$
 $$\big\|\,{\sup_n}^+ M^\a_n(x) \big\|_p\le C_{p,\b}\,
 \exp(c_\b\n\g^2)\|x\|_p\ ,\quad \forall\; x\in L_p(\M).$$
\end{lem}

\pf Without loss of generality, we assume  $x\ge0$. Let
$(y_n)\subset L_{p'}^+(\M)$. Using Lemma \ref{A estimate} and
(\ref{An}), we have
 \be
 |\t(M_n^\a(x)y_n)|
 &\le& (n+1)^{-\b}\sum_{k=0}^{n}|A_{n-k}^{\a-1}|\,\t[T^k(x)y_n]\\
 &\le& C_\b\,\exp({c_\b\g^2})\,(n+1)^{-\b}
 \sum_{k=0}^{n}(n-k+1)^{\b-1}\,\t[T^k(x)y_n]\\
 &\le& C_\b\,\exp({c_\b\g^2})(n+1)^{-1}\sum_{k=0}^{n}\,
 \t[T^k(x)y_n]=C_\b\,\exp({c_\b\g^2})\t[M_n(x)y_n ].
 \ee
Therefore,
 $$\sum_{n\ge 0} |\t(M_n^\a(x)y_n)|\le
 C_\b\,\exp({c_\b\,\g^2})\,\sum_{n\ge 0}\t[M_n(x)y_n ].$$
Taking the supremum over all $(y_n)\subset L_{p'}^+(\M)$ such that
$\|\sum_ny_n\|_{p'}\le 1$ and using Proposition \ref{vectorLp0},
iii) and Theorem \ref{max}, we deduce the assertion. \cqd

\medskip

Our next step is to prove a similar maximal inequality in
$L_2(\M)$ for $M_n^\a$ with $\a=-m+1/2$ and $m\in\nat$. To this
end we will need the following inequality on the discrete
Littlewood-Paley square function.  Let
 $$B_k^m=k(k-1)\cdots(k-m+1)\ \mbox{for}\; m\le k \quad\mbox{and}\quad
  B_k^m=0\ \mbox{for}\; m> k.$$

\begin{lem}\label{littlewood}
Let $m\in\nat$. Then for  every selfadjoint operator $x\in
L_2(\M)$
 $$\t\big[\sum_{k\ge m}k\big(B_{k-1}^{m-1}\,\D_k^m(x)\big)^2\big]
 \le C_m\, \t(x^2).$$
\end{lem}

\pf Note that if $x$ is selfadjoint, so is $\D_k^m(x)$ (recalling
that $T$ is positive). Moreover, $\D_k^m$, considered as an
operator on $L_2(\M)$, is also selfadjoint by virtue of (0.IV).
Fix a unit selfadjoint $x\in L_2(\M)$. We have
 \begin{eqnarray*}
 \t\big[\sum_{k\ge m}k\big(B_{k-1}^{m-1}\,\D_k^m(x)\big)^2\big]
 &=&\sum_{k\ge m}k\,(B_{k-1}^{m-1})^2\,\|\D_k^m(x)\|_2^2\\
 &=&\sum_{k\ge m}k\, (B_{k-1}^{m-1})^2\,
 \langle x,\,(\D_k^m)^2(x)\rangle,
 \end{eqnarray*}
where $\langle\,,\,\rangle$ stands for the scalar product on
$L_2(\M)$. Observe the following easily checked formula
 $$\D_k^m=(T-1)^m\,T^{k-m}\;,\quad \forall\;k\ge m.$$
Let
 $\displaystyle T=\int_0^1\l\,de_\l$
be the spectral resolution of $T$ on $L_2(\M)$ (recalling that $T$
is a selfadjoint positive contraction on $L_2(\M)$). Then
 $$\langle x,\,(\D_k^m)^2(x)\rangle
 =\int_0^1(1-\l)^{2m}\,\l^{2k-2m}\,d\langle x,\,e_\l x\rangle\ .$$
Since $d\langle x,\,e_\l x\rangle$ is a probability measure on
$[0, 1]$, it remains to estimate
 \begin{eqnarray*}
 \sum_{k\ge m}k\, (B_{k-1}^{m-1})^2\,(1-\l)^{2m}\,\l^{2k-2m}
 &\le& C_m + C_m\,(1-\l)^{2m}\,\sum_{k\ge 2m-1}
 B_{k}^{2m-1}\,\l^{2k-2m}\\
 &\le& C_m + C_m\,\l^{2m-2}(1-\l)^{2m}\,\sum_{k\ge 2m-1}
 B_{k}^{2m-1}\,(\l^2)^{k-2m+1}\\
 &\le & C_m + C_m\,\l^{2m-2}(1-\l)^{2m}\,(1-\l^2)^{-2m}\le C_m.
 \end{eqnarray*}
Therefore the lemma is proved.\cqd

\begin{lem}\label{max L2}
 Let $m\in\nat$. Then
 $$\big\|\,{\sup_n}^+\big((n+1)^m\,\D_n^m(x)\big) \big\|_2\le
 C_m\,\|x\|_2,\quad\forall\;x\in L_2(\M).$$
\end{lem}

\pf It suffices to show this for a positive $x\in L_2(\M)$. To
this end let us first observe the following formula
 $$\sum_{k=m}^n B_{k}^{m+1}\,(a_k-a_{k-1})= B_n^{m+1}\,a_n -
 (m+1)\,\sum_{k=m}^{n-1} B_{k}^{m}\,a_k.$$
Applying this to $a_k=\D_k^m(x)$, we deduce
 \begin{equation}\label{reduction}
 B_{n}^{m+1}\D_n^m(x)=(m+1)\,\sum_{k=m}^{n-1} B_{k}^{m}\,\D_k^m(x)
 +\sum_{k=m}^{n} B_{k}^{m+1}\,\D_k^{m+1}(x).
 \end{equation}
Now let $(y_n)\subset L^+_2(\M)$. Using the convexity of the
operator-valued function $x\mapsto |x|^2$, we have (recalling that
$\D_k^m(x)$ is selfadjoint)
\begin{eqnarray*}
 \Big|\t\big(\frac{1}{n+1}\sum_{k=m}^{n-1}
 B_{k}^{m}\,\D_k^m(x)y_n\big)\Big|
 &\le&\t\big(\frac{1}{n+1}\sum_{k=m}^{n-1}B_{k}^{m}\,|\D_k^m(x)|\,y_n\big)\\
 &\le&\t\big[\big(\sum_{k=m}^{n-1}(k+1)^{-1}\big|B_{k}^{m}\,
 \D_k^m(x)\big|^2\big)^{\frac{1}{2}}\,y_n\big]\\
 &\le&\t\big[\big(\sum_{k=m}^{\8}(k+1)^{-1}\big|B_{k}^{m}\,
 \D_k^m(x)\big|^2\big)^{\frac{1}{2}}\,y_n\big].
\end{eqnarray*}
Hence, by Lemma \ref{littlewood},
\begin{eqnarray*}
 \Big|\sum_{n}\t\big(\frac{1}{n+1}\sum_{k=m}^{n-1}
  B_{k}^{m}\,\D_k^m(x)y_n\big)\Big|
 &\le&\t\big[\big(\sum_{k=m}^{\8}(k+1)^{-1}\big|B_{k}^{m}\,
 \D_k^m(x)\big|^2\big)^{\frac{1}{2}}\,\sum_n y_n\big]\\
 &\le&\big\|\big(\sum_{k=m}^{\8}(k+1)^{-1}\big|B_{k}^{m}\,
 \D_k^m(x)\big|^2\big)^{\frac{1}{2}}\big\|_2\;
 \big\|\sum_n y_n\big\|_2\ \\
 &\le& C_m\,\|x\|_2\, \|\sum_n y_n\|_2.
\end{eqnarray*}
Similarly,
$$
 \Big|\sum_{n}\t\big(\frac{1}{n+1} \sum_{k=m}^{n}
  B_{k}^{m+1}\,\D_k^{m+1}(x)y_n\big)\Big|
 \le C_m\,\|x\|_2\, \big\|\sum_n y_n\big\|_2\ .
$$
Combining these inequalities with (\ref{reduction}), we deduce
$$
 \Big|\sum_{n}\t\big(\frac{1}{n+1}\, B_{n}^{m+1}\D_n^m(x)y_n\big)\Big|
 \le C_m\,\|x\|_2\, \big\|\sum_n y_n\big\|_2\ ;
$$
whence
$$\big\|\,{\sup_n}^+\big(\frac{1}{n+1}\, B_{n}^{m+1}\,\D_n^m(x)\big)\big\|_2\le
 C_m\,\|x\|_2\ .$$
This is clearly equivalent to the desired inequality.\cqd

\begin{lem}\label{max combination}
 Let  $x=(x_n)\in L_p(\M;\el_\8)$ and $(z_{n,k})_{n,k}\subset \comp$. Then
 $$\big\|\,{\sup_n}^+\sum_kz_{n,k}x_k\big\|_p\le
 \sup_n\big(\sum_k|z_{n,k}|\big)\,\big\|\,{\sup_k}^+x_k\big\|_p\ .$$
\end{lem}

\pf This is easy. Indeed, given a factorization of $x$ as
$x_k=ay_kb$, we have
 $$\sum_kz_{n,k}x_k=a\,\big(\sum_kz_{n,k}y_k\big)\,b\ .$$
Thus
\begin{eqnarray*}
 \big\|\,{\sup_n}^+\sum_kz_{n,k}x_k\big\|_p
 &\le&\|a\|_{2p}\,\|b\|_{2p}\,\sup_n\big\|\sum_kz_{n,k}y_k\big\|_\8\\
 &\le&\|a\|_{2p}\,\|b\|_{2p}\,\sup_k\|y_k\|_\8\,\sup_n\sum_k|z_{n,k}|\
 .
 \end{eqnarray*}
This implies the desired inequality. \cqd

\begin{lem}\label{max L2bis}
 Let $\a=\b+i\,\g$ such that $\b=-m+1/2$ with $m\in\nat$. Then
 $$\big\|\,{\sup_n}^+M_n^\a(x)\big\|_2\le
 C_m\exp(5\g^2)\,\|x\|_2,\quad\forall\;x\in L_2(\M)\ .$$
\end{lem}

\pf Given $n\in\nat$ choose $d_n\in\nat$ such that $n/2 +1\le
d_n<n/2 +3$. Then by successive use of the Abel summation, we
obtain
 \begin{eqnarray*}
 S_n^\a
 &=&\sum_{k=0}^{n}A_{n-k}^{\a-1}\,T^k
 =\sum_{k=0}^{d_n-1}A_{n-k}^{\a-1}\,T^k
 +\sum_{k=d_n}^{n}A_{n-k}^{\a-1}\,T^k\\
 &=&\sum_{k=0}^{d_n-1}A_{n-k}^{\a-1}\,T^k
 +A_{n-d_n}^{\a}\,T^{d_n-1}
 +\sum_{k=d_n}^{n}A_{n-k}^{\a}\,\D_k\\
 &=&\sum_{k=0}^{d_n-1}A_{n-k}^{\a-1}\,T^k
 +A_{n-d_n}^{\a}\,T^{d_n-1}+A_{n-d_n}^{\a+1}\,\D_{d_n-1}
 +\sum_{k=d_n}^{n}A_{n-k}^{\a+1}\,\D_k^2\\
 &\vdots&\\
 &=&\sum_{k=0}^{d_n-1}A_{n-k}^{\a-1}\,T^k
 +\sum_{j=1}^mA_{n-d_n}^{\a +j-1}\,\D_{d_n-1}^{j-1}
 +\sum_{k=d_n}^{n}A_{n-k}^{\a+m-1}\,\D_k^m\ .
 \end{eqnarray*}
Therefore, by triangle inequality and Lemma \ref{max combination},
we get
 $$\big\|\,{\sup_n}^+M_n^\a(x)\big\|_2\le I\times II,$$
where
 $$I=\sum_{j=0}^m\big\|\,{\sup_n}^+(n+1)^j\D_n^j(x)\big\|_2$$
and
 $$II=\sup_n\frac{1}{(n+1)^{\b}}\max\left\{
 \sum_{k=0}^{d_n-1}|A_{n-k}^{\a-1}|,\;
 \max_{1\le j\le m}\,\frac{|A_{n-d_n}^{\a +j-1}|}{d_n^{j-1}}\,,\;
 \sum_{k=d_n}^{n}\frac{|A_{n-k}^{\a+m-1}|}{(k+1)^{m}}
 \right\}.$$
By Lemma \ref{max L2}, $I\le C_m\,\|x\|_2$. On the other hand, by
Lemma \ref{A estimate}, i) ($c_{1/2}\le 5$ with $r=1/2$ there)
 $$II\le C_m\,\exp({5\g^2})\,\sup_n\frac{1}{(n+1)^{\b}}\max\left\{
 \sum_{k=0}^{d_n-1}|A_{n-k}^{\b-1}|,\;
 \max_{1\le j\le m}\,\frac{|A_{n-d_n}^{\b +j-1}|}{d_n^{j-1}}\,,\;
 \sum_{k=d_n}^{n}\frac{|A_{n-k}^{\b+m-1}|}{(k+1)^{m}}
 \right\}.$$
Now using the following easily verified estimate
 $$|A_{k}^{\d}|\le C_\d\, (k+1)^\d$$
for real $\d$ (see also \cite[Vol.I, Chapiter III.1]{zy}) and by
the choice of $d_n$, we get
 $$\sum_{k=d_n}^{n}\frac{|A_{n-k}^{\b+m-1}|}{(k+1)^{m}}\le
 \frac{C_m}{(n+1)^{m}}\sum_{k=1}^n\frac{1}{\sqrt k}
 \le C_m\,(n+1)^{\b}\ .$$
This gives the desired estimate on the last term in the brackets
above. The other two terms can be estimated similarly. Therefore,
$II\le C_m$. Putting together all preceding inequalities yields
the lemma.\cqd

\medskip

Now we are in a position to prove Theorem \ref{smaxbis}.

\medskip

\n{\em Proof of Theorem \ref{smaxbis}.} Write $\a=\b+i\g$ with
$\b, \g\in\real$. Choose $\th\in(0,1),\; q\in(1,\8)$, $m\in\ent$
and $b>\max(\b, 1)$ such that
 $$\frac{1}{p}=\frac{1-\th}{2} + \frac{\th}{q}
 \quad\mbox{and}\quad
 \b=(1-\th)a +\th\,b\quad\mbox{with}\ a=m+\frac{1}{2}.$$
Let $x\in L_p(\M)$ and $y=(y_n)$ be a finite sequence in
$L_{p'}(\M)$ with $\|x\|_p<1$ and $\|y\|_{L_{p'}(\M; \el_1)}<1$.
Define
 $$f(z)=u\,\big|x\big|^{\frac{p(1-z)}{2}+ \frac{pz}{q}}\, ,
 \quad z\in\comp,$$
where $x=u|x|$ is the polar decomposition of $x$. On the other
hand, by Proposition \ref{c-interp-vectorLp}, there is a function
$g=(g_n)_n$ continuous on the strip $\{z\in\comp\;:\; 0\le{\rm
Re}(z)\le 1\}$ and analytic in the interior such that $g(\th)=y$
and
 $$\sup_{t\in\real}\,\max\Big\{\|g(i\,t)\|_{L_{2}(\M; \el_1)}\ ,\
 \|g(1+i\,t)\|_{L_{q'}(\M; \el_1)}\Big\}<1.$$
Now define
 $$F(z)=\exp\big(\d(z^2-\th^2))\,\sum_n\t\big[M_n^{(1-z)a + zb+i\,\g}(f(z))\,
 g_n(z)\big],$$
where $\d>0$ is a constant to be specified. $F$ is a function
analytic in the open strip $\{z\in\comp\;:\; 0<{\rm Re}(z)< 1\}$.
Applying Lemma \ref{max big 1} when $m\ge1$ and Lemma \ref{max
L2bis} when $m\le0$, we have
 \begin{eqnarray*}|F(it)|
 &\le& \exp\big(\d(-t^2-\th^2))\,\big\|\big(M_n^{a +i(-ta +t
 b+\g)}(f(it))\big)_n\big\|_{L_{2}(\M; \el_\8)}
 \,\big\|g(i\,t)\big\|_{L_{2}(\M; \el_1)}\\
 &\le& C_\a\, \exp\big((-\d+c_{\b, b,\g})t^2-\d\th^2)\,\|f(it)\|_2
 \le C_\a\, \exp\big((-\d+c_{\b, b,\g})t^2-\d\th^2)\ .
 \end{eqnarray*}
Similarly, by Lemma \ref{max big 1},
 $$|F(1+it)|\le C_{\a, q}\, \exp\big((-\d+c'_{\b, b,\g})t^2+\d(1-\th^2))\ .$$
Choosing $\d$ bigger than $\max(c_{\b, b,\g},\; c'_{\b, b,\g})$,
we get
 $$\sup_{t\in\real}\,\max\Big\{|F(i\,t)|\ ,\
 |F(1+i\,t)|\Big\}\le C_{\a, p}.$$
Therefore, by the maximum principle,
 $|F(\th)|\le  C_{p,\b,b,\g}.$
Namely,
 $$\big|\sum_n\t\big[M_n^{\a}(x)\,
 y_n\big]\big|\le C_{\a, p}\ .$$
This yields (\ref{smaxbis1}) and thus the theorem is proved. \cqd

\medskip

We end this section with some direct consequences of Theorem
\ref{smax}.

\begin{cor}\label{smaxc}
Let $T$ be as in Theorem \ref{smax} and $2<p<\8$. Then
 $$\big\|\big(T^n(x)\big)_n\big\|_{L_p(\M; \el_\8^c)}
 \le \sqrt{C'_p}\; \|x\|_p\;,\quad \forall\; x\in L_p(\M).$$
\end{cor}

\pf Based on Theorem \ref{smax}, this corollary is proved in the same way as
Corollary \ref{maxc}.\cqd

\medskip

By iteration, we get the following

\begin{cor}\label{smaxc-multi} Let $T_1, ..., T_d$
satisfy {\rm (0.I) - (0.IV)}.
Then for $1<p<\8$
 $$\big\|\mathop{{\sup}^+}_{n_1,..., n_d}\,T_d^{n_d}\,\cdots\,T_1^{n_1}(x)\big\|_p
 \le (C'_{p})^d\, \|x\|_p\;,\quad \forall\; x\in L_p(\M)$$
and for $2<p<\8$
 $$\big\|\big(T_d^{n_d}\,\cdots\,T_1^{n_1}(x)\big)_{n_1,..., n_d}
 \big\|_{L_p(\M; \el_\8^c(\nat^d))}
 \le (C'_{p})^{\frac{d}{2}}\, \|x\|_p\;,\quad \forall\; x\in L_p(\M).$$
 \end{cor}

By discretization, the previous maximal inequalities on
contractions imply similar ones on semigroups.

\begin{cor}\label{smaxc-sg}
{\rm i)} Let $(T_t)_{t\ge 0}$ be a semigroup verifying the
conditions {\rm (0.I) - (0.IV)}. Then for $1<p<\8$
 $$\big\|\mathop{{\sup}^+}_{t\ge 0}T_t(x)\big\|_p
 \le C'_{p}\, \|x\|_p\;,\quad \forall\; x\in L_p(\M)$$
and for $2<p<\8$
 $$\big\|\big(T_t(x)\big)_{t}
 \big\|_{L_p(\M; \el_\8^c(\real_+))}
 \le \sqrt{C'_{p}}\, \|x\|_p\;,\quad \forall\; x\in L_p(\M).$$
 {\rm ii)} A similar statement holds for $d$ such semigroups.
 \end{cor}

\section{Individual ergodic theorems}\label{individual}

In this section we apply the maximal inequalities proved in the
two previous sections to study the pointwise ergodic convergence.
To this end we first need an appropriate analogue for the
noncommutative setting of the usual almost everywhere convergence.
This is the almost uniform convergence introduced by Lance
\cite{la-erg} (see also \cite{ja1}).

\begin{defi} Let $\M$ be a von Neumann algebra equipped with
a semifinite normal
faithful trace $\t$. Let $x_n, x\in L_0(\M)$.
 \begin{enumerate}[{\rm i)}]
 \item $(x_n)$ is said to converge {\rm bilaterally almost uniformly}
$(${\rm b.a.u.} in short$)$ to $x$ if for every $\e>0$ there is a
projection $e\in \M$ such that
 $$\t(e^\perp)<\e\quad \mbox{and}\quad
 \lim_{n\to\8}\|e(x_n-x)e\|_\8=0.$$
 \item $(x_n)$ is said to converge {\rm almost uniformly}
$(${\rm a.u.} in short$)$ to $x$ if for every $\e>0$ there is a
projection $e\in \M$ such that
 $$\t(e^\perp)<\e\quad \mbox{and}\quad
 \lim_{n\to\8}\|(x_n-x)e\|_\8=0.$$
 \end{enumerate}
\end{defi}

In the commutative case, both convergences in the
definition above are equivalent to
the usual almost everywhere convergence by virtue of
Egorov's theorem. However they
are different in the noncommutative setting.
Similarly, we introduce these notions of
convergence for functions with values in $L_0(\M)$
and for nets in $L_0(\M)$.

\medskip
In order to deduce the individual ergodic theorems from  the
corresponding  maximal ergodic theorems, it is convenient to use
the subspace $L_p(\M; c_0)$ of $L_p(\M; \el_\8)$.  $L_p(\M; c_0)$
is defined as the space of all sequences $(x_n)\subset L_p(\M)$
such that there are $a,b\in L_{2p}(\M)$ and $(y_n)\subset\M$
verifying
 $$x_n=ay_nb\quad\mbox{and}\quad \lim_{n\to\8} \|y_n\|_\8=0.$$
It is easy to check that $L_p(\M; c_0)$ is a closed subspace of
$L_p(\M; \el_\8)$ and
 $$\big\|\,{\sup_n}^+ x_n\big\|_p=\inf\big\{\|a\|_{2p}\,
 \sup_{n\ge0}\|y_n\|_\8\,\|b\|_{2p}\big\},$$
where the infimum runs over all factorizations of $(x_n)$ as
above. We define similarly the subspace $L_p(\M; c_0^c)$ of
$L_p(\M; \el^c_\8)$. Note that $L_p(\M; c_0)$ (resp. $L_p(\M;
c_0^c)$) is just the closure in $L_p(\M; \el_\8)$ (resp. $L_p(\M;
\el_\8^c)$) of finite sequences in $L_p(\M)$ for $1\le p<\8$. The
definition of these spaces is readily extended to any index set
instead of $\nat$.

\medskip

The following simple lemma from \cite{dj} will be useful for our
study of individual ergodic theorems. We include a proof for
completeness.

\begin{lem}\label{con-c0}
{\rm i)} If $(x_n)\in L_p(\M; c_0)$ with $1\le p<\8$, then $x_n$
converges b.a.u. to $0$.

{\rm ii)} If $2\le p<\8$ and $(x_n)\in L_p(\M; c_0^c)$, then $x_n$
converges a.u. to $0$.
\end{lem}

\pf  i) Let $(x_n)\in L_p(\M; c_0)$. Then there are $a, b\in
L_{2p}(\M)$ and $y_n\in \M$ such that
 $$x_n=ay_nb\quad\mbox{and}\quad \|a\|_{2p}<1,\ \|b\|_{2p}<1,\;
  \lim_{n\to\8} \|y_n\|_\8=0.$$
We can clearly assume $a, b\ge0$.  Let $e'$ be a spectral
projection of $a$ such that
 $$\t({e'}^\perp)<\frac{\e}{2}\quad\mbox{and}\quad
 \|e'a\|_\8\le \big(\frac{2}{\e}\big)^{\frac{1}{2p}}\;.$$
Similarly, we find a spectral projection $e''$ of $b$. Set
$e=e'\wedge e''$. Then
 $$\t(e^\perp)\le \t({e'}^\perp) + \t({e''}^\perp)<\e$$
and
 $$\|ex_ne\|_\8\le \|ea\|_\8\,\|y_n\|_\8\, \|be\|_\8
 \le\|y_n\|_\8\,\|e'a\|_\8\, \|be''\|_\8
 \le \big(\frac{2}{\e}\big)^{\frac{1}{p}}\,\|y_n\|_\8\;.$$
Thus $\lim_n\|ex_ne\|_\8=0$, and so $x_n\to 0$ b.a.u.

ii) The proof of this part is similar and left to the reader. \cqd

\medskip

Now let $T$ be a linear map satisfying the conditions (0.I) -
(0.III). Let $(M_n)_n$ denote the ergodic averages of $T$. Recall
that $\F_p$ denotes the fixed point subspace of $T$ in $L_p(\M)$
and $F$ the projection from $L_p(\M)$ onto $\F_p$ (see section~1).

\begin{thm}\label{maxc0}
Let $T$ be a map on $\M$ satisfying {\rm(0.I) - (0.III)}. Let $1<
p<\8$ and $x\in L_p(\M)$. Then $\big(M_n(x)-F(x)\big)_n\in
L_p(\M;c_0)$ Moreover, if $p>2$,   $\big(M_n(x)-F(x)\big)_n\in
L_p(\M;c_0^c)$.
\end{thm}

\pf  Let $x\in L_p(\M)$. Since $(I-T)\big(L_1(\M)\cap
L_\8(\M)\big)$ is dense in $\F_p^\perp$, there are
$x_k\in(I-T)\big(L_1(\M)\cap L_\8(\M)\big)$ such that
 \be
 \lim_{k\to\8}\|x-F(x)-x_k\|_p=0.
 \ee
By Theorem \ref{max},
 $$\big\|\big(M_n(x)-F(x)-M_n(x_k)\big)_n\big\|_{L_p(\M;\el_\8)}
 \le C_p\,\|x-F(x)-x_k\|_p\,.$$
Thus
 $$\lim_{k\to\8}\big(M_n(x_k)\big)_n=\big(M_n(x)-F(x)\big)_n
 \quad\mbox{in}\quad L_p(\M;\el_\8).$$
Since $L_p(\M;c_0)$ is closed in $L_p(\M;\el_\8)$, it suffices to
show $\big(M_n(x_k)\big)_n\in L_p(\M;c_0)$ for every $k$. To this
end consider an arbitrary $z\in (I-T)\big(L_1(\M)\cap
L_\8(\M)\big)$. Let  $y\in L_1(\M)\cap L_\8(\M)$ such that
$z=y-T(y)$. Then
 $$M_n(z)=\frac{1}{n+1}\,\big[y-T^{n+1}(y)\big].$$
Since $z\in L_q(\M)$ for any $1<q<\8$ we deduce from Theorem
\ref{max} that $\big(M_n(z)\big)_n$ belongs to $L_q(\M;\el_\8)$.
Choose a $q\in(1, p)$. Then by Proposition
\ref{c-interp-vectorLp}, for any $m<n$
 \be
 \big\|\mathop{{\sup}^+}_{m\le j\le n}M_j(z)\big\|_p
 &\le&\sup_{m\le j\le n}\big\|M_j(z)\big\|_\8^{1-\frac{q}{p}}\,
 \big\|\mathop{{\sup}^+}_{m\le j\le n}M_j(z)\big\|_q^{\frac{q}{p}}\\
 &\le&\big[\frac{2\|y\|_\8}{m+1}\big]^{1-\frac{q}{p}}\,
 \big\|\mathop{{\sup}^+}_{j\ge1}M_j(z)\big\|_q^{\frac{q}{p}}\,.
 \ee
Let $\vec{z}^{(k)}$ denote the finite sequence $(M_0(x),...,
M_k(x), 0,...)$. The inequality above shows that the sequence
$(\vec{z}^{(k)})_{k\ge 0}$ converges to $(M_n(z))_n$ in
$L_p(\M;\el_\8)$ as $k\to \8$. Thus $\big(M_n(z)\big)_n\in
L_p(\M;c_0)$, as wanted.

The second part can be similarly proved. Now we use Corollary
\ref{maxc} and the analogue for the spaces $L_p(\M;c_0^c)$ of
Proposition \ref{c-interp-vectorLp}.  \cqd

\medskip

The following is an extension of Yeadon's noncommutative
individual ergodic theorem \cite{ye1} to all $L_p(\M)$ with
$1<p<\8$.

\begin{cor}\label{con}
Let $T$ be a map satisfying the conditions {\rm (0.I) - (0.III)}.
Let $1< p<\8$ and $x\in L_p(\M)$. Then $(M_n(x))_n$ converges to
$F(x)$ b.a.u. for $1< p\le 2$ and a.u. for $2< p<\8$.
\end{cor}

\pf This is an immediate consequence of Lemma \ref{con-c0} and
Theorem \ref{maxc0}.\cqd

\medskip
\n{\bf Remark.}  Corollary \ref{con} can also be proved by using
Yeadon's theorem. This is however not the case for the multiple
individual ergodic theorem below. We refer to \cite{Ska} for
multiple ergodic theorems for commuting operators.

\begin{rk}\label{au convergence for p=2}
Again using Yeadon's theorem, one can prove that the convergence in
Theorem \ref{con} is a.u. for $p=2$.
\end{rk}

\pf  Fix $x\in L_2(\M)$ and $\e>0$. By decomposing $x$ into its
real and imaginary part, we can assume $x$ selfadjoint. Let
$(\e_n)$ and $(\d_n)$ be two sequences of small positive numbers.
Then  for each $m\ge1$ there are $w_m\in L_2(\M)\cap L_\8(\M)$ and
$z_m\in L_2(\M)$ such that
 $$x=F(x)+ y_m+z_m \quad\mbox{with}\quad y_m=w_m-T(w_m)\quad\mbox{and}
 \quad \|z_m\|_2<\d_m.$$
Since $x$ is selfadjoint,  $w_m, y_m$ and $z_m$ can be chosen
selfadjoint too. We have
 $$M_n(x)-F(x)=M_n(y_m) + M_n(z_m)$$
and
 $$\|M_n(y_m)\|_\8\le\frac{2}{n+1}\,\|w_m\|_\8.$$
Now we apply Yeadon's weak type (1,1) inequality (Lemma
\ref{yeadon}) to $z_m^2$. Thus there is a projection $e_m$ such
that
 $$\sup_n\big\|e_mM_n(z_m^2)e_m\big\|_\8\le \e_m^2\quad\mbox{and}\quad
 \t(e_m^\perp)<\e_m^{-2}\,\t(z_m^2)\le\e_m^{-2}\,\d_m^2\;.$$
By Kadison's Cauchy-Schwarz inequality \cite{kad-cs}, we get
 $$\big\|M_n(z_m)e_m\big\|^2_\8\le \big\|e_mM_n(z_m^2)e_m\big\|_\8\;.$$
Thus
 $$\sup_n\big\|M_n(z_m)e_m\big\|_\8\le \e_m\;.$$
Let $e=\bigwedge_m\,e_m$. Then
 $$\t(e^\perp)\le \sum_{m\ge1}\e_m^{-2}\, \d_m^2<\e $$
provided $\e_m$ and $\d_m$ are appropriately chosen. On the other
hand, by the preceding inequalities, we deduce
 \be
 \big\|e[M_n(x)-F(x)]e\big\|_\8
 &\le& \frac{2}{n+1}\|w_m\|_\8 + \|eM_n(z_m)e\|_\8\\
 &=&\frac{2}{n+1}\|w_m\|_\8 + \|e[e_mM_n(z_m)e_m]e\|_\8
 \le \frac{2}{n+1}\|w_m\|_\8 + \e_m.
 \ee
 It then follows that
 $${\mathop{\rm lim\,sup}_{n\to\8}}\,\big\|e[M_n(x)-F(x)]e\big\|_\8
 \le \e_m.$$
Since $\lim_m\e_m=0$, we get that $\lim_{n\to0}
\|e[M_n(x)-F(x)]e\|_\8=0$. Hence, $M_n(x)$ converges  to $F(x)$
b.a.u.. \cqd

\medskip

We pass to  the multiple version of Theorem \ref{maxc0} and
Corollary \ref{con}. Let $T_1, ...\,,T_d$ be $d$ maps satisfying
{\rm (0.I) - (0.III)}. As before, set
 $$M_{n_1,...\,,n_d}=\big[\prod_{j=1}^d\frac{1}{n_j+1}\big]
 \sum_{k_d=0}^{n_d}
 \;\cdots\;\sum_{k_1=0}^{n_1} T_d^{k_d}\,\cdots\,T_1^{k_1}\,.$$
Let $F_j$ be the projection onto the fixed point subspace of
$T_j$.

\begin{thm}\label{con-multi}
Let $T_1, ...\,,T_d$ be $d$ maps satisfying {\rm (0.I) - (0.III)}.
Let $1<p<\8$ and $x\in L_p(\M)$. Then
 $$\big(M_{n_1,\,...,\,n_d}(x)-F_d\,\cdots\,
 F_1(x)\big)_{n_1, \,...,\,n_d\ge 1}\in L_p(\M;
 c_0(\nat^d))$$
and if $p>2$,
 $$\big(M_{n_1,\,...,\,n_d}(x)-F_d\,\cdots
 \,F_1(x)\big)_{n_1, \,...,\,n_d\ge 1}\in L_p(\M;
 c_0^c(\nat^d)).$$
Consequently, $M_{n_1,...\,,n_d}(x)$ converges b.a.u. to
$F_d\,\cdots\, F_1(x)$ as $n_1, ..., n_d$ tend to $\8$. Moreover,
the convergence is a.u.  in the case of $p>2$.
\end{thm}

\pf This proof is similar to that of Theorem \ref{maxc0}, modulo
an iteration argument. We consider only the typical case $d=2$.
Note that
 $$M_{n_1,n_2}=M_{n_2}(T_2)M_{n_1}(T_1).$$
Fix $x\in L_p(\M)$ and decompose $x$ as $x=F_1(x)+y_k+u_k$ with
 $$y_k \in (I-T_1)\big(L_1(\M)\cap L_\8(\M)\big),\quad
 u_k\in L_p(\M),\quad  \|u_k \|_p\le \frac{1}{k}.$$
Similarly, we decompose $F_1(x)$ with respect to $T_2$:
 $F_1(x)=F_2(F_1(x))+z_k+v_k$ with
 $$z_k\in (I-T_2)\big((L_1(\M)\cap L_\8(\M)\big),\quad
 v_k\in L_p(\M),\quad  \|v_k\|_p\le \frac{1}{k}.$$
Applying successively $M_{n_1}(T_1)$ to $x$ and $M_{n_2}(T_2)$ to
$F_1(x)$, we get
 $$M_{n_1,n_2}(x)-F_2  F_1(x)
 =M_{n_1,n_2}(y_k)+M_{n_1,n_2}(u_k)
 +M_{n_2}(T_2)(z_k)+M_{n_2}(T_2)(v_k)\,.$$
By Corollary \ref{max-multi},
 $$\big\|\mathop{{\sup}^+}_{n_1, n_2}\,M_{n_1,n_2}(u_k)\big\|_p\le C_p^2\,\|u_k\|_p
 \le \frac{C_p^2}{k}\ \to 0\ \mbox{as}\ k\to\8.$$
Similarly,
  $$\lim_{k\to\8}\big\|\,{\sup_{n_2}}^+M_{n_2}(T_2)(v_k)\big\|_p=0.$$
Therefore, as in the proof of Theorem \ref{maxc0}, we need only to
show
 $$\big(M_{n_1,n_2}(y_k)\big)_{n_1, n_2\ge 1}\in L_p(\M;
 c_0(\nat^2))\quad\mbox{and}\quad
  \big(M_{n_2}(z_k)\big)_{n_2\ge 1}\in L_p(\M;
 c_0).$$
Theorem \ref{maxc0} implies the latter. The former is proved by
the arguments in the proof of Theorem \ref{maxc0}. Thus the first
part of the theorem is proved.  The second on $p>2$ is left to the
reader.

Then applying Lemma \ref{con-c0} to multiple sequences, we deduce
the announced pointwise  multiple ergodic convergence. \cqd

\medskip

We have the following stronger convergence result  for symmetric $T$.

\begin{thm}\label{scon}
Let $T$ be a map satisfying {\rm (0.I) - (0.IV)}. Assume further
that $T$ is positive as an operator on $L_2(\M)$.  Let $1< p<\8$
and $x\in L_p(\M)$. Then $(T^n(x)-F(x))_n$ belongs to $L_p(\M;
c_0)$ and to $L_p(\M; c_0^c)$ if additionally $p>2$. Consequently,
$T^n(x)$ converges to $F(x)$ b.a.u. for $1<p\le 2$ and a.u. for
$2<p<\8$.
\end{thm}

\pf Let us first treat the case $p=2$. Write the spectral
decomposition of $T$:
 $$T=\int_0^1 \l\, d\,e_\l\;.$$
Note that for any $x\in\overline{(I-T)(L_2(\M))}$
 $$\lim_{\l\to 1}e_{\l}(x)=x\quad \mbox{in}\quad L_2(\M).$$
Given $x\in L_2(\M)$ choose  $x_k\in (I - T)(L_2(\M))$ such that
$\lim\|x-F(x)-x_k\|_2=0$. Then
 $\lim_{\l\to 1}e_{\l}(x_k)=x_k.$
Thus replacing $x_k$ by $e_{\l_k}(x_k)$ with an appropriate
$\l_k\in(0,1)$, we can assume that $x_k=e_{\l_k}(y_k)$ for some
$y_k\in L_2(\M)$. Then
 $$T^n(x_k)=\int_0^{\l_k} \l^n\,d\,e_\l(y_k)\;;\quad
 \mbox{whence}\quad
 \|T^n(x_k)\|_2\le\l_k^n\,\|y_k\|_2\;.$$
It then follows that $\big(T^n(x_k)\big)_n\in L_2(\M; c_0)$ for
every $k$, and so by Theorem \ref{smax},
$\big(T^n(x)-F(x)\big)_n\in L_2(\M; c_0)$.

To treat the general case  we first claim that
 $$\lim_n\|T^n(x)-F(x)\|_p=0,\quad \forall\; x\in L_p(\M).$$
Indeed,  the preceding argument shows that this is true for $p=2$.
Now  let $2<p<\8$ and  $x\in L_1(\M)\cap\M$. By interpolation,
 $$\|T^n(x)-F(x)\|_p\le \|T^n(x)-F(x)\|_\8^{1-\frac{2}p}\,
 \|T^n(x)-F(x)\|_2^{\frac{2}p}\,;$$
whence $\lim_n\|T^n(x)-F(x)\|_p=0$. Then our claim in the case
$p>2$ follows from the density of $L_1(\M)\cap\M$ in $L_p(\M)$.
The case $p<2$ is proved similarly.

Now we can easily finish the proof of the theorem. Let $x\in
L_p(\M)$. Fix $k\in\nat$. Then  Theorem \ref{smax} and the claim
above imply
 $$\lim_{k\to\8}\big\|\big(T^n\big(T^k(x)-F(x)\big)\big)_n
 \big\|_{L_p(\M; \el_\8)}=0.$$
Note that $T^n\big(T^k(x)-F(x)\big)=T^{n+k}(x)-F(x)$, and so the
sequence $\big(T^n\big(T^k(x)-F(x)\big)\big)_{n\ge0}$ can be
considered as the rest of $\big(T^n(x)-F(x)\big)_{n\ge0}$ starting
from the $k$-th coordinate. It follows that
$\big(T^n(x)-F(x)\big)_{n\ge0}\in L_p(\M; c_0)$.

In a similar way, using Corollary \ref{smaxc}, we show that
$\big(T^n(x)-F(x)\big)_n\in L_p(\M; c_0^c)$ for any $x\in L_p(\M)$
with $p>2$. \cqd

\medskip
\n{\bf Remark.} If we remove the additional assumption that $T$ is
a positive operator on $L_2(\M)$ in Theorem \ref{scon}, then for
any $x\in L_p(\M)$ the two subsequences $(T^{2n}(x))_n$ and
$(T^{2n+1}(x))_n$ still converge b.a.u.; however, their limits are
not equal in general.

\medskip

We end this section with the pointwise ergodic theorems for
semigroups. Let $(T_t)_{t\ge0}$ be a semigroup satisfying {\rm
(0.I) - (0.III)}. We denote again by $F$ the projection from
$L_p(\M)$ onto the fixed point subspace of $(T_t)_{t\ge0}$.

\begin{samepage}\begin{thm}\label{con-sg} Let $(T_t)_{t\ge0}$ be a semigroup
with {\rm (0.I) - (0.III)}. Let $(M_t)_{t>0}$ denote the ergodic
averages of $(T_t)_{t\ge0}$. Let $1< p<\8$ and $x\in L_p(\M)$.
 \begin{enumerate}[{\rm i)}]
  \item  Then \newline
   {\rm a)} $M_t(x)$ converges to $F(x)$  b.a.u. for $1< p<2$ and a.u.
   for $2\le p<\8$ when $t\to \8$.\newline
   {\rm b)} $M_t(x)$ converges to $x$ b.a.u. for $1< p<2$ and a.u.
   for $2\le p<\8$ when $t\to 0$.
 \item Assume in addition that $(T_t)_{t\ge0}$ satisfies {\rm (0.IV)}.
 Then\newline
   {\rm a)} $T_t(x)$ converges to $F(x)$ b.a.u. for $1< p\le 2$ and a.u.
   for $2< p<\8$ when $t\to \8$.\newline
  {\rm b)} $T_t(x)$ converges to $x$ b.a.u. for $1< p\le 2$ and a.u.
  for $2< p<\8$ when $t\to 0$.
 \end{enumerate}
\end{thm}\end{samepage}

\pf The two statements a) can be proved similarly as in the
discrete case, using Theorem \ref{max-sg}. The main step here is
to obtain the semigroup analogue of Theorem \ref{maxc0} as
$t\to\8$, namely, to show that the family
$\big(M_t(x)-F(x)\big)_{t\ge 1}$ belongs to $L_p(\M; c_0([1,
\8)))$ or $L_p(\M; c_0^c([1, \8)))$. Note however that the a.u.
convergence for $p=2$ in the first statement a) is shown similarly
as Remark \ref{au convergence for p=2}, using now Remark
\ref{yeadon-sg}. We leave this part of the proof to the reader,
and will show the two statements b). (The first of them is the
noncommutative analogue of the classical Wienner local pointwise
ergodic theorem.)

\medskip

Let us first consider i), b). Let $x\in L_p(\M)$. By Lemma
\ref{con-c0}, it suffices to prove $\big(M_t(x)-x\big)_{0<t\le1}$
belongs to $L_p(\M; c_0((0,1]))$ (with respect to $t\to 0$). By
the mean ergodic theorem, we have $M_t(x)\to x$ when $t\to 0$.
Thus by a limit argument as in the proof of Theorem \ref{maxc0},
we may assume $x=M_{t_0}(y)$ for some $0<t_0<1$ and $y\in
L_p(\M)$. On the other hand, by the density of $L_1(\M)\cap\M$ in
$L_p(\M)$, we can further assume that $y\in L_1(\M)\cap\M$. Let
$0<s\le t<t_0$. Then
 $$T_s(x)-x= \frac{1}{t_0}\, \Big[\int_{t_0}^{t_0+s}T_u(y)du -
 \int_{0}^{s}T_u(y)du\Big].$$
It then follows that
 $$\|T_s(x)-x\|_\8 \le \frac{2s\|y\|_\8}{t_0}$$
and so
 $$\|M_t(x)-x\|_\8 \le \frac{2t\|y\|_\8}{t_0}\ \to 0\ \mbox{as}\ t\to 0.$$
Thus by the interpolation argument used in the proof of Theorem
\ref{maxc0}, we deduce that the family
$\big(M_t(x)-x\big)_{0<t\le1}$ belongs to $L_p(\M; c_0((0,1]))$.
This proves the first part of the statement i), b).

The second part for $p>2$ can be shown in the same way. The case
$p=2$ is dealt with similarly as Remark \ref{au convergence for
p=2} in virtue of Remark \ref{yeadon-sg}.

 ii), b) is proved similarly by using Corollary \ref{smaxc-sg}
and Lemma \ref{con-c0}. Indeed, for $x=M_{t_0}(y)$ as above, we
have already proved
 $$\lim_{t\to 0}\|T_t(x)-x\|_\8 =0.$$
Therefore, $\big(M_t(x)-x\big)_{0<t\le1}\in L_p(\M; c_0((0,1]))$.
Thus the proof of the theorem is finished.\cqd

\medskip

\n{\bf Remarks.} i)  Both statements in Part i) of Theorem
\ref{con-sg} also hold for $p=1$ because of Remark
\ref{yeadon-sg}. On the other hand, both Theorem \ref{scon} and
Theorem \ref{con-sg} admit multiple versions, similar to Theorem
\ref{con-multi}.

ii) Using Corollary \ref{max-poisson}, one sees that the ergodic
averages $M_t$ in Theorem \ref{con-sg}, i) can be replaced by
$P_t$, where $(P_t)$ is a semigroup subordinated to $(T_t)$. This
is also true for $p=1$.

\section{The non tracial case}\label{haagerup}

So far we have restricted our attention to the semifinite case
only. In this section we will extend the  previous results to
arbitrary von Neumann algebras. Despite the obvious similarity
between the statements in the semifinite and the non tracial
cases, we want to point out that the situation for type III von
Neumann algebras is more complicated. This is due to the fact that
for a state $\f$ the equality $\f(e\vee f)\le \f(e)+\f(f)$ is no
longer valid, and therefore many (Egorov type) arguments from the
previous section do not apply in this general setting. Our tool
for maximal ergodic inequalities in Haagerup noncommutative
$L_p$-spaces is an important unpublished theorem due to Haagerup,
which consists in reducing the general case to the tracial one.
For clarity, we divide this section into several subsections.

\bigskip\n{\bf 7.1  Haagerup noncommutative $L_p$-spaces}

\medskip\n
The general noncommutative $L_p$-spaces used below will be those
constructed by Haagerup \cite{haag-Lp}. Our reference is
\cite{terp}. Throughout this section $\M$ will be a von Neumann
algebra equipped with a distinguished normal faithful state $\f$,
unless explicitly stated otherwise. $L_p(\M)$ denotes the
associated noncommutative $L_p$-space ($0< p\le\8$). Recall that
$L_\8(\M)$ is just $\M$ itself and $L_1(\M)$ is the predual of
$\M$. The duality between $\M$ and $L_1(\M)$ is realized via the
distinguished tracial functional $\tr$ on $L_1(\M)$:
 $$\la x, y\ra=\tr(xy), \quad y\in L_1(\M),\  x\in \M.$$
As a normal positive functional on $\M$, $\f$ corresponds to a
positive element in $L_1(\M)$. In the sequel this element will be
always denoted by $D$, called the density of $\f$ in $L_1(\M)$.
Then $\f$ can be recovered from $D$ through the preceding duality:
 $$\f(x)=\tr(xD)=\tr(Dx),\quad x\in \M.$$
We will often use the density of $D^{\frac{1-\th}{p}}\M
D^{\frac{\th}{p}}$ in $L_p(\M)$ for any $p\in(0,\8)$ and
$0\le\th\le1$. Moreover, $D^{\frac{1-\th}{p}}\M_a
D^{\frac{\th}{p}}$ is also dense in $L_p(\M)$, where $\M_a$ is the
family of all elements in $\M$ analytic with respect to the
modular group $\s_t^\f$ of $\f$ (see \cite[Lemma 1.1]{jx-burk}).

An important link between the spaces $L_p(\M)$ is the following
external product. Let $\frac{1}{r}=\frac{1}{p}+\frac{1}{q}$ and
$x\in L_p(\M),\; y\in L_q(\M)$. Then $xy\in L_r(\M)$ and
 $$\|xy\|_r\le \|x\|_p\,\|y\|_q\;.$$
Namely,  the usual H\"older inequality extends to  Haagerup
$L_p$-spaces too. In particular, the dual space of $L_p(\M)$ is
$L_{p'}(\M)$ for $1\le p<\8$, and we have
 $$\tr(xy)=\tr(yx), \quad x\in L_p(\M),\ y\in L_{p'}(\M).$$

The definition of all vector-valued $L_p$-spaces extends verbatim
to the present setting. These include $L_p(\M; \el_\8)$, $L_p(\M;
\el_\8^c)$, $L_p(\M; c_0)$ and $L_p(\M; c_0^c)$. For instance,
$L_p(\M; \el_\8)$ consists of all sequences $x=(x_n)$ in $L_p(\M)$
which admit a factorization of the following type: there are $a,
b\in L_{2p}(\M)$ and a bounded sequence $(y_n)\subset \M$ such
that $x_n=ay_nb$ for all $n$. The norm of $x$ is then defined as
 $$\|x\|_{L_p(\M; \el_\8)}=\inf\big\{\|a\|_{2p}\,\sup_n\|y_n\|_\8\,
 \|b\|_{2p}\big\},$$
where the infimum runs over all factorizations as above. We adopt
the convention introduced in section \ref{Lpel} and denote again
this norm by $\big\|\,{\sup_n}^+ x_n\big\|_p$.  Note that $L_p(\M;
c_0)$ (resp. $L_p(\M; c_0^c)$) is again a closed subspace of
$L_p(\M; \el_\8)$ (resp. $L_p(\M; \el_\8^c)$). Similarly, given an
index set $I$ we define the analogues of these spaces for families
indexed by $I$.

All properties in section \ref{Lpel} continue to hold in the
present setting. However, for an inclusion $L_p(\N)\subset
L_p(\M)$ we will now require the existence of a normal state
preserving conditional expectation from $\M$ onto $\N$. Under this
hypothesis Remark \ref{vectorLp1bis} still holds.  Also note that
for the interpolation result in Proposition
\ref{c-interp-vectorLp} we use Kosaki interpolation theorem
\cite{kos-int}.

\bigskip\n{\bf 7.2 An extension result}

\medskip\n
 Let $T$ be a map on $\M$. We will assume that $T$ satisfies
conditions similar to (0.I) - (0.IV). More precisely, we will
consider the following properties of $T$:
 \begin{itemize}
\item[{\rm (7.I)}] $T$ is a contraction on $\M$.
\item[{\rm (7.II)}] $T$ is completely positive.
\item[{\rm (7.III)}] $\f\circ T\le \f$.
\item[{\rm (7.IV)}] $T\circ\s_t^\f=\s_t^\f\circ T$ for all
$t\in\real$
\item[{\rm (7.V)}] $T$ is symmetric with respect to $\f$,
i.e. $\f(T(y)^*x)=\f(y^*T(x))$ for all $x,y\in\M$.
 \end{itemize}

\medskip

In order to consider maximal ergodic inequalities in $L_p(\M)$, we
need first to extend a map $T$ with the properties above to a
contraction on $L_p(\M)$ for all $1\le p<\8$. The following is the
non tracial analogue of Lemma \ref{extension}.

\begin{lem}\label{extension1}
Let $T$ be a map on $\M$ satisfying {\rm (7.I)} - {\rm (7.III)}.
Define
 \be
 T_{p}\ :\ D^{\frac{1}{2p}}\,\M\,D^{\frac{1}{2p}}
 &\to&~D^{\frac{1}{2p}}\,\M\,D^{\frac{1}{2p}}\\
  D^{\frac{1}{2p}}\,x\,D^{\frac{1}{2p}}~
 &\mapsto& D^{\frac{1}{2p}}\,T(x)\,D^{\frac{1}{2p}}\;.
 \ee
Then $T_{p}$ extends to a positive contraction on $L_{p}({\M})$
for all $1 \le p < \infty$. Moreover, $T$ is normal. If
additionally $T$ verifies {\rm (7.V)}, then the extension of $T_2$
is selfadjoint on $L_2(\M)$.
\end{lem}

In fact, the complete positivity assumption can be weakened to
positivity. This result comes from \cite{jx-red}. Its proof is
much more involved than that of Lemma \ref{extension}. The main
difficulty is to show the extension property of $T_1$. This
extension is essentially a reformulation of Lemma 1.2 from
\cite{haag-normalweight} into the present setting. We refer to
\cite{jx-red} for more details. The following observation is
easily checked:

\medskip

\n{\bf Remark.} Let $T$ and $T_p$ be as in Lemma \ref{extension1}
($T_p$ also denoting the extension). Let $S=T_1^*$. Then $S$
satisfies {\rm (7.I)} - {\rm (7.III)} too. Moreover,
$S_p^*=T_{p'}$ for all $1\le p<\8$, where $S_p$ is the extension
of $S$ on $L_p(\M)$, guaranteed by Lemma \ref{extension1}. This
shows, in particular, that $T$ is normal.

\medskip

The extension in Lemma \ref{extension1} is  symmetric with respect
to the injection of $\M$ into $L_p(\M)$. We could also consider
the left extension: $xD^{1/p}\mapsto T(x)D^{1/p}$ ($x\in\M$). More
generally, for any $0\le\th\le1$ we can define
 \be
 T_{p,\th}\ :\ D^{\frac{1-\theta}{p}}\,\M\, D^{\frac{\theta}{p}}
 &\to& ~D^{\frac{1-\theta}{p}}\,\M\,D^{\frac{\theta}{p}}\\
 D^{\frac{1-\theta}{p}}\,x\,D^{\frac{\theta}{p}}~
 &\mapsto& D^{\frac{1-\theta}{p}}\,T(x)\, D^{\frac{\theta}{p}}\,.
 \ee
Note that $T_{p,1/2}$ is exactly the $T_p$ defined in Lemma
\ref{extension1}. Assume in addition that $T$ satisfies {\rm
(7.IV)}. Using the equality
 $$D^{\frac{1-\th}{p}}\,\M_a\,
 D^{\frac{\th}{p}} = \M_a\,D^{\frac{1}{p}}\;,$$
one easily checks that
 $$T_{p,\th}\Big|_{\M_a\,D^{\frac{1}{p}}}
 =T_p\Big|_{\M_a\,D^{\frac{1}{p}}}\;.$$
Thus $T_{p,\theta}$ does not depend on $\th$ (at least when
restricted to analytic elements). Consequently, $T_{p,\theta}$
extends to a contraction on $L_p(\M)$. Since we will use this
observation later, we formulate it explicitly.

\begin{rk}\label{extension3}
Let $T$ satisfy {\rm (7.I)} - {\rm (7.IV)}. Then $T_{p,\theta}$
does not depend on $\th\in [0,1]$ and extends to a positive
contraction on $L_p(\M)$.
\end{rk}

\n{\bf Convention.} In the sequel, we will denote, by the same
symbol $T$, all the maps $T_{p}$ and $T_{p,\th}$ as well as their
extensions to the $L_{p}$-spaces in Lemma \ref{extension1} and
Remark \ref{extension3}, whenever no confusion can occur.

\medskip

Let $T$ be a map on $\M$ with  (7.I) - (7.III). We will consider
again the ergodic averages:
 $$M_n\equiv M_n(T)=\frac{1}{n+1}\,\sum_{k=0}^n T^k\,.$$
All discussions in section \ref{preli} concerning the mean ergodic
theorem are still valid now. Thus $T$ is ergodic on $L_p(\M)$ for
all $1\le p\le\8$ (relative to the w*-topology for $p=\8$). We
still have the decomposition
 $$L_p(\M)=\F_p(T)\oplus \F_p(T)^\perp,$$
with $\F_p(T)=\{x\in L_p(\M):T(x)=x\}$. In the previous sections,
we used several times the fact that $(I-T)(L_1(\M)\cap\M)$ is
dense in $\F_p(T)^\perp$. Now this fact should be changed to the
following: $D^{\frac{1}{2p}}\,(I-T)(\M)\, D^{\frac{1}{2p}}$ is
dense in $\F_p(T)^\perp$. If $T$ further satisfies (7.IV), this
dense subspace can be replaced by $D^{\frac{1-\th}{p}}
(I-T)(\M_a)\, D^{\frac{\th}{p}}$ for any $\th\in[0,1]$, which is
equal to $(I-T)(\M_a)\, D^{\frac{1}{p}}$ too. The easy
verification of these facts is left to the reader. As before in
the tracial case, the projection from $L_p(\M)$ onto $\F_p(T)$
will be denoted by $F$ for any $1\le p\le\8$. Again $F$ is normal
as a map on $L_\8(\M)$.

\medskip
The discussion above is readily transferred to semigroups. Let
$(T_t)$ be a semigroup of maps on $\M$ satisfying {\rm (7.I)} -
{\rm (7.III)} (i.e. each $T_t$ satisfying {\rm (7.I)} - {\rm
(7.III)}). Then $(T_t)$ extends to a semigroup of positive
contractions on $L_p(\M)$ for $1\le p<\8$. We will assume that
$(T_t)$ is w*-continuous and $T_0$ is the identity.  Then $(T_t)$
is strongly continuous on $L_p(\M)$ for every $p<\8$. Again the
fixed point projection of $(T_t)$ is denoted by $F$. Then the mean
ergodic theorem asserts that $M_t(x)$ converges to $F(x)$ as
$t\to\8$ for all $x\in L_p(\M)$ (relative to the w*-topology for
$p=\8$), where $M_t$ denotes the ergodic averages of $(T_t)$.

\medskip

The following extends a well-known result in the commutative case
to the present situation.

\begin{rk}
 Let $T$ be a map on $\M$ verifying {\rm (7.I)} -
{\rm (7.III)}. Assume in addition that $\f\circ T=\f$. Then
$\F_\8(T)$ is a von Neumann subalgebra of $\M$ and $F$ is the
normal conditional expectation from $\M$ onto $\F_\8(T)$ such that
$\f\circ F=\f$.
\end{rk}

\pf First note that under the assumptions above, both $T$ and $F$
are unital and $F$ preserves the state $\f$. Thus $F$ is a normal
unital completely positive projection from $\M$ onto $\F_\8(T)$.
Consequently, $\F_\8(T)$ contains the unit of $\M$ and is closed
under involution, and so $\F_\8(T)$ is a w*-closed operator
system. Therefore, it remains to show that $\F_\8(T)$ is closed
under the product of $\M$.

To that end we will use the following formula from
\cite{choi-Ef} (formula (3.1) there),
 \beq\label{chef}
 F(aF(x))=F(ax) \quad\mbox{and}\quad  F(F(x)a)=F(xa),\quad
 \forall\; a\in \F_\8(T),\; x\in\M.
 \eeq
Let us consider the preadjoint of $F$, $F_*: \M_*\to \M_*$. We
claim that
 $$F_*(x\f)=F(x)\f,\quad\forall\; x\in \M.$$
Indeed, since $\f\circ F=\f$, given $y\in \M$, by (\ref{chef}) we
have
 \be
 F_*(x\f)(y)
 &=&x\f(F(y))=\f(F(y)x)=\f[F(F(y)x)]=\f[F(F(y)F(x))]\\
 &=&\f[F(y F(x))]=\f(yF(x))=[F(x)\f](y).
 \ee
Now let $a, b\in\F_\8(T)$. Then
 $F_*(ab\f)=F(ab)\f.$
On the other hand, for any $x\in\M$
 \be
 F_*(ab\f)(x)
 &=& \f(F(x)ab)=\f[F(F(x)ab)]=\f\big[F\big(F[F(x)a]b\big)\big]\\
 &=&\f[F(F(xa)b)]=\f[F(xab)]=\f(xab)=[ab\f](x).
 \ee
Hence, $F_*(ab\f)=ab\f$. It thus follows that $F(ab)\f=ab\f$. Then
the faithfulness of $\f$ implies that $F(ab)=ab$, and so
$\F_\8(T)$ is stable under multiplication, as desired.\cqd

\bigskip\n{\bf 7.3 Maximal ergodic inequalities}

\medskip\n
The following is the extension of Theorems \ref{max} and
\ref{smax} to the non tracial case.

\begin{thm}\label{haagmax}
{\rm i)} Let $T$ satisfy {\rm (7.I) - (7.IV)}. Let $(M_n)$ denote
the ergodic averages of $T$.  Then for any $1<p<\8$
 $$\big\|\,{\sup_n}^+M_n(x)\big\|_p\le C_p\,\|x\|_p\;,\quad x\in
 L_p(\M).$$
\indent{\rm ii)} If $T$ further satisfies {\rm (7.V)}, then
 $$\big\|\,{\sup_n}^+ T^n(x)\big\|_p\le C'_p\,\|x\|_p\;,\quad x\in
 L_p(\M).$$
Here $C_p$ and  $C_p'$ are respectively the  constants in
{\rm(\ref{max1})} and {\rm (\ref{smax1})}.
\end{thm}

\n{\bf Remark.} Compared with Theorems  \ref{max} and \ref{smax}
in the tracial case, the assumption in Theorem \ref{haagmax} is a
little bit stronger, namely, the positivity of $T$ in those
theorems is now reinforced to the complete positivity (7.II). It
is likely that this is not really needed.

\medskip

As in the tracial case, Theorem \ref{haagmax} immediately yields
the following two corollaries.

\begin{cor}\label{haagmaxc}
 Let $T$ satisfy {\rm (7.I) - (7.IV)} and $2<p<\8$. Then
 $$\big\|\big(M_n(x)\big)_{n\ge0}\big\|_{L_p(\M; \el_\8^c)}
 \le \sqrt{C_{p/2}}\,\|x\|_p\;, \quad\forall\;x\in L_p(\M).$$
If additionally $T$ has {\rm (7.V)}, then
 $$\big\|\big(T^n(x)\big)_{n\ge0}\big\|_{L_p(\M; \el_\8^c)}
 \le \sqrt{{C'}_{p/2}}\,\|x\|_p\;, \quad\forall\;x\in L_p(\M).$$
 \end{cor}

\medskip

\begin{cor}\label{haagmax-sg}
Let $(T_t)$ be a w*-continuous semigroup of maps on $\M$
satisfying {\rm (7.I) - (7.IV)}. Let
 $$M_t=\frac{1}{t}\,\int_0^t T_s\,ds,\quad t>0.$$
Then for any $1<p<\8$
 $$\big\|\,{\sup_t}^+M_t(x)\big\|_p\le C_p\,\|x\|_p\;,\quad x\in
 L_p(\M)$$
and for $p>2$
 $$\big\|\big(M_t(x)\big)_{t>0}\big\|_{L_p(\M; \el_\8^c(\real^+))}
 \le \sqrt{C_{p/2}}\,\|x\|_p\;, \quad\forall\;x\in L_p(\M).$$
If additionally each $T_t$ satisfies {\rm (7.V)}, then
 $$\big\|\,{\sup_t}^+T_t(x)\big\|_p\le C'_p\,\|x\|_p\;,\quad x\in
 L_p(\M)$$
and for $p>2$
 $$\big\|\big(T_t(x)\big)_{t>0}\big\|_{L_p(\M; \el_\8^c(\real^+))}
 \le \sqrt{{C'}_{p/2}}\,\|x\|_p\;, \quad\forall\;x\in L_p(\M).$$
\end{cor}

Although they are not stated here, all other inequalities in
sections \ref{maxergodic} and \ref{symergodic} continue to hold
for Haagerup noncommutative $L_p$-spaces. We omit the details. The
rest of this subsection is devoted to the proof of Theorem
\ref{haagmax}. It relies in a crucial way on  Haagerup's reduction
theorem \cite{haag-red}. We will need the precise form of
Haagerup's construction that we recall very briefly below.

\medskip

Let $G$ denote the discrete subgroup $\bigcup_{m \geq 1} 2^{-m}
\ent$ of $\real$. We consider the crossed product $\R = \M
\rtimes_{\s^{\f}} G$. Here the modular automorphism group
$\sigma^{\f}$ is also regarded as an automorphic representation of
$G$ on $\M$. As usual, $\M$ is viewed as a von Neumann subalgebra
of $\R$. Let $\wh\f$ denote the dual weight of $\f$. Since $G$ is
discrete, $\wh\f$ is a normal faithful state on $\R$ and its
restriction to $\M$ coincides with $\f$. Moreover, there is a
normal faithful conditional expectation $\Phi$ from $\R$ onto $\M$
such that
 $$\wh\f \circ \Phi= \wh\f \quad
 \mbox{and } \quad
 \s_{t}^{\wh\f} \circ \Phi = \Phi \circ
 \s^{\wh\f}_{t}, \quad t \in {\real}.$$
Then Haagerup's reduction theorem can be stated as follows.

\begin{thm}\label{haag} {\bf (Haagerup)} With the notations above, there is an
increasing sequence $(\R_m)_{m \ge 1}$ of von Neumann subalgebras
of $\R$ satisfying the following properties
 \begin{enumerate}[{\rm i)}]
 \item each $\R_m$ is equipped with a normal faithful tracial state $\t_m$;
 \item $\bigcup_{m \ge 1} \,\R_m$ is w*-dense in $\R$;
 \item there is a normal faithful conditional expectation
 $\Phi_{m}$ from $\R$ onto $\R_m$ such that
 $$\wh\f \circ \Phi_{m} = \wh\f \quad
 \mbox{and } \quad
 \s_{t}^{\wh\f} \circ \Phi_{m} = \Phi_{m} \circ
 \s^{\wh\f}_{t}, \quad t \in {\real}.$$
 \end{enumerate}
\end{thm}

We refer to \cite{haag-red} for the proof. \cite{jx-red}
reproduces Haagerup's proof and  presents several applications of
Theorem \ref{haag}.

In the situation above, $L_p(\M)$ and $L_p(\R_m)$ can be regarded
naturally and isometrically as  subspaces of $L_p(\R)$. Moreover,
the conditional expectation $\Phi$ (resp. $\Phi_m$) extends to a
positive contractive projection from $L_p(\R)$ onto $L_p(\M)$
(resp. $L_p(\R_m)$) (see \cite{jx-burk}; this is also a particular
case of Lemma \ref{extension1}). On the other hand, $\bigcup_{m
\ge 1} \,L_p(\R_m)$ is dense in $L_p(\R)$ for $p<\8$ and the
sequence $(\Phi_m)$ is increasing. Thus $(\R_m)$ gives rise to a
martingale structure on $\R$, and consequently, given $x\in
L_p(\R)$ with $1\le p<\8$,  $\Phi_m(x)$ converges to $x$ in
$L_p(\R)$ as $m\to\8$.

Let us also observe that by Remark \ref{vectorLp1bis} applied to
Haagerup spaces, $L_p(\M; \el_\8)$ and $L_p(\R_m; \el_\8)$ are
isometrically subspaces of $L_p(\R; \el_\8)$.

\medskip

For the proof of Theorem \ref{haagmax} we will further need the
following result from \cite{jx-red}.

\begin{lem}\label{crossed extension}
Let $T$ be as in Theorem \ref{haagmax}.
 \begin{enumerate}[{\rm i)}]
 \item Then $T$ has an extension $\wh T$ to $\R$ which satisfies {\rm (7.I) - (7.IV)}
relative to $(\R, \wh\f)$. Moreover, if $T$ verifies {\rm  (7.V)},
so does $\wh T$ relative to $\wh\f$.
 \item $\wh T(\R_m)\subset \R_m$ and $\t_m\circ \wh T\le \t_m$
for all $m\ge1$.
 \end{enumerate}
\end{lem}

Now we are ready to show Theorem \ref{haagmax}.

\smallskip

\n{\em Proof of Theorem \ref{haagmax}.} Fix $1<p<\8$ and $x\in
L_p(\M)$. We consider $x$ as an element in $L_p(\R)$ and then
apply the conditional expectation $\Phi_m$ to it:
 $\displaystyle x_m{\mathop =^{\rm def}}\Phi_m(x)\in L_p(\R_m)$.
Note that $\wh T\big|_{\R_m}$ satisfies the conditions (0.I) -
(0.III) relative to $\t_m$. So we can apply Theorem \ref{max} to
$\wh T$ on $\R_m$ and get
 $$\big\|\,{\sup_n}^+ M_n(\wh T)(x_m)\big\|_p\le C_p\,\|x\|_p\;,\quad
 \forall\; m\in\nat.$$
By the martingale convergence theorem recalled previously,
 $$\lim_{m\to\8}x_m=x\quad \mbox{in}\ L_p(\R).$$
Consequently,
 $$\lim_{m\to\8}\wh T\,^k(x_m)=x\quad \mbox{in}\ L_p(\R),
 \quad\forall\; k\ge0.$$
On the other hand, it is clear that the norm of $L_p(\R;\el_\8^n)$
is equivalent to that of $\el_\8^n\big(L_p(\R)\big)$ for each
fixed $n$. We then deduce that
 $$\lim_{m\to\8}\big\|\mathop{{\sup}^+}_{1\le k\le n} M_k(\wh T)(x_m)\big\|_p
 =\big\|\mathop{{\sup}^+}_{1\le k\le n} |M_k(\wh T)(x)|\big\|_p\;.$$
However, since $x\in L_p(\M)$,
 $$M_k(\wh T)(x)=M_k(T)(x).$$
Therefore, we deduce
 $$\big\|\mathop{{\sup}^+}_{1\le k\le n} M_k(T)(x)\big\|_p\le C_p\,\|x\|_p\;,
 \quad\forall\; n\in\nat.$$
Thus, by Proposition \ref{vectorLp0} and Remark
\ref{vectorLp1bis}, we have
 $$\big\|\,{\sup_n}^+ M_n(T)(x)\big\|_p\le C_p\,\|x\|_p\;.$$
This shows the first part of Theorem \ref{haagmax}. The second
part is proved similarly.\cqd

\bigskip\n{\bf 7.4 Individual ergodic theorems}

\medskip\n In this subsection we consider  individual ergodic theorems
in Haagerup's noncommutative $L_p$-spaces. As mentioned earlier,
the situation is more complicated than that in the tracial case.
One of the reasons is that the elements in $L_p(\M)$ are no longer
closed densely defined operators affiliated with $\M$ but
affiliated with a larger von Neumann algebra, namely the crossed
product $\M\rtimes_{\s^\f} {\mathbb R}$. We first need to
introduce an appropriate analogue of the almost everywhere
convergence for sequences in $L_p(\M)$. There are several such
generalizations. Here we adopt the almost sure convergence
introduced by Jajte \cite{ja2} (following ideas from \cite{dj2}).
In the $L_\8$-case, we continue to use Lance's almost uniform
convergence.

\begin{defi} {\rm i)} Let $x_n, x\in\M$. $x_n$ is said to
converge {\rm almost uniformly (a.u.} in short$)$ to $x$ if for
every $\e>0$ there is a projection $e\in\M$ such that
 $$\f(e^\perp)<\e\quad\mbox{and}\quad
 \lim_{n\to\8}\|(x_n-x)e\|_\8=0.$$
{\rm ii)} Let $x_n, x\in L_p(\M)$ with $p<\8$. The sequence
$(x_n)$ is said to converge {\rm almost surely (a.s.} in short$)$
to $x$ if for every $\e>0$ there is a projection $e\in\M$ and a
family $(a_{n,k})\subset \M$ such that
 $$\f(e^\perp)<\e\quad\mbox{and}\quad
 x_n-x=\sum_{k\ge1}a_{n,k}\,D^{\frac{1}{p}}\;,\quad
 \lim_{n\to\8}\big\|\sum_{k\ge1}(a_{n,k}\,e)\big\|_\8=0,$$
where the two series converge in norm in $L_p(\M)$ and $\M$,
respectively.

{\rm iii)} Similarly, we define {\rm bilateral almost uniform
(b.a.u.)} convergence and {\rm bilateral almost sure (b.a.s.)}
convergence. Note that for the latter we use the symmetric
injection of $\M$ into $L_p(\M):$ $a\mapsto
D^{\frac{1}{2p}}\,a\,D^{\frac{1}{2p}}$.
\end{defi}

The following non tracial analogue of Lemma \ref{con-c0} is
obtained in \cite{dj2}. For the sake of completeness we provide a
simplified proof.

\begin{lem}\label{haagcon-c0}
{\rm i)} If $(x_n)\in L_p(\M; c_0)$ with $1\le p<\8$, then $x_n$
converges b.a.s. to $0$.

{\rm ii)} If $2\le p<\8$ and $(x_n)\in L_p(\M; c_0^c)$, then $x_n$
converges a.s. to $0$.
\end{lem}

\pf  Suppose $(x_n)\in L_p(\M; c_0)$. Then there are $a, b\in
L_{2p}(\M)$ and $y_n\in \M$ such that
 $$x_n=ay_nb\quad\mbox{and}\quad \|a\|_{2p}<1,\ \|b\|_{2p}<1,\;
  \lim_{n\to\8} \|y_n\|_\8=0.$$
By the density of $D^{\frac{1}{2p}}\M$ in $L_{2p}(\M)$, there are
$a_k\in \M$ such that
 $$a=\sum_{k\ge1} D^{\frac{1}{2p}} a_k \quad\mbox{and}\quad
 \big\|D^{\frac{1}{2p}} a_k\big\|_{2p}< 2^{-k}\;.$$
Similarly,
 $$b=\sum_{k\ge1} b_k D^{\frac{1}{2p}}\quad\mbox{and}\quad
 \big\|b_k D^{\frac{1}{2p}}\big\|_{2p}< 2^{-k}\;.$$
Thus
 $$x_n=\sum_{j, k}
 D^{\frac{1}{2p}} a_jy_nb_k D^{\frac{1}{2p}}\,,\quad\mbox{convergence in}\;
 L_p(\M).$$
By the H\"older inequality,
 $$\f(a_ka_k^*)=\big\|D^{\frac{1}{2}} a_ka_k^* D^{\frac{1}{2}}\big\|_{1}
 \le\big\|D^{\frac{1}{2p}} a_ka_k^* D^{\frac{1}{2p}}\big\|_{p}
 <  2^{-2k}\;.$$
In the same way,
 $\f(b_k^*b_k)<  2^{-2k}$.
Now let $\e>0$. Then by \cite[Corollary 2.2.13]{ja1}, there is a
projection $e\in\M$ such that
 $$\f(e^\perp)<\e\quad\mbox{and}\quad
 \max\big\{\|ea_ka_k^*e\|_\8,\; \|eb_k^*b_ke\|_\8\big\}\le
 8\,\e^{-1}\,2^{-k}\;,\quad\forall\;k\ge1.$$
Therefore,
 $$\sum_{j,k\ge1}\|ea_jy_nb_ke\|_\8\le 8\,\e^{-1}\,
 \|y_n\|_\8\big[\sum_k 2^{-k/2}\big]^2\;;$$
whence the double series $\sum_{j,k}(ea_jy_nb_ke)$ converges
absolutely in $\M$ and
 $$\lim_{n\to\8}\sum_{j,k}ea_jy_nb_ke=0.$$
Hence $x_n\to  0$ b.a.s.. The second part is proved similarly.
\cqd

\begin{thm}\label{haagmaxc0}
{\rm i)} Let $T$ be a map on $\M$ satisfying {\rm(7.I) - (7.IV)}.
Then $\big(M_n(x)-F(x)\big)_n\in L_p(\M;c_0)$ for $1< p<\8$ and
$x\in L_p(\M)$. More generally, let $T_1,\,...,\, T_d$ be $d$ such
maps and let
 $$M_{n_1,\,...,\,n_d}=M_{n_d}(T_d)\,\cdots
 \, M_{n_1}(T_1)\,.$$
Let $F_k$ be the projection on the fixed point subspace of $T_k$.
Then
 $$\big(M_{n_1,\,...,\,n_d}(x)-F_d\,\cdots
 \, F_1(x)\big)_{n_1, \,...,\,n_d\ge 1}\in L_p(\M;
 c_0(\nat^d)),\quad \forall\; x\in L_p(\M),\ 1<p<\8$$
and
 $$\big(M_{n_1,\,...,\,n_d}(x)-F_d\,\cdots
 \, F_1(x)\big)_{n_1, \,...,\,n_d\ge 1}\in L_p(\M;
 c_0^c(\nat^d)),\quad \forall\; x\in L_p(\M),\ 2<p<\8.$$
 \indent{\rm ii)} If the $T_k$ further verify {\rm(7.V)}
and are positive operators on $L_2(\M)$, then in the statement
above the iterated ergodic averages $M_{n_1,\,...,\,n_d}$ can be
replaced by the iterated powers $T_d^{n_d}\,\cdots T_1^{n_1}$.
\end{thm}

\pf i) Let $1<p<\8$ and $x\in L_p(\M)$. By the discussion
following Remark \ref{extension3}, we can find $y_k\in\M$  and
$x_k=D^{1/2p}\big(y_k-T(y_k)\big)D^{1/2p}$  such that
 \be
 \lim_{k\to\8}\|x-F(x)-x_k\|_p=0\;.
 \ee
We have
 $$M_n(x_k)=\frac{1}{n+1}\,D^{\frac{1}{2p}}\big[y_k-T^{n+1}(y_k)\big]
 D^{\frac{1}{2p}}$$
and so $\big(M_n(x_k)\big)_n\in L_p(\M;c_0)$. Then as in the proof
of Theorem \ref{maxc0}, we deduce that $\big(M_n(x)-F(x)\big)_n\in
L_p(\M;c_0)$.

Now assume $2<p<\8$. Then by Remark \ref{extension3}, the $x_k$
above can be defined by
 $$x_k=(y_k-T(y_k))D^{\frac{1}{p}}\quad\mbox{with}\
 y_k\in\M_a\;.$$
Then $\big(M_n(x_k)\big)_n\in L_p(\M;c_0^c)$, and so
$\big(M_n(x)-F(x)\big)_n\in L_p(\M;c_0^c)$.

Iterating the arguments above as in the proof of Theorem
\ref{con-multi} and using the Haagerup space analogue of Corollary
\ref{max-multi}, we obtain the result for the multiple ergodic
averages.
\medskip

ii) In the case of one contraction, this part is proved in the
same way as Theorem \ref{scon}. The general case is dealt with by
iteration. We omit the details. \cqd

\begin{cor}\label{haagcon}
With the assumption and notations in Theorem \ref{haagmaxc0} {\rm
i)}, for any $1<p<\8$ and $x\in L_p(\M)$
 $$\lim_{n_1\to\8,\,...,\,n_d\to\8}\,M_{n_1,\,...,\,n_d}(x)
 =F_d\,\cdots\, F_1(x)\quad\mbox{b.a.s.;}$$
if $p>2$, the convergence above is a.s.. With the same assumption
as in Theorem \ref{haagmaxc0} {\rm ii)}, we have
 $$\lim_{n_1\to\8,\,...,\,n_d\to\8}\,T_d^{n_d}\,\cdots\, T_1^{n_1}(x)
 =F_d\,\cdots\, F_1(x)\quad\mbox{b.a.s.}$$
for $x\in L_p(\M)$ and $1<p<\8$. Again the convergence is {a.s.}
for $p>2$.
\end{cor}

\n{\bf Remarks.} i) Combining the preceding arguments with those
in the tracial case in section \ref{individual}, we easily show
that Theorem \ref{con-sg} continues to hold in the present setting
of Haagerup $L_p$-spaces with a.s. convergence in place of a.u.
convergence, as in Corollary \ref{haagcon} above.

ii) Using Goldstein's maximal weak type $(1,1)$  inequality
(\cite{gold}; see also \cite[Theorem 2.2.12]{ja2}), we can show
that the first part of Corollary \ref{haagcon} remains true for
$p=1$ and $d=1$ (i.e. for only one contraction).

iii) Jajte \cite{ja2} states a multiple individual ergodic theorem
in $L_2(\M)$ (Theorem 2.3.4 there), which corresponds to the first
part of Corollary \ref{haagcon} in the case of $p=2$. His proof
uses in an essential way his previous Theorem 2.2.4.  Based upon
an iteration using Goldstein's maximal weak type $(1,1)$
inequality, the proof of the latter theorem seems, however,  to
present a serious gap.

\medskip

Corollary \ref{haagcon} excludes the case $p=\8$, so does not
allow to recover all previous results by Lance \cite{la-erg},
K\"ummerer \cite{kum}, etc. This situation can be easily remedied.
This is done by virtue of the following simple lemma (see also
\cite{dj2}):

\begin{lem}\label{con-c0bis}
Let $1\le p<\8$ and $x_n\in\M$. Then
 \be
 (D^{\frac{1}{2p}}x_nD^{\frac{1}{2p}})_n\in L_p(\M;c_0) &\Longrightarrow&
 x_n\to 0\ {b.a.u.}\\
 (x_nD^{\frac{1}{2p}})_n\in L_p(\M;c_0^c)&\Longrightarrow&
 x_n\to 0\ {a.u.}
 \ee
\end{lem}

\pf Assume $(D^{\frac{1}{2p}}x_nD^{\frac{1}{2p}})_n\in
L_p(\M;c_0)$. Choose $a, b, y_n, a_k$ and $b_k$ exactly as in the
proof of Lemma \ref{haagcon-c0} (with
$D^{\frac{1}{2p}}x_nD^{\frac{1}{2p}}=ay_nb$). Next for each $n$
choose an integer $k_n$ such that
 $$\big\|D^{\frac{1}{2p}}x_nD^{\frac{1}{2p}} -\sum_{j, k=1}^{k_n}
 D^{\frac{1}{2p}}a_jy_nb_kD^{\frac{1}{2p}}\big\|_p<4^{-n}\;.$$
Set
 $$z_n=x_n\, - \,\sum_{j, k=1}^{k_n}a_jy_nb_k.$$
Then
 $$\big\|D^{\frac{1}{2}}z_nD^{\frac{1}{2}}\big\|_1\le
 \big\|D^{\frac{1}{2p}}z_nD^{\frac{1}{2p}}\big\|_p<4^{-n}\;.$$
Let $u_n$ and $v_n$ be respectively the real and imaginary part of
$z_n$. Then the inequality above holds with $u_n$ and $v_n$
instead of $z_n$. Now we apply \cite[Lemma 1.2]{haag-normalweight}
already quoted previously and reformulated in our setting as in
\cite{jx-red}. We then find $u_n', u_n''\in\M_+$ such that
$u_n=u_n' - u_n''$ and
 $$\big\|D^{\frac{1}{2}}u_nD^{\frac{1}{2}}\big\|_1
 =\big\|D^{\frac{1}{2}}u'_nD^{\frac{1}{2}}\big\|_1
 +\big\|D^{\frac{1}{2}}u''_nD^{\frac{1}{2}}\big\|_1
 =\f(u_n')+\f(u_n'').$$
Similarly, we have $v_n'$ and $v_n''$ for $v_n$. Thus
 $$\f(u_n')+\f(u_n'')<4^{-n}\;,\quad
 \f(v_n')+\f(v_n'')<4^{-n}\;.$$

Now given $\e>0$, applying \cite[Corollary 2.2.13]{ja1} to the
family
 $$\big\{a_na_n^*,\ b_n^*b_n,\ u_n',\ u_n'',\ v_n',\ v_n''\ :\
n\in\nat\big\},$$
we get a projection $e\in\M$ such that
$\f(e^\perp)<\e$ and
 $$\max\big\{\|ea_na_n^*e\|_\8, \, \|eb_n^*b_ne\|_\8,\,
 \|eu_n'e\|_\8, \, \|eu_n''e\|_\8, \, \|ev_n'e\|_\8,
 \, \|ev_n''e\|_\8 \big\}
 <16\,\e^{-1}\, 2^{-n}$$
for all $n\in\nat$. Therefore,
 \be
 \|ex_ne\|_\8
 &\le& \|ez_ne\|_\8 + \big\|\sum_{j, k=1}^{k_n}ea_jy_nb_ke\big\|_\8\\
 &\le& \|e(u_n'- u_n'')e\|_\8 +\|e(v_n'- v_n'')e\|_\8 +
 \sum_{j, k=1}^{k_n} \|ea_jy_nb_ke\|_\8\\
 &\le& 64\,\e^{-1}\, 2^{-n} + 16\,\e^{-1}\,
 \|y_n\|_\8\big[\sum_{k\ge1} 2^{-k/2}\big]^2\ \to0\ \mbox{as}\
 n\to\8.
 \ee
Thus $x_n\to0$ b.a.u.. The proof of the second part on the a.u.
convergence is similar and even easier (without appealing to
Haagerup's Lemma). Thus we omit the details. \cqd

\medskip

The first part of the following is well-known (cf., e.g.
\cite{ja1}).

\begin{cor}\label{haagconL8}
let $T_1,\,...,\, T_d$  satisfy {\rm(7.I) - (7.IV)} and let
 $$M_{n_1,\,...,\,n_d}=M_{n_d}(T_d)\,\cdots
 \, M_{n_1}(T_1)\,.$$
Then for any $x\in \M$
 $$\lim_{n_1\to\8,\,...,\,n_d\to\8}\,M_{n_1,\,...,\,n_d}(x)
 =F_d\,\cdots\, F_1(x)\quad\mbox{a.u.}$$
If $T_1,\,...,\, T_d$ additionally have {\rm(7.V)}, then
 $$\lim_{n_1\to\8,\,...,\,n_d\to\8}\,T_d^{n_d}\,\cdots\, T_1^{n_1}(x)
 =F_d\,\cdots\, F_1(x)\quad\mbox{a.u.}$$
\end{cor}

\pf This immediately follows from Theorem \ref{haagmaxc0} and
Lemma \ref{con-c0bis}.\cqd

\medskip
\n{\bf Remark.} In the case of $d=1$,  the first part of Corollary
\ref{haagconL8} permits to recover Lance's theorem. However,
compared with K\"ummerer's theorem, our hypothesis is stronger for
K\"ummerer assumed only that $T$ is a positive contraction
verifying (7.III). We do not know whether all ergodic theorems in
this section hold for such contractions or not. In particular, is
Theorem \ref{haagmax} true for a positive contraction $T$
satisfying (7.III) (and (7.V))?

\medskip
\n{\bf Remark.} As in the tracial case, all the preceding
individual ergodic theorems admit semigroup analogues.

\section{Examples}

We will give some natural examples to which the results in the
previous sections can be applied.
\bigskip

\n{\bf 8.1 Modular groups}

\medskip\n
The very first examples are modular automorphism groups. Let $\f$
be a normal faithful state on a von Neumann algebra $\M$. Let
$\s_t^\f$ be the modular group of $\f$. Then $T_t=\s_t^\f$
satisfies the properties (7.I) - (7.IV).  On the other hand, (7.V)
is equivalent to $\f(\s_t^\f(y)x)=\f(y\s_{-t}^\f(x))$ for all $x,
y\in\M$ and $t\in\real$. Thus applying Corollary \ref{haagmax-sg},
we get that for $1<p<\8$
 $$\big\|\,{\sup_t}^+\frac{1}{t}\int_0^t \s_s^\f(x)ds\big\|_p
 \le C_p\, \|x\|_p\,,\quad\forall\; x\in L_p(\M).$$
Note that the fixed point subspace $\F_\8$ of $(\s_t^\f)$
coincides with the centralizer ${\M}_\f$ of $\f$. Consequently,
$\F_p$ coincides with $L_p({\M}_\f)$, considered as a subspace of
$L_p(\M)$. Thus applying the results in subsection 7.4, we deduce
that the ergodic averages
 $\dis\frac{1}{t}\int_0^t \s_s^\f(x)ds$
converge b.a.u. to $x$ (resp. $F(x)$) as $t\to 0$ (resp. $t\to\8$)
for all $x\in L_p(\M)$ and $1\le p\le\8$. Moreover, the
convergence is a.u. in the case of $p\ge2$. Let us consider a
state $\f(x)=\lambda x_{11}+\mu x_{22}$ on the matrix algebra
$\mathbb M_2$ of $2\times 2$ matrices, where $0<\l\neq\mu<1$. Then
we see that $\sigma_t^{\f}(e_{12})=e^{it(\lambda-\mu)}\,e_{12}$ is
not convergent for $t\to \infty$. At least in this case it is
obvious that the symmetry condition (7.V) is really necessary.

\bigskip

\n{\bf 8.2 Semi-noncommutative case}

\medskip\n
Let $(\O, \F, \mu)$ be a $\s$-finite measure space. Let $\N$ be a
von Neumann algebra equipped with a semifinite normal faithful
trace $\nu$. Let $(\M, \t)=(L_\8(\O),\mu)\bar \ot(\N, \nu)$ be the
von Neumann algebra tensor product. (Note that we consider $\mu$
as a trace on $L_p(\O)$ via integration.) Given $p<\8$ the
corresponding noncommutative $L_p(\M)$  is just $L_p(\O;
L_p(\N))$, the usual $L_p$-space of strongly measurable
$p$-integrable functions on $\O$ with values in $L_p(\N)$. Now let
$(S_t)$ be a semigroup on $L_p(\O)$ satisfying the conditions
(0.I) - (0.III) (with $\M=L_\8(\O)$ there). Then $T_t=I\ot S_t$ is
a semigroup on $L_p(\M)$ verifying the same conditions. Moreover,
if $S_t$ is symmetric, so is $T_t$. Thus we can transfer all
classical semigroups to this semi-noncommutative setting and
obtain the corresponding ergodic theorems. In particular, applying
this procedure to the usual Poisson semigroup $(P_t)$ on the unit
circle  $\T$ or on $\real^n$, by Corollary \ref{max-poisson}, we
get
 $$ \big\|\,{\sup_t}^+I\ot P_t(x)\big\|_p\le C_p\,\|x\|_p,
 \quad x\in L_p(\real^n; L_p(\N)),\ 1<p<\8.$$
For $p=1$ we also have a weak type inequality (see Remark
\ref{yeadon-sg}). These results were also proved by Mei \cite{mei}
using a different method. Moreover, he obtained the non-tangential
analogue (for the upper half plane) of the inequality above. Note
that in this discussion, the usual Poisson semigroup on $\real^n$
can be replaced by the Poisson semigroup subordinated to the
Ornstein-Uhlenbeck semigroup on $\real^n$.

\medskip
The situation above  readily extends to the non tracial case.
Assume that $\mu$ is a probability measure and $\N$ a von Neumann
algebra equipped with a normal faithful state $\psi$. Then the
tensor product $\M$ is equipped with the tensor state
$\f=\mu\ot\psi$. This allows to apply the ergodic results in
section \ref{haagerup} to this semi-noncommutative setting.

\bigskip


 \n{\bf 8.3  Schur multipliers}

\medskip\n
Let $\M=B(\el_2)$. Then the associated noncommutative $L_p$-spaces
are the Schatten classes $S_p$. The elements in $S_p$ are
represented as infinite matrices. Let $\phi$ be a function on
$\nat\times \nat$. Recall that $\phi$ is a Schur multiplier on
$S_p$  if the map $M_\phi: x\mapsto (\phi_{jk}x_{jk})$, defined
for finite matrices $x$, extends to a bounded map  on $S_p$ (which
is still denoted by $M_\phi$).

Let us consider a function  $f: \nat\to H$, where $H$ is a real
Hilbert space, and the associated kernel
 $$K(j,k)=\|f(j)-f(k)\|,\quad j,k\in\nat.$$
We are interested in the semigroups $(T_t)$ and $(P_t)$ of Schur
multipliers, which are determined by
 $$T_t(e_{jk})=e^{-t K(j,k)^2}e_{jk}\quad\mbox{and}\quad
 P_t(e_{jk})=e^{-tK(j,k)}e_{jk}\,,$$
where the $e_{jk}$'s stand for the canonical matrix units of
$B(\el_2)$. It is well-known that these are completely positive
contractive semigroups on $B(\el_2)$. Indeed, let $\mu$ be a
Gaussian measure on $H$, i.e. a probability space $(\O, \mu)$
together with a measurable function $w: \O\to H$ such that
 $$\exp\big(- \|h\|^2 \big)=\int_\O\exp\big(i\la h,\,
 w(\o)\ra\big)\,d\mu(\o)\,,\quad h\in H.$$
Given $\o\in\O, t>0$  let $D_t(\o)$ be the diagonal matrix with
diagonal entries  $\exp\big(i\sqrt t\,\la f(j),\,w(\o)\ra\big)$,
$j\in\nat$. Then it is easy to see that
 \beq\label{schur}
 T_t(x)= \int_\O D_t(\o)xD_t(\o)^*\,d\mu(\o)\,, \quad x\in B(\el_2).
 \eeq
Since $D_t(\o)$ is unitary,  this formula shows that $T_t$ is a
completely positive contraction on $B(\el_2)$. In fact,
(\ref{schur}) is the Stinespring representation of $T_t$. The
semigroup $(T_t)$ satisfies all properties (0.I) - (0.IV) with
$\M=B(\el_2)$ and $\t$ being the usual trace on $B(\el_2)$. Since
$(P_t)$ is the Poisson semigroup subordinated to $(T_t)$ via
(\ref{subordination}), $(P_t)$ has the same properties. Thus these
semigroups extend to symmetric positive contractive semigroups on
$S_p$ for $1\le p<\8$.

Thus we have the maximal inequalities in Theorem \ref{max} and
Theorem \ref{smax}  for $(T_t)$ as well as (\ref{max-poisson1})
for $(P_t)$. Note that in this situation the a.u. convergence
reduces to the uniform convergence in $B(\el_2)$.

\bigskip \n{\bf 8.4 Hamiltonians}

\medskip\n
In this subsection, $\M$ is semifinite and  equipped with a normal
faithful semifinite trace $\t$. Let $L\in L_0(\M)$ be selfadjoint.
We consider the Hamiltonian semigroup given by the generator ${\rm
ad}\,L$:
 $${\rm ad}\,L (x)=Lx-xL,\quad x\in\M.$$
Note that
 $$({\rm ad}\,L)^2 (x)=L^2x+xL^2 -2LxL.$$
Set
 $$T_t=e^{-t({\rm ad}\,L)^2}\quad\mbox{and}\quad
 P_t=e^{-t\,|{\rm ad}\,L |}\,.$$
It is again well-known that these are completely positive
contractive semigroups on $\M$ (see \cite[Example
30.1]{parth-bk}). Since $(P_t)$ is the Poisson semigroup
subordinated to $(T_t)$, it suffices to show this for $(T_t)$.  In
fact, $(T_t)$ admits a Stinespring representation similar to
(\ref{schur}):
 \beq\label{hamiltonian}
 T_t(x)={\mathbb E}\big[e^{i \sqrt t\, gL} xe^{-i \sqrt t\,
 gL}\big]\,,\quad x\in\M,
 \eeq
where $g$ is a Gaussian variable with mean zero and variance
$\sqrt2$ and ${\mathbb E}$ denotes the expectation with respect to
$g$. To check this, let us first write the spectral resolution of
$L$:
 $$L=\int_{-\8}^\8 \l\, de_\l.$$
Let $R>0$ and $e$ be the spectral projection of $L$ corresponding
to the interval $[-R, R]$. Consider $x\in\M$ such that $x=exe$.
Then
 \beq\label{hamiltonian0}
 \|L^jxL^k\|\le R^{j+k} \|x\|, \quad \forall\; j, k\ge 0.
 \eeq
A simple induction shows
 $$({\rm ad}\,L)^n(x)=(-1)^n\sum_{k=0}^n (-1)^k C_n^k\, L^k x
 L^{n-k}\,.$$
Now consider the formal power series representation
 \be
 {\mathbb E}\big[e^{i \sqrt t\, gL} xe^{-i \sqrt t\, gL}\big]
 &=&{\mathbb E}\big[\sum_{j, k=0}^\8
 \frac{(i\sqrt t)^j(-i\sqrt t)^k}{j!\,k!}\,g^{j+k} L^jxL^k\big]\\
 &=&\sum_{n=0}^\8 \frac{(-1)^n t^n}{n!}\,
 \sum_{j+k=2n}\frac{(2n)!}{j!\,k!} (-1)^k L^jxL^k\\
 &=& \sum_{n=0}^\8 \frac{(-1)^n t^n}{n!} ({\rm ad}\,L)^{2n}(x)
 =e^{-t({\rm ad}\,L)^2}(x).
 \ee
Note that the series above are absolutely convergent due to
(\ref{hamiltonian0}). Thus (\ref{hamiltonian}) is proved for all
$x\in\M$ such that $x=exe$. However,  the left hand side of
(\ref{hamiltonian}) defines a normal contraction on $\M$ for
$\exp(i \sqrt t\, g(\o)L)$ is a unitary in $\M$ for every $\o$. On
the other hand, $\lim_{R\to\8} \un_{[-R,\, R]}(L)=1$ weakly in
$\M$. By the w*-continuity, we see that (\ref{hamiltonian}) is
true for all $x\in\M$. (\ref{hamiltonian}) also shows that $T_t$
preserves the trace $\t$. On the other hand, since $({\rm
ad}\,L)^2$ is positive on $L_2(\M)$, $T_t$ is symmetric. Thus the
semigroup $(T_t)$ verifies (0.I) - (0.IV).

\medskip
\n{\bf Remark.} Let us consider a particular case where
$\M=B(\el_2)$ and $L$ is a real diagonal matrix with diagonal
entries $(\l_0, \l_1,\,\cdots)$. Then
 $$|{\rm ad}\, L|(x)= \big(|\l_j-\l_k|x_{jk}\big)_{j,k}\,.$$
Thus $|{\rm ad}\, L|$ becomes a Schur multiplier and so $(T_t)$
reduces to the semigroup already considered in the previous
example with $H=\real$ and $f(j)=\l_j$.

\bigskip

\n{\bf 8.5 Free product}

\medskip\n
Let $(\M_i, \f_i)_{i\in I}$ be a family of von Neumann
algebras, each equipped with a normal faithful state $\f_i$. Let
 $$(\M, \f)=\ast_{i\in I} \big(\M_i, \f_i\big)$$
be the von Neumann algebra reduced free product (cf.
\cite{voi-sym} and \cite{VDN}). Recall that $\f$ is a normal
faithful state on $\M$. If all $\f_i$ are tracial, so is $\f$. Now
for every $i\in I$ let be given a w*-continuous semigroup
$(T_t^i)_{t\ge 0}$ on $\M_i$ satisfying the following conditions:
 \begin{enumerate}[{\rm i)}]
 \item $T_t^i$ is unital;
 \item $\f\circ T_t^i=\f_i$;
 \item $T_t^i$ is completely positive.
\end{enumerate}
As usual, we always assume $T_0^i={\rm id}_{\M_i}$. Then by
\cite{blandyk} (see also \cite{choda-red}) it follows that for
each $t$ the family $\{T_t^i\}_{i\in I}$ defines a completely
positive unital map $T_t$ on $\M$, preserving the state $\f$.
$T_t$ is uniquely determined by its action on the monomials:
 $$T_t(x_1\,\cdots\, x_n)=T_t^{i_1}(x_1)\,\cdots\, T_t^{i_n}(x_n)$$
for any $x_1,\,...,\, x_n$ with $x_k\in \M_{i_k}^\circ$ and
$i_1\not=i_2\not=\cdots \not=i_n$, where $\M_i^\circ=\{x\in\M_i:
\f_i(x)=0\}$. $T_t$ is called  the free product of  the family
$\{T_t^i\}_{i\in I}$ and denoted by $T_t=\ast_{i\in I}\, T_t^i$.
Then it is easy to see that $(T_t)$ is a w*- continuous semigroup
on $\M$. Thus this semigroup satisfies the conditions (7.I) -
(7.III). By Lemma \ref{extension1}, $(T_t^i)$ and $(T_t)$ extend
to norm continuous semigroups respectively on $L_p(\M_i)$ and
$L_p(\M)$ for all $1\le p<\8$.

Recall that the modular group $\s_t^\f$ is the free product of the
modular groups $\s_t^{\f_i}$, $i\in I$ (cf. \cite{dyk-crelle}).
Thus if each $T_t^i$ satisfies (7.IV), so does $T_t$. On the other
hand, it is clear that the property (I.V) is also stable under
free product.

\medskip

Let us consider one special case. Note that $\M_i= \comp 1_{\M_i}
\oplus \M_i^\circ$.  Let $T_t^i: \M_i\to\M_i$ be defined by
 $$T_t^i\big|_{\comp 1_{\M_i}}={\rm id}_{\comp 1_{\M_i}}\quad
 \mbox{and}\quad
 T_t^i\big|_{\M_i^\circ}=e^{-t}\,{\rm id}_{\M_i^\circ}\,,\quad t\ge0.$$
Then it is easy to check that $(T_t^i)$ verifies the conditions i)
- iii) above; moreover, $T_t^i$ is symmetric relative to $\f_i$.
The corresponding free product semigroup $(T_t)$ is uniquely
determined by
 $$T_t(x_1\,\cdots\, x_n)=e^{-nt}x_1\,\cdots\, x_n$$
for any $x_1,\,...,\, x_n$ with $x_k\in \M_{i_k}^\circ$ and
$i_1\not=i_2\not=\cdots \not=i_n$ with $n\in\nat$. This is the
free analogue of the classical Poisson semigroup on the unit
circle. It plays an important role in \cite{rixu}.

The fixed point subspace of $(T_t)$ above is simply $\comp
1_{\M}$. Let us briefly discuss the pointwise convergence in this
case. Every element $x\in \M$ admits the following formal
development
 $$x=\f(x) +\sum_{n\ge1}\sum_{i_1\not=\,\cdots\n\not=i_n}
 x_1\,\cdots\, x_n\,,$$
where $x_k\in \M_{i_k}^\circ$. Then
 $$T_t(x)=\f(x) +\sum_{n\ge1}e^{-nt}\,\sum_{i_1\not=\,\cdots\n\not=i_n}
 x_1\,\cdots\, x_n\,.$$
Thus by the results in subsection~7.4,
 $$\lim_{t\to0} T_t(x)=x\quad\mbox{and}\quad
 \lim_{t\to\8} T_t(x)=\f(x)\quad\mbox{a.u.}$$
A similar result also holds  for $x\in L_p(\M)$ with $1<p<\8$.

\bigskip \n{\bf 8.6 Group von Neumann algebras}

\medskip\n
Let $G$ be a discrete group. Let $VN(G)$ denote the group von
Neumann algebra of $G$. Recall that $VN(G)$ is a von Neumann
algebra on $\el_2(G)$ generated by the left regular representation
$\l$. Let $\t_G$ be the canonical faithful tracial state on
$VN(G)$, i.e. $\t_G$ is the vector state given by the unit basis
vector $\d_e$, where $e$ is the identity of $G$ and where
$\{\d_g\}_{g\in G}$ denotes the canonical basis of $\el_2(G)$.

Now we assume that $G$ is equipped with a length function, denoted
by $|\cdot|$. More precisely, $|\cdot|$ is a positive function on
$G$ satisfying the following conditions:
\begin{enumerate}[i)]
 \item $|e|=0$;
 \item $|g^{-1}|=|g|$ for any $g\in G$;
 \item if $d(f, g)=\frac{1}{2}(|f|+|g|-|fg^{-1}|)$, then for all
 $f, g, h\in G$
  $$d(f, g)\ge\min\{d(f, h),\; d(h, g)\}.$$
 \end{enumerate}
Bo{\.z}ejko \cite{boz-length}  proved that $g\mapsto e^{-t|g|}$ is
a positive definite function on $G$ (see also \cite{boz-lect}).
Thus the associated Herz-Schur multiplier $T_t$ is a normal
completely positive unital map on $VN(G)$. More precisely, $T_t$
is given on polynomials by
 $$T_t\big(\sum_g a_g\,\l(g)\big)=\sum_g e^{-t|g|}\,a_g\,\l(g).$$
Moreover, $T_t$ preserves the trace $\t_G$. Thus by Lemma
\ref{extension}, $(T_t)$ extends to a semigroup on $L_p(VN(G))$
for all $1\le p<\8$.  Note that if $G=\ent$, then $VN(G)=L_\8(\T)$
and $T_t$ becomes the usual Poisson semigroup on $\T$.

More generally, it is proved in \cite{boz-length} that for any
$0<\a<1$ the function $g\mapsto e^{-t|g|^\a}$ is  positive
definite on $G$. It follows that
 $$P_t \big(\sum_g a_g\,\l(g)\big)=\sum_g e^{-t|g|^\a}\,a_g\,\l(g)$$
defines a completely positive unital trace preserving semigroup on
$VN(G)$. This last statement also follows from the previous for
$(P_t)$ is subordinated to $(T_t)$ by (\ref{subordinationbis}).

\medskip

Now let us specify the situation above to free groups. Let $G$ be
a free group, say, $G=\FF_n$, a free group on $n$ generators
$\{g_1, \,...,\, g_n\}$ ($n$ can be infinite). Let $|\cdot |$ be
the length function with respect to $\{g_1, \,...,\, g_n\}$. Then
the fact that $e^{-t|\cdot |}$ is a positive definite function on
$\FF_n$ goes back to Haagerup \cite{haag-free}.  Note that this is
also a special case of the free product in the previous example.
Indeed, writing $\FF_n$ as the reduced free product of $n$ copies
of $\ent$, we have
 $$\big(VN(\FF_n),\;\t_{\FF_n}\big) \ast_{1\le k\le n}\big(L_\8(\T),\; \t_{\ent}\big).$$
Then the semigroup on $\FF_n$ appears as the free product of $n$
copies of the usual Poisson semigroup on $\T$. Applying our
ergodic theorems to this case, we get Theorem \ref{free}.

More generally, let $\{G_i\}_{i\in I}$ be a family of discrete
groups, each equipped with a length function. Let $T_t^i$ be the
associated semigroup on $G_i$ defined previously. Let
$G=\ast_{i\in I}\,G_i$ be the reduced free product. Then by
\cite{boz-length} (or the previous example), the free product
$T_t=\ast_{i\in I} T_t^i$ yields a symmetric completely positive
contractive semigroup on $G$.

\bigskip\n{\bf 8.7 $q$-Ornstein-Uhlenbeck semigroups}

\medskip\n
Let $H_\real$ be a real Hilbert space and $H_\comp$ its
complexification. For $-1\le q\le 1$ let $\F_q(H_\comp)$ be the
$q$-Fock space based on $H_\comp$ constructed by Bo\.zejko and
Speicher (see \cite{bos-example} and \cite{bos-interp}). Note that
$\F_1(H_\comp)$, $\F_{-1}(H_\comp)$ and $\F_0(H_\comp)$ are
respectively the symmetric, anti-symmetric and full (=free) Fock
spaces. Given a vector $h\in H_\comp$, let $c(h)$ denote the
associated (left) creation operator on $\F_q(H_\comp)$. $c(h)$ is
a bounded operator for $q<1$ and a closed densely defined operator
for $q=1$. Its adjoint $c(h)^*$ is the annihilation operator
associated to $h$ and denoted by $a(h)$. Let
 $$g_q(h)=c(h)+a(h), \quad h\in H_\real.$$
$g_q(h)$ is a so-called $q$-Gaussian variable. Note that $g_1(h)$
is a usual Gaussian variable, $g_0(h)$ a semi-circular variable in
Voiculescu's sense (cf. \cite{voi-sym} and \cite{VDN}), and
finally $g_{-1}(h)$ corresponds to a Fermion. The $q$-von Neumann
algebra $\Gamma_q(H_\real)$ is the von Neumann algebra on
$\F_q(H_\comp)$ generated by all $q$-Gaussians, namely,
 $$\Gamma_q(H_\real)=\{g_q(h)\;:\; h\in H_\real\}''\subset B(\F_q(H_\comp)).$$

Let $\O$ be the vacuum vector in $\F_q(H_\comp)$ and $\t_q$ the
associated vector state. Then $\t_q$ is faithful and tracial.
Hence $\Gamma_q(H_\real)$ is a type II$_1$ von Neumann algebra for
$q<1$. ($\Gamma_1(H_\real)$ is commutative.) Moreover, it is a non
injective factor if $-1<q<1$ and $\dim H\ge 2$. We refer to
\cite{boks}, \cite{nou-inj} and  \cite{ricard-fac} for more
information.

\medskip

Now let $S$ be a contraction on $H_\real$. Then $S$ extends to a
contraction on $H_\comp$.  The second quantization $\Gamma(S)$ is
a normal completely positive unital trace preserving map on
$\Gamma_q(H_\real)$. To give the definition of $\Gamma(S)$, we
recall the Wick product. Since $\O$ is separating for
$\Gamma_q(H_\real)$, the map $x\in \Gamma_q(H_\real)\mapsto x(\O)$
is injective. Its image is a dense subspace of $\F_q(H_\comp)$
(for $\O$ is cyclic). It is easy to see that all elementary
tensors belong to this image. The inverse map (defined on the
image) is called the Wick product, denoted by $W$. Thus if $\xi$
is a linear combination of elementary tensors, $W(\xi)$ is the
unique operator in $\Gamma_q(H_\real)$ such that $W(\xi)\O=\xi$.
Note that the collection of all such $W(\xi)$'s  forms a w*-dense
$\ast$-subalgebra of $\Gamma_q(H_\real)$. Then $\Gamma(S)$ is
uniquely determined by
 $$\Gamma(S)\big(W(h_1\ot\,\cdots\,\ot h_n)\big) W(Sh_1\,\ot\cdots \,\ot Sh_n),\quad h_1, \,...,\, h_n\in H_\comp.$$

Applying this construction to $S=e^{-t}\,{\rm id}_{H_\real}$ for
$t\ge0$, we get a normal completely positive unital trace
preserving map $T_t=\Gamma(e^{-t}\,{\rm id}_{H_\real})$. The
action of $T_t$ on the Wick products is given by
 $$T_t\big(W(h_1\ot\,\cdots\,\ot h_n)\big)
 =e^{-nt}\,W(h_1\ot\,\cdots\,\ot h_n).$$
Then $(T_t)$ is a semigroup on $\Gamma_q(H_\real)$ satisfying all
conditions (0.I) - (0.IV). This is the $q$-Ornstein-Uhlenbeck
semigroup associated with $H_\real$. The negative of its
infinitesimal generator is the so-called number operator. The case
$q=1$ and $q=-1$ corresponds respectively to the classical and
Fermionic Ornstein-Uhlenbeck semigroup. $(T_t)_t$ in these two
special cases have been extensively studied. See \cite{boz-ultra},
\cite{carlenlieb} and \cite{biane-hyper} for related results.

\medskip
The preceding discussion also applies to the quasi free case. Then
the corresponding von Neumann algebras are of type III. See
\cite{shlya-quasifree} for the case of $q=0$  and \cite{hiai} for
the general case. In particular, for $q=-1$, we have the classical
Araki-Woods factors. In this case, the resulting semigroup is the
extension of the previous Fermionic Ornstein-Uhlenbeck semigroup
to the type III setting.


\begin{thebibliography}{VDN}


\bibitem[BeL]{bl}
Bergh, J. and L{\"o}fstr{\"o}m, J.~
\newblock {\em Interpolation spaces.}
\newblock Springer-Verlag, Berlin, 1976.

\bibitem[Bi]{biane-hyper}
Biane, Ph.~
\newblock Free hypercontractivity.
\newblock {\em Comm. Math. Phys.}, 184:457--474, 1997.

\bibitem[BlD]{blandyk}
Blanchard, E. and Dykema, K.~
\newblock Embeddings of reduced free products of operator algebras.
\newblock {\em Pacific J. Math.}, 199:1--19, 2001.

\bibitem[Bo1]{boz-length}
Bo{\.z}ejko, M.~
\newblock Positive-definite kernels, length functions on groups and a
  noncommutative von {N}eumann inequality.
\newblock {\em Studia Math.}, 95:107--118, 1989.

\bibitem[Bo2]{boz-lect}
Bo{\.z}ejko, M.~
\newblock Positive and negative definite kernels on discrete groups.
\newblock Lectures at Heidelberg University, 1987.


\bibitem[Bo3]{boz-ultra}
Bo{\.z}ejko, M.~
\newblock Ultracontractivity and strong {S}obolev inequality for
  {$q$}-{O}rnstein-{U}hlenbeck semigroup {$(-1<q<1)$}.
\newblock {\em Infin. Dimens. Anal. Quantum Probab. Relat. Top.},
  2:203--220, 1999.

\bibitem[BKS]{boks}
Bo{\.z}ejko, M.,  K{\"u}mmerer, B. and Speicher, R.~
\newblock {$q$}-{G}aussian processes: non-commutative and classical aspects.
\newblock {\em Comm. Math. Phys.}, 185:129--154, 1997.


\bibitem[BS1]{bos-example}
Bo{\.z}ejko, M. and Speicher, R.~
\newblock An example of a generalized {B}rownian motion.
\newblock {\em Comm. Math. Phys.}, 137:519--531, 1991.

\bibitem[BS2]{bos-interp}
Bo{\.z}ejko, M. and Speicher, R.~
\newblock Interpolations between bosonic and fermionic relations given by
  generalized {B}rownian motions.
\newblock {\em Math. Z.}, 222:135--159, 1996.

\bibitem[CaL]{carlenlieb}
Carlen, E. and Lieb, E.~
\newblock Optimal hypercontractivity for {F}ermi fields and related
noncommutative integration inequalities.
\newblock {\em Comm. Math. Phys.}, 155:27--46, 1993.

\bibitem[Ch]{choda-red}
Choda, M.~
\newblock Reduced free products of completely positive maps and entropy
for free product of automorphisms.
\newblock {\em Publ. Res. Inst. Math. Sci.}, 32:371--382, 1996.

\bibitem[ChE]{choi-Ef}
Choi, M-D and Effros, Ed.~
\newblock Injectivity and operator spaces.
\newblock {\em J. Funct. Anal.}, 24:156--209, 1977.

\bibitem[CoN]{cdn}
Conze, J.-P.  and Dang-Ngoc, N.~
\newblock Ergodic theorems for noncommutative dynamical systems.
\newblock {\em Invent. Math.}, 46:1--15, 1978.

\bibitem[Cu]{cu}
Cuculescu, I.~
\newblock Martingales on von {N}eumann algebras.
\newblock {\em J. Multivariate Anal.}, 1:17--27, 1971.

\bibitem[Da]{dngoc}
Dang-Ngoc, N.~
\newblock Pointwise convergence of martingales in von {N}eumann algebras.
\newblock {\em Israel J. Math.}, 34:273--280, 1979.

\bibitem[DJ1]{dj}
Defant, A.~ and Junge, M.~
\newblock Maximal theorems of {M}enchoff-{R}ademacher type in non-commutative
  {$L\sb q$}-spaces.
\newblock {\em J. Funct. Anal.}, 206:322--355, 2004.

\bibitem[DJ2]{dj2}
Defant, A.~ and Junge, M.~
\newblock Classical summation methods in noncommutative
 probability.
\newblock In preparation.

\bibitem[DS]{dunford}
Dunford, N. and Schwartz, J.T.~
\newblock {\em Linear {O}perators. {I}. {G}eneral {T}heory}.
\newblock Applied Mathematics, Vol. 7. Interscience Publishers, Inc., New York,
1958.

\bibitem[Dy]{dyk-crelle}
Dykema, K.~
\newblock Factoriality and {C}onnes' invariant {$T({\M})$} for free
products of von {N}eumann algebras.
\newblock {\em J. Reine Angew. Math.}, 450:159--180, 1994.

\bibitem[FK]{fk}
Fack, Th. and Kosaki, H.~
\newblock Generalized {$s$}-numbers of {$\tau$}-measurable operators.
\newblock {\em Pacific J. Math.}, 123:269--300, 1986.

\bibitem[Go]{gold}
Goldstein, M.~S.~
\newblock Theorems on almost everywhere convergence in von {N}eumann algebras.
\newblock {\em J. Operator Theory}, 6:233--311, 1981.

\bibitem[GoG]{gog}
Goldstein, M.~S. and Grabarnik, G.~Y.~
\newblock Almost sure convergence theorems in von Neumann algebras.
\newblock {\em Israel J. Math.}, 76:161--182, 1991.

\bibitem[H1]{haag-normalweight}
Haagerup, U.~
\newblock Normal weights on {$W\sp{\ast} $}-algebras.
\newblock {\em J. Functional Analysis}, 19:302--317, 1975.

\bibitem[H2]{haag-Lp}
Haagerup, U.~
\newblock {$L\sp{p}$}-spaces associated with an arbitrary von {N}eumann
  algebra.
\newblock In {\em Alg\`ebres d'op\'erateurs et leurs applications en physique
  math\'ematique (Proc. Colloq., Marseille, 1977)}, volume 274 of {\em Colloq.
  Internat. CNRS},  175--184. CNRS, Paris, 1979.

\bibitem[H3]{haag-red}
Haagerup, U.~
\newblock Non-commutative integration theory.
\newblock See also Haagerup's Lecture given at the Symposium in Pure
  Mathematics of the Amer. Math. Soc., Queens University, Kingston, Ontario, ~
  1980.

\bibitem[H4]{haag-free}
Haagerup, U.~
\newblock An example of a nonnuclear {$C\sp{\ast} $}-algebra, which has the
  metric approximation property.
\newblock {\em Invent. Math.}, 50:279--293, 1978/79.

\bibitem[Hi]{hiai}
Hiai, F.~
\newblock {$q$}-deformed {A}raki-{W}oods algebras.
\newblock In {\em Operator algebras and mathematical physics (Constan\c ta,
  2001)}, 169--202. Theta, Bucharest, 2003.

\bibitem[Ho]{holm}
Holmstedt, T.~
\newblock Interpolation of quasi-normed spaces.
\newblock {\em Math. Scand.}, 26:177--199, 1970.

\bibitem[Ja1]{ja1}
Jajte, R.~
\newblock {\em Strong limit theorems in noncommutative probability}.
\newblock {\em Lect. Notes in Math.}, 1110, Springer-Verlag, Berlin, 1985.

\bibitem[Ja2]{ja2}
Jajte, R.~
\newblock {\em Strong limit theorems in noncommutative {$L\sb
2$}-spaces}.
 \newblock {\em Lect. Notes in Math.}, 1477, Springer-Verlag, Berlin, 1991.

\bibitem[Ju]{ju-doob}
Junge, M.~
\newblock Doob's inequality for non-commutative martingales.
\newblock {\em J. Reine Angew. Math.}, 549:149--190, 2002.

\bibitem[JX1]{jx-maxnote}
 Junge, M.~ and Xu, Q.~
 \newblock Th\'eor\`emes ergodiques maximaux dans les espaces {$L\sb p$}
              non commutatifs.
\newblock {\em C. R. Math. Acad. Sci. Paris}, 334:773--778, 2002.

\bibitem[JX2]{jx-burk}
Junge, M.~ and Xu, Q.~
\newblock Noncommutative {B}urkholder/{R}osenthal inequalities.
\newblock {\em Ann. Probab.}, 31:948--995, 2003.

\bibitem[JX3]{jx-const}
Junge, M.~ and Xu, Q.~
\newblock On the best constants in certain noncommutative martingale
  inequalities.
\newblock {\em Bull. London Math. Soc.}, 37:243--253, 2005.

\bibitem[JX4]{jx-ros}
Junge, M.~ and Xu, Q.~
\newblock Noncommutative Burkholder/Rosenthal inequalities: Applications.
\newblock To appear.

\bibitem[JX5]{jx-red}
Junge, M.~ and Xu, Q.~
\newblock Haagerup's reduction on noncommutative $L_p$-spaces and aplications.
\newblock In preparation.

\bibitem[Ka]{kad-cs}
Kadison, R.V.~
\newblock A generalized {S}chwarz inequality and algebraic invariants for
  operator algebras.
\newblock {\em Ann. of Math.}, 56:494--503, 1952.

\bibitem[Ko]{kos-int}
Kosaki, H.~
\newblock Applications of the complex interpolation method to a von {N}eumann
  algebra: noncommutative {$L\sp{p}$}-spaces.
\newblock {\em J. Funct. Anal.}, 56:29--78, 1984.

\bibitem[K{\"u}]{kum}
K{\"u}mmerer, B.~
\newblock A non-commutative individual ergodic theorem.
\newblock {\em Invent. Math.}, 46:139--145, 1978.

\bibitem[L]{la-erg}
Lance, E.C.~
\newblock Ergodic theorems for convex sets and operator algebras.
\newblock {\em Invent. Math.}, 37(3):201--214, 1976.

\bibitem[M]{mei}
Mei, T.~
\newblock Operator-valued Hardy spaces.
\newblock {\em Memoirs Amer. Math. Soc.}, to appear.

\bibitem[Mu]{Mu}
Musat, M.~
\newblock Interpolation between non-commutative {BMO} and
              non-commutative $L_p$-spaces
\newblock {\em J. Funct. Anal.}, 202:195--225, 2003.

\bibitem[N]{nou-inj}
Nou, A.~
\newblock Non-injectivity of the $q$-deformed.
\newblock {\em Math. Ann.}, 330:17--38, 2004.

\bibitem[Par]{parth-bk}
Parthasarathy, K.R.~
\newblock {\em An introduction to quantum stochastic calculus},
\newblock Birkh\"auser Verlag, Basel, 1992.

\bibitem[Pau]{pa-cb}
Paulsen, V.~
\newblock {\em Completely bounded maps and operator algebras},
\newblock Cambridge Univ. Press, 2002.

\bibitem[P]{pis-ast}
Pisier, G.~
\newblock Non-commutative vector valued {$L\sb p$}-spaces and completely
  {$p$}-summing maps.
\newblock {\em Ast\'erisque}, 247, 1998.

\bibitem[PX1]{px-BG}
Pisier, G.~ and Xu, Q.~
\newblock Non-commutative martingale inequalities.
\newblock {\em Comm. Math. Phys.}, 189:667--698, 1997.

\bibitem[PX2]{px-survey}
Pisier, G.~ and Xu, Q.~
\newblock Non-commutative {$L\sp p$}-spaces.
\newblock In {\em Handbook of the geometry of Banach spaces, Vol.\ 2}, pages
  1459--1517. North-Holland, Amsterdam, 2003.

\bibitem[Ra]{ran-mtrans}
Randrianantoanina, N.~
\newblock Non-commutative martingale transforms.
\newblock {\em J. Funct. Anal.}, 194:181--212, 2002.

\bibitem[Ri]{ricard-fac}
Ricard, E.~
\newblock Factoriality of q-Gaussian von Neumann algebras.
\newblock {\em Comm. Math. Phys.}, 257:659--665, 2005.

\bibitem[RX]{rixu}
Ricard, E.~ and Xu, Q.~
\newblock Khintchine type inequalities for free product and applications.
\newblock {\em J. Reine Angew. Math.}, to appear.

\bibitem[Sh]{shlya-quasifree}
Shlyakhtenko, D.~
\newblock Free quasi-free states.
\newblock {\em Pacific J. Math.}, 177:329--368, 1997.

\bibitem[Ska]{Ska}
Skalsi, A.~
\newblock On a classical scheme in noncommutative multiparameter
ergodic theory.
\newblock perprint.

\bibitem[Sta]{sta}
Starr, N.~
\newblock Operator limit theorems.
\newblock {\em Trans. Amer. Math.Soc.}, 121:90--111, 1966.

\bibitem[St1]{st-erg}
Stein, E.M.~
\newblock On the maximal ergodic theorem.
\newblock {\em Proc. Nat. Acad. Sci. U.S.A.}, 47:1894--1897, 1961.

\bibitem[St2]{st-lp}
Stein, E.M.~
\newblock {\em Topics in harmonic analysis related to the {L}ittlewood-{P}aley
  theory.}
\newblock Annals of Mathematics Studies, No. 63. Princeton Univ. Press,
  Princeton, N.J., 1985.

\bibitem[Te]{terp}
Terp, M.~
\newblock $L_p$ spaces associated with von Neumann algebras.
\newblock Notes, Math. Institute, Copenhagen Univ., 1981.

\bibitem[V]{voi-sym}
Voiculescu, D.~
\newblock Symmetries of some reduced free product {$C\sp \ast$}-algebras.
\newblock In {\em Operator algebras and their connections with topology and
  ergodic theory (Bu\c steni, 1983)}.  {\em Lect. Notes in
  Math.} 1132:556--588, 1985.

\bibitem[VDN]{VDN}
Voiculescu, D., Dykema, K. and Nica, A.~
\newblock {\em Free random variables}, vol.1 of {\em CRM Monograph Series}.
\newblock Amer. Math. Soc., Providence, RI, 1992.

\bibitem[X]{xu-martsurv}
Xu, Q.~
\newblock Recent development on noncommutative martingale inequalities.
\newblock {\em Proceedinds of International Conference on function
 space and its applications}, Wuhan 2003; pp. 283--313.
 Research Information Ltd
 UK, 2004.


\bibitem[Ye]{ye1}
Yeadon, F.J.~
\newblock Ergodic theorems for semifinite von {N}eumann algebras. {I}.
\newblock {\em J. London Math. Soc. (2)}, 16(2):326--332, 1977.

\bibitem[Yo]{yosida}
Yosida, K.~
\newblock {\em Functional analysis}.
\newblock Second edition.  Springer-Verlag New York, 1968.

\bibitem[Z]{zy}
Zygmund, A.~
\newblock {\em Trigonometric series. {V}ol. {I}, {II}}.
\newblock Cambridge Univ. Press, third edition, 2002.

\end{thebibliography}

\bigskip

\begin{itemize}
\item[M. J.:]  Department of Mathematics, University of Illinois,
Urbana, IL 61801 - USA\\
junge@math.uiuc.edu
\item[Q. X.:] Laboratoire de Math\'{e}matiques, Universit\'{e} de
Franche-Comt\'{e}, 25030 Besan\c con, cedex - France\\
qx@math.univ-fcomte.fr \end{itemize}

\end{document}